\documentclass[a4paper, reqno]{amsart}
\usepackage{url} 
\usepackage{rotating} 

\usepackage[utf8]{inputenc}
\usepackage[english]{babel}
\usepackage{amsmath,mathrsfs,amsfonts,amssymb,amsxtra,latexsym,amscd,amsthm,marvosym,multirow,mathtools}
\usepackage[a4paper,inner=2.3cm, outer=2.3cm, top=2.3cm, bottom=2.3cm]{geometry}
\usepackage{setspace}
\usepackage{subcaption} 
\usepackage{eucal}
\usepackage{enumerate}
\usepackage{csquotes}
\usepackage[shortlabels]{enumitem}
\usepackage{cases}
\usepackage{cleveref}
\usepackage{soul}

\usepackage{tikz}
\usetikzlibrary{patterns}
\newtheorem{theorem}{Theorem}[section]
\newtheorem{lemma}[theorem]{Lemma}

\newtheorem{remark}[theorem]{Remark}
\newtheorem{proposition}[theorem]{Proposition}
\newtheorem{corollary}[theorem]{Corollary}

\renewcommand{\Re}{\mathrm{Re}\,}
\renewcommand{\Im}{\mathrm{Im}\,}

\DeclareMathOperator*{\Ker}{Ker \;}

\DeclareMathOperator{\sgn}{sgn}

\def\beq{\begin{eqnarray*}}\def\eeq{\end{eqnarray*}}
\def\bq{\begin{equation}}\def\eq{\end{equation}}

\newcommand{\N}{\mathbb{N}}

\newcommand{\R}{\mathbb{R}}
\newcommand{\C}{\mathbb{C}}

\newcommand{\dd}{\mathrm{d}}
\newcommand{\D}{\mathcal{D}}

\newcommand{\coloneqq}{=}

\makeatletter

\@addtoreset{equation}{section}  
\makeatother

\begin{document}
\title[Dislocated Dirac operators]{Spectral analysis of Dirac operators \\ for dislocated potentials \\ with a purely imaginary jump}
\author{Lyonell Boulton}
\address[Lyonell Boulton]{Department of Mathematics and Maxwell Institute for Mathematical Sciences, Heriot-Watt University, Edinburgh, EH14 4AS, United Kingdom.}
\email{L.Boulton@hw.ac.uk}
\author{David Krej\v{c}i\v{r}\'{i}k}
\address[David Krej\v{c}i\v{r}\'{i}k]{Department of Mathematics, Faculty of Nuclear Sciences and Physical Engineering,
Czech Technical University in Prague, Trojanova 13, 12000 Prague 2, Czech Republic}
\email{david.krejcirik@fjfi.cvut.cz}
\author{Tho Nguyen Duc}
\address[Tho Nguyen Duc]{Faculty of Mathematics and Statistics, Ton Duc Thang University, Ho Chi Minh City, Vietnam.}
\email{nguyenductho@tdtu.edu.vn}

\begin{abstract}
In this paper we present a complete spectral analysis of Dirac operators with non-Hermitian matrix potentials of the form $i\operatorname{sgn}(x)+V(x)$ where $V\in L^1$. For $V=0$ we compute explicitly the matrix Green function. This allows us to determine the spectrum, which is purely essential, and its different types. It also allows us to find sharp  enclosures for the pseudospectrum and its complement, in all parts of the complex plane. Notably, this includes the instability region, corresponding to the interior of the band that forms the numerical range. Then, with the help of a Birman-Schwinger principle, we establish in precise manner how the spectrum and pseudospectrum change when $V\not=0$, assuming the hypotheses $\|V\|_{L^1}<1$ or $V\in L^1\cap L^p$ where $p>1$. We show that the essential spectra remain unchanged and that the $\varepsilon$-pseudospectrum stays close to the instability region for small $\varepsilon$. We determine sharp asymptotics for the discrete spectrum, whenever $V$ satisfies further conditions of decay at infinity. Finally, in one of our main findings, we give a complete description of the weakly-coupled model. 
 \end{abstract}

\date{8th April 2025}
\maketitle
\tableofcontents
\newpage


\section{Introduction}
\subsection{Context and motivations}
The significance of the Dirac equation lies on the fact that it describes, 
not only relativistic particles in quantum mechanics,
but also advanced materials such as graphene.
Mathematically, the study of this equation is notoriously difficult, partly 
because of the spinorial structure
of the underlying Hilbert space and partly because of the lack of  
semi-boundedness and positivity-preserving properties. 
In this respect, the spectral analysis of the underlying matrix differential operator is considerably more challenging 
than that of its non-relativistic counterpart, the scalar Schr{\"o}dinger operator.

Recently, we have seen unprecedented interest in the study of
non-selfadjoint electromagnetic perturbations of the Dirac equation  
\cite{Cuenin-Laptev-Tretter14,Cuenin_2014,Cuenin_2017,Cuenin-Siegl18,
CIKS,Krejcirik-Nguyen22,HK2,
AFKS,Ancona-Fanelli-Schiavone_2022,Mizutani-Schiavone_2022,
Heriban-Tusek_2022,
KK10}.
This is justified, on the one hand, by the availability of 
new proper connections with quantum mechanics 
\cite{GHS,Krejcirik-Siegl-Tater-Viola15}, and, on the other hand, by
the new conceptual paradigms posed by the notions of pseudospectra of linear operators.

This paper is motivated by the current lack of understanding 
of the structure of the pseudospectra of non-selfadjoint Dirac operators. 
Indeed, the only available work in this direction seems to be 
the recent non-semiclassical construction of pseudomodes 
given in~\cite{Krejcirik-Nguyen22}.
\medskip

In order to provide an insight into the pseudospectral properties
of the Dirac equation in general, and at the same time set a benchmark that differs from the available scalar models given by the classical non-selfadjoint Schr{\"o}dinger operator, the present paper is devoted to introducing and examining in close detail the non-selfadjoint Dirac operator, 
\begin{equation}\label{our}
  \mathscr{L}_{m,V}
  \coloneqq 
  \begin{pmatrix}
  m & -\partial_x \\
  \partial_x & -m  
  \end{pmatrix}
  +   \begin{pmatrix}
  i \sgn(x) & 0 \\
  0 & i\sgn(x)
  \end{pmatrix}
  + V(x),
\end{equation}
on a suitable domain of $L^2(\R,\C^2)$ 
where $V:\R\longrightarrow\C^{2 \times 2}$ is a long-range possibly non-Hermitian matrix potential.
The unperturbed operator $\mathscr{L}_{m}=\mathscr{L}_{m,0}$ 
is a relativistic non-selfadjoint version of the quantum mechanical infinite square well. 

The choice is motivated 
by the scalar Schr{\"o}dinger case   on $L^2(\R,\C)$ considered in~\cite{Henry-Krejcirik17},
\begin{equation}\label{Henry}
  -\partial_x^2 + i\sgn(x)
\end{equation}
  with potential perturbations, and it serves as a link with the analysis conducted in~\cite{Cuenin-Laptev-Tretter14}.
The scalar model \eqref{Henry} was instrumental
for the non-semiclassical construction of pseudomodes 
established in~\cite{Krejcirik-Siegl19}, which eventually solved a notorious open problem raised during 
a workshop at the American Institute of Mathematics in 2015, \emph{cf.} \cite[Open Problem~10.1]{AIM-2015}.
The relativistic variant of this construction has now been reported in~\cite{Krejcirik-Nguyen22}.

The  simplicity of the linear operator~\eqref{Henry} is deceiving, as it hides  
a non-trivial structure on its  pseudospectrum. 
Our findings reveal two main distinctions between \eqref{our} and this scalar model. One is about this structure 
and one is about the available tools to analyse it. The latter is to be expected, given the higher degree of complexity of a matrix versus a scalar operator.  But the former is rather surprising. As we shall see in the next section, the linear operator \eqref{our}  
has a resolvent norm that grows quadratically, rather than linearly, at infinity inside the band forming the spectral instability region. 
 
\medskip
\subsection{Structure of the paper}
The organisation of the paper is as follows. 
Section~\ref{Sec Main Results} 
is devoted to a proper mathematical description of our main contributions. 
The body of the paper is 
Sections~\ref{Sec Res and Sp}--\ref{Sec L1 Perturbation}, where we establish the proofs of these results. After this, Appendix~\ref{Appendix Abstract Birman Schwinger} includes complete details of the Birman-Schwinger-type principle formulated 
in \cite{Gesztesy-Latushkin-Mitrea-Zinchenko20} which is crucial to our analysis.

The core Sections~\ref{Sec Res and Sp}--\ref{Sec L1 Perturbation} 
consist of two parts.
\begin{itemize}
\item \textbf{The unperturbed operator,} comprising 
Sections~\ref{Sec Res and Sp} and~\ref{Sec Resolvent Estimate}.
We begin our study of~\eqref{our} 
in Section~\ref{Sec Res and Sp} 
by finding the explicit expression of the matrix Green function of
the unperturbed operator~$\mathscr{L}_{m}$. 
This allows to determine 
the spectrum of~$\mathscr{L}_{m}$.  
As it turns out, the latter is purely essential and it comprises four symmetric segments on the boundary of the band $\Sigma \coloneqq \{|\Im z|< 1\}$ for $m>0$, while by contrast it is equal to $\overline\Sigma$ for $m=0$. In Section~\ref{Sec Resolvent Estimate} we find enclosures for the pseudospectrum and its complement, in all parts of the complex plane. In one of our main contributions, we compute the explicit asymptotic constants, 
up to order~0, of the resolvent norm of~$\mathscr{L}_{m}$  
inside the instability region~$\Sigma$. 
The latter coincides with the interior of the numerical range. 
\item \textbf{Perturbations}, comprising Section~\ref{Sec L1 Perturbation}. 
The second part of the paper is devoted to the case $V\not=0$. By applying a non-selfadjoint version of the classical Birman-Schwinger principle, we establish in precise manner how spectrum and pseudospectrum change, under two general hypotheses, 
$\|V\|_{L^1}<1$ or $V\in L^1(\R,\C^{2\times 2})\cap L^p(\R,\C^{2\times 2})$ for some $p>1$. 
We show that the essential spectra remain unchanged in both cases and that the $\varepsilon$-pseudospectrum stays close to the instability region 
for small $\varepsilon > 0$ whenever $\|V\|_{L^1}<1$. 
Then, we formulate our two other main contributions of this part. On the one hand, we determine sharp asymptotics for the discrete spectrum, whenever $V$ satisfies further conditions of decay at infinity. On the other hand, we give a complete description of the weakly-coupled model, corresponding to potential $\epsilon V$ in the regime $\epsilon\to0$. 
\end{itemize} 

\subsection{Notation used throughout the work}
The following specific conventions will be used throughout this paper. 
\begin{itemize}
\item $\R_{+}=[0,+\infty)$ and $\R_{-}=(-\infty,0]$.
\item $[[m,n]]\coloneqq  \{k\in\mathbb{N}\,:\, m\leq k \leq n\}$ for $m,n\in \R$.
\item For a closed operator $\mathscr{L}:\operatorname{Dom}(\mathscr{L})\longrightarrow \mathcal{H}$, the \emph{spectrum} is denoted by
\[ \operatorname{Spec}(\mathscr{L})=\operatorname{Spec}_{\mathrm{p}}(\mathscr{L}) \cup \operatorname{Spec}_{\mathrm{r}}(\mathscr{L}) \cup \operatorname{Spec}_{\mathrm{c}}(\mathscr{L}),\]
where the \emph{point}, \emph{residual} and \emph{continuous} spectrum denote, as usual,
\begin{align*}
&\operatorname{Spec}_{\mathrm{p}}(\mathscr{L})\coloneqq   \{z\in \C: \mathscr{L}-z \text{ is not injective} \},\\
&\operatorname{Spec}_{\mathrm{r}}(\mathscr{L})\coloneqq   \{z\in \C: \mathscr{L}-z \text{ is injective and } \overline{\mathrm{Ran}(\mathscr{L}-z)} \subsetneq \mathcal{H} \},\\
&\operatorname{Spec}_{\mathrm{c}}(\mathscr{L})\coloneqq   \{z\in \C: \mathscr{L}-z \text{ is injective and } \overline{\mathrm{Ran}(\mathscr{L}-z)}=\mathcal{H} \text{ and } \mathrm{Ran}(\mathscr{L}-z)\subsetneq \mathcal{H} \}.
\end{align*}
The \emph{resolvent set}  is denoted by $\rho(\mathscr{L})$.
\item We use the classical definitions of the \emph{five types of essential spectrum}, \cite[Sec. IX]{Edmunds-Evans18} or \cite[Sec. 5.4]{Krejcirik-Siegl15}. We write $\operatorname{Spec}_{\mathrm{ej}}(\mathscr{L})$ for $\mathrm{j}\in [[1,5]]$, where $\mathscr{L}-z$
\begin{itemize}
\item[-] is not semi-Fredholm for $\operatorname{e1}$,
\item[-] possess a singular Weyl sequence for $\operatorname{e2}$,
\item[-] is not Fredholm for $\operatorname{e3}$,
\item[-] is not Fredholm of index 0 for $\operatorname{e4}$,
\item[-] is such that either $z\in \operatorname{Spec}_{\mathrm{e1}}(\mathscr{L})$ or the resolvent set does not intersect the connected component of the complement of $\operatorname{Spec}_{\mathrm{e1}}(\mathscr{L})$ where $z$ lies, for $\operatorname{e5}$. 
\end{itemize}
\item The \emph{discrete spectrum} is denoted as 
\[ \operatorname{Spec}_{\mathrm{dis}}(\mathscr{L})=\C\setminus \operatorname{Spec}_{\mathrm{e5}}(\mathscr{L}),\]
and it is the set of isolated eigenvalues such that the Riesz projector has a finite-dimensional range. 
\item For $\varepsilon > 0$, the $\varepsilon$-\emph{pseudospectrum}  is denoted by
\begin{equation}\label{Def Pseudospectrum}
\begin{aligned} \operatorname{Spec}_{\varepsilon}(\mathscr{L})&= \operatorname{Spec}(\mathscr{L})\cup \{z\in \rho(\mathscr{L}): \Vert (\mathscr{L}-z)^{-1}\Vert> \varepsilon^{-1} \}\\&=\bigcup_{\Vert V\Vert<\varepsilon} \operatorname{Spec}(\mathscr{L}+V).\end{aligned}
\end{equation}
By convention, we set $\operatorname{Spec}_0(\mathscr{L})=\operatorname{Spec}(\mathscr{L})$.
\item  $\vert \cdot \vert_{\C^{2}}$ is the Euclidean norm for vectors. 
\item $\vert \cdot \vert_{\C^{2\times 2}}$ is the operator norm induced by $\vert \cdot \vert_{\C^{2}}$ for matrices. 
\item $\vert \cdot \vert_{\mathrm{F}}$ is the Frobenius norm for matrices.
\item For a measurable matrix-valued function $V:\R \to \C^{2\times 2}$, we write $V\in L^p(\R,\C^{2\times 2})$ for $p\in [1,\infty]$ to indicate that $|V|_{\C^{2\times 2}} \in L^p(\R,\C)$. We denote the $L^p$ norm of $V$ by
\[ \Vert V \Vert_{L^p} = \Vert |V|_{\C^{2\times 2}} \Vert_{L^p} .\]
\item For a positive function $\nu'(x)$, we define the measure $\dd \nu(x) = \nu'(x) \dd x$. We say that $V\in L^1(\R,\C^{2\times 2};\dd \nu)$ if
\[ \Vert V \Vert_{L^1(\dd \nu)} = \int_{\R} \vert V(x) \vert_{\C^{2\times 2}} \nu'(x)\, \dd x<\infty.\]
\end{itemize}


\section{Mathematical framework and summary of results} \label{Sec Main Results}
Below, we set
\begin{equation}\label{Set S_m}
S_{m}\coloneqq  \left\{
\begin{aligned}
&\{z\in \C: |\Re z| \geq m, |\Im z|=1 \},\qquad &&\text{if } m>0,\\
&\{z\in \C : |\Im z| \leq 1 \},\qquad &&\text{if } m=0,
\end{aligned}
\right.
\end{equation}
and we write $L^{2}(\R,\C^2)=L^2(\R,\C)\oplus L^2(\R,\C)$ with the inner product
\[ \left\langle f,g \right\rangle \coloneqq \int_{\R} f_{1}(x) \overline{g_{1}(x)}+f_{2}(x) \overline{g_{2}(x)}\, \dd x \qquad \text{ for  components } f_j,\,g_j \in L^2(\R,\C).\]
As usual, 
\[\sigma_{1}\coloneqq  \begin{pmatrix}
0 && 1\\
1 && 0
\end{pmatrix},\qquad\sigma_{2}\coloneqq  \begin{pmatrix}
0 && -i\\
i && 0
\end{pmatrix} \quad \text{and} \quad\quad\sigma_{3}\coloneqq  \begin{pmatrix}
1 && 0\\
0 && -1
\end{pmatrix}, \]
denote the Pauli matrices.
For mass $m \geq 0$, let  $\mathcal{D}_{m}:H^{1}(\R,\C^2)\longrightarrow L^{2}(\R,\C^2)$ be the one-dimensional free particle Dirac operator, given by
\begin{equation}\label{Free Dirac}
\begin{aligned}
&(\mathcal{D}_{m} f)(x) \coloneqq  \left(-i\partial_{x} \right)\sigma_2f(x) + m\sigma_{3}f(x).
\end{aligned}
\end{equation}
It is well-known that $\D_m$ is a selfadjoint linear operator and that its spectrum is 
\[\mathrm{Spec}(\D_{m})= (-\infty,-m] \cup [m, +\infty). \] Moreover, $\D_{m}$ is unitarily equivalent to the selfadjoint operator $\widetilde{\D}_{m}=\left(-i\partial_{x} \right)\sigma_{1} +m \sigma_{3}:H^{1}(\R,\C^2)\longrightarrow L^{2}(\R,\C^2)$ via transformation by $\begin{pmatrix}
i & 0\\
0 & -1
\end{pmatrix}.$
This other realisation of the free particle Dirac operator is the one considered in \cite{Cuenin-Laptev-Tretter14,Cuenin-Siegl18,Krejcirik-Nguyen22}. In the present paper we prefer the formulation \eqref{Free Dirac}, noting that the results we report below, map in a straightforward manner onto $\widetilde{\D}_{m}$.  

\subsection{The Dirac operator with a dislocation}
In the present paper we examine the linear operator $\mathscr{L}_{m}:H^{1}(\R,\C^2)\longrightarrow L^{2}(\R,\C^2)$ given by
\begin{equation}\label{Dirac Operator}
\begin{aligned}
&(\mathscr{L}_{m}f)(x) \coloneqq  \mathcal{D}_{m}f(x) +i\sgn(x)  f(x) 
\end{aligned}
\end{equation}
and its perturbations by a potential.
Since multiplication by $i\sgn(x)  I $ is a bounded operator, $\mathscr{L}_{m}$ is also closed. From the explicit expression of the adjoint, \[ (\mathscr{L}_{m}^*f)(x) \coloneqq  \mathcal{D}_{m}f(x) -i\sgn(x)  f(x),\]
it follows that $\mathscr{L}_{m}$ is not a normal operator. As we shall see below, the spectral and pseudospectral properties of $\mathscr{L}_{m}$ and its perturbations, are rather unusual and very interesting. 

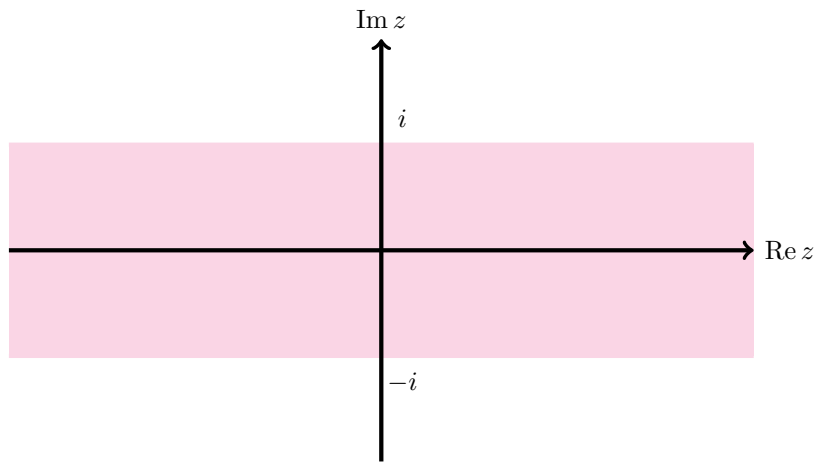
\begin{figure}[h]
\centering
\begin{tikzpicture}[scale=0.7]
\fill[ fill=magenta!20](-7,-2) rectangle (7,2);

\draw[-,ultra thick, ultra thick, magenta!20] (-7,2)--(7,2);
\draw[-,ultra thick, ultra thick, magenta!20] (-7,-2)--(7,-2);
\draw[->,ultra thick] (-7,0)--(7,0) node[right]{$\Re z$};
\draw[->,ultra thick] (0,-4)--(0,4) node[above]{$\Im z$};

\draw (0.4,2.5) node {$i$};
\draw (0.4,-2.5) node {$-i$};
\end{tikzpicture}
\caption{The numerical range of $\mathscr{L}_{m}$.}\label{Fig:Numerical Range}
\end{figure}

The numerical range of $\mathscr{L}_{m}$ is the infinite closed strip
\[\textup{Num}(\mathscr{L}_{m})= \R+i[-1,1].\]
 See Figure \ref{Fig:Numerical Range} and the beginning of Section~\ref{Sec Res and Sp}.
Unlike the Schr{\"o}dinger operator with dislocated potential analysed in \cite{Henry-Krejcirik17}, $\mathscr{L}_{m}$ is neither sectorial nor $\mathcal{PT}$-symmetric. However,  $\mathscr{L}_{m}$ is $\mathcal{T}$-selfadjoint and this implies that the spectrum respects some of the symmetries that the Schr{\"o}dinger model possesses.

 Here and everywhere below a distinction of the case $m=0$, is evident, unavoidable and it is typical of families of non-selfadjoint operators dependent on a parameter, \emph{cf.} the examples in \cite{B24} and references therein.  
Moreover, setting an unitary operator $Sf=g$ such that $g(x)=a f(a^2 x)$ for fixed non-zero $a\in \mathbb{R}\cup i\mathbb{R}$, it is readily seen that $S$ preserves the domains and
\[
(S^{*}\mathscr{L}_{a^2m}Sf)(x)=a^2\left(\mathcal{D}_{m}f(x)+\frac{i}{a^2}\sgn\left( \frac{x}{a^2}\right) f(x)\right).
\]
This shows that both the case $m<0$ and the case of a dislocated potential of the form $c\sgn(x)I$, where $c>0$ is a constant, are covered by our results below. 

\begin{figure}[h]
\centering
\begin{tikzpicture}[scale=0.7]

\draw[-,ultra thick, red] (2,2)--(7,2);
\draw[-,ultra thick,  red] (2,-2)--(7,-2);
\draw[-,ultra thick,  red] (-2,2)--(-7,2);
\draw[-,ultra thick,  red] (-2,-2)--(-7,-2);
\draw[-,ultra thick,  thick] (2,-0.3) -- (2,0.3);
\draw[-,ultra thick,  thick] (-2,-0.3) -- (-2,0.3);
\draw[->,ultra thick] (-7,0)--(7,0);
\draw[->,ultra thick] (0,-4)--(0,4);

\draw (-0.4,2.1) node {$i$};
\draw (-0.5,-2.1) node {$-i$};
\draw (2, 0) node[below=4pt] {$m$};
\draw (-2, 0) node[below=4pt] {$-m$};
\end{tikzpicture}
\caption{The red lines represent the spectrum of the operator $\mathscr{L}_{m}$ when $m>0$.}\label{Fig:Spectrum m>0}
\end{figure}
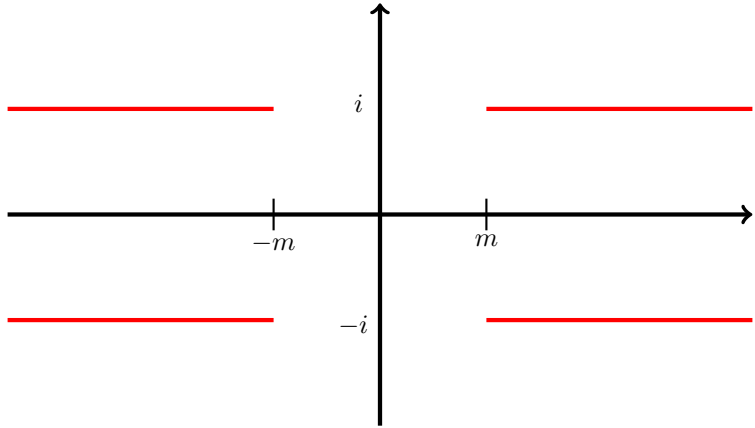

As we shall see in Section~\ref{Sec Res and Sp}, the spectrum of $\mathscr{L}_{m}$ can be determined analytically using arguments based on resolvent construction and Weyl sequences. For all $m>0$, the spectrum is purely continuous
\[ 
\operatorname{Spec}(\mathscr{L}_{m})=\operatorname{Spec}_{\textup{c}}(\mathscr{L}_{m})=S_{m}.
\]
Furthermore, all essential spectra coincide,
\[ \operatorname{Spec}_{\mathrm{ej}}(\mathscr{L}_{m})=\operatorname{Spec}(\mathscr{L}_{m}), \qquad \forall \mathrm{j}\in [[1,5]].\]
Note that $\mathcal{T}$-selfadjointness implies equality of the essential spectra for $\mathrm{j}\in [[1,4]]$ only, so the proof of equality of the full set requires extra arguments given below. An illustration of the spectrum is included in  Figure~\ref{Fig:Spectrum m>0}. 

\begin{figure}[h]
\centering
\begin{tikzpicture}[scale=0.7]
\fill[ fill=red!20](-7,-2) rectangle (7,2);

\draw[-,ultra thick, red] (-7,2)--(7,2);
\draw[-,ultra thick,  red] (-7,-2)--(7,-2);

\draw[->,ultra thick] (-7,0)--(7,0);
\draw[->,ultra thick] (0,-4)--(0,4);

\draw (-0.4,2.3) node {$i$};
\draw (-0.5,-2.3) node {$-i$};

\end{tikzpicture}
\caption{The red lines show the continuous spectrum of $\mathscr{L}_{0}$. These lines coincide with the essential spectra $\mathrm{e}1,\ldots,\mathrm{e}4$. The inner part in light colour shows the point spectrum. The union of both these regions form the essential spectrum $\mathrm{e}5$ and also the full $\operatorname{Spec}(\mathscr{L}_{0})$.}\label{Fig:Spectrum m=0}
\end{figure}
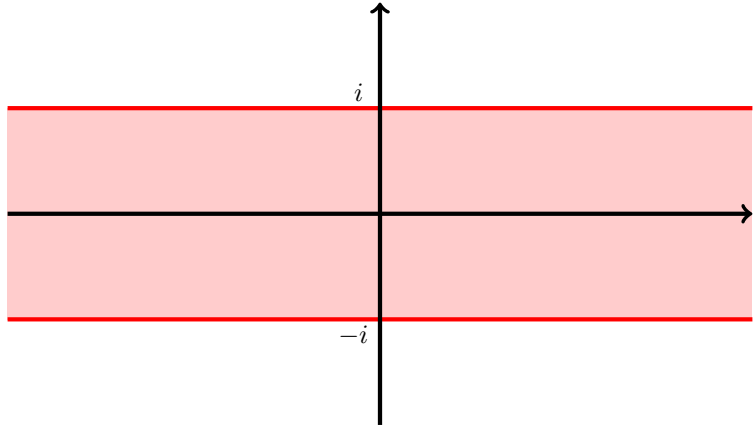

In contrast, for $m=0$, the spectrum is characterized by
\begin{gather*}
\operatorname{Spec}_{\mathrm{p}}(\mathscr{L}_{0})=\{z\in \C: |\Im z|<1\}
\qquad \text{and} \qquad
\operatorname{Spec}_{\mathrm{c}}(\mathscr{L}_{0})=\{z\in \C: |\Im z|=1 \}.
\end{gather*}
Unlike the case $m>0$, the essential spectra are not identical: 
\begin{gather*}
	\operatorname{Spec}_{\mathrm{ej}}(\mathscr{L}_{0})=\{z\in \C:|\Im z|=1\}, \qquad \forall \mathrm{j} \in [[1,4]] \qquad\text{and} \qquad \operatorname{Spec}_{\mathrm{e5}}(\mathscr{L}_{0})=\operatorname{Spec}(\mathscr{L}_{0}).
\end{gather*}
See the illustration in Figure~\ref{Fig:Spectrum m=0}.

For all $m\geq0$, the comparison with the spectrum of $\mathcal{D}_{m}$ is obvious and striking. Despite the seemingly simple structure of the dislocation $i\sgn(x)  I $, the perturbation splits the spectrum into two halves. In the case $m>0$, this phenomenon mimics the findings of  \cite[Thm. 2.1]{Henry-Krejcirik17} for the Schr{\"o}dinger case. For $m=0$, rather than a clear-cut split, the spectrum smears throughout the full band that forms the numerical range. In Lemma~\ref{Lem Point Spectrum} we show that, at any $z$ in the point spectrum, $(\mathscr{L}_{0}-z)$ is Fredholm with index zero. Moreover, the geometric multiplicity of $z$ is one while its algebraic multiplicity is infinite. This is neither obvious, nor a consequence of general principles.

\medskip

Let us now describe the pseudospectrum of $\mathscr{L}_{m}$.
We adopt the following terminology, which is convenient in order to articulate  our findings. We say that a closed operator $\mathscr{L}$ has a \emph{trivial pseudospectrum}, if there exists a constant $C\geq 1$ such that
\begin{equation}\label{Def Trivial Pseudospectrum}
\operatorname{Spec}_{\varepsilon}(\mathscr{L}) \subseteq \{ z\in \C: \mathrm{dist}(z,\operatorname{Spec}\mathscr{L})\leq C \varepsilon\}
\end{equation}
for all $\varepsilon>0$.
Any normal operator has a trivial pseudospectrum with $C=1$, since in this case the $\varepsilon$-pseudospectrum coincides with the $\varepsilon$-neighbourhood of the spectrum. Any closed operator which is similar to a normal operator, via a bounded invertible similarity transformation $S$, has also trivial pseudospectrum with $C= \|S\|\|S^{-1}\|$, the condition number of $S$. 

According to the well-known inequalities
\begin{equation*}
   \frac{1}{\operatorname{dist}\left(z,\operatorname{Spec}\mathscr{L}_{m}\right)}\leq \|(\mathscr{L}_{m}-z)^{-1}\|\leq \frac{1}{\operatorname{dist}\left(z,\operatorname{Num}\mathscr{L}_{m}\right)}
\end{equation*}
for all $z\not\in\operatorname{Num}(\mathscr{L}_{m})$, it is readily seen that for $|\Im z|>1$,
\begin{enumerate}[label=\textup{\textbf{\arabic*)}}]
\item and for $\vert \Re z \vert\geq m$,
\begin{equation} \label{Eq Im E>1}
\begin{aligned}
\Vert (\mathscr{L}_{m}-z)^{-1} \Vert =\frac{1}{\vert\Im z\vert-1};
\end{aligned}
\end{equation}
\item while for $\vert \Re z \vert< m$,
\begin{equation*}
\begin{aligned}
\frac{1}{\sqrt{\left( \vert \Re z \vert-m\right)^2+\left(\vert \Im z \vert-1\right)^2}}\leq\Vert (\mathscr{L}_{m}-z)^{-1} \Vert \leq\frac{1}{\vert\Im z\vert-1}.
\end{aligned}
\end{equation*}

\end{enumerate}
In the next theorem we see that, in stark contrast, the resolvent norm increases at a quadratic rate as $z\to\infty$ inside the numerical range.
This statement is our first main contribution.

\begin{theorem}\label{Theo Resolvent 2}
Let $m> 0$ and let $\mathscr{L}_{m}$ be the operator given by \eqref{Dirac Operator}. Assume that $z \in  \C$ is such that $ \vert \Im z \vert <1$. Then,
\begin{equation}\label{Accurate 1}
	\left\|(\mathscr{L}_{m}-z)^{-1}\right\| = \frac{|\Re z|^2}{m(1-|\Im z|^2)}\left(1+\mathcal{O}\big(|\Re z|^{-2}\big)\right)
\end{equation}
as $|\Re z| \to \infty$. The $\mathcal{O}$ term is uniform in $|\Im z|$ and locally uniform in $m$ \textup{(}\emph{i.e.}, the constants involved can be chosen independent of $z\in \C$ for all $|\Im z|<1$ and $m$ on any fixed compact subset of $(0,\infty)$\textup{)}.
\end{theorem}

\begin{figure}[h]
\includegraphics[width=.45\textwidth]{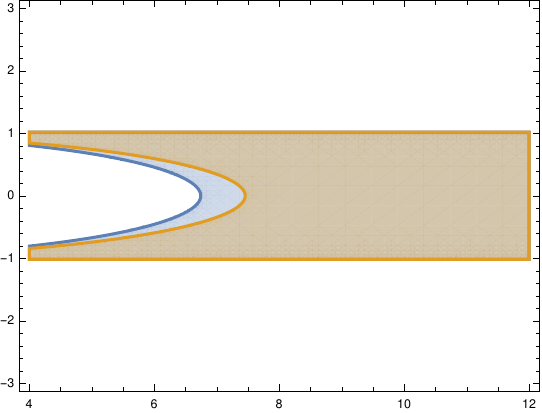}\qquad 
\includegraphics[width=.45\textwidth]{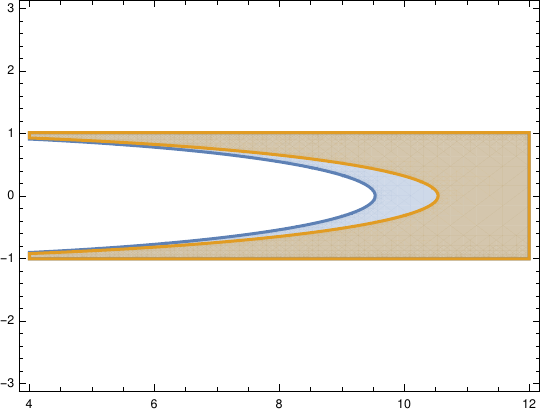}
\caption{Regions $\Lambda^{0.1}_{+}(\varepsilon)$ and $\Lambda^{0.1}_{-}(\varepsilon)$ in the box as shown, for  $m=1$, $\varepsilon=0.02$ (left) and $\varepsilon=0.01$ (right). The part in blue corresponds to $\Lambda^{0.1}_{+}(\varepsilon)\setminus \Lambda^{0.1}_{-}(\varepsilon)$ and $\Lambda^{0.1}_{-}(\varepsilon)$ is shows in brick.}\label{Level set}
\end{figure}

According to this theorem, when we consider the level set \[ \left\{z\in \rho(\mathscr{L}_{m}):\left\|(\mathscr{L}_{m}-z)^{-1}\right\|=\frac{1}{\varepsilon} \right\}\]  for some $\varepsilon>0$, we obtain curves whose points $z\in\mathbb{C}$ satisfy 
\[ |\Im z|^2= 1- \frac{\varepsilon |\Re z|^2}{m}\left(1+\mathcal{O}\left(\frac{1}{|\Re z |^2} \right)\right)\]
in the asymptotic regime $|\Re z| \to +\infty$. Motivated by this, for $\alpha\in(0,1)$ and $\varepsilon>0$, we set regions
\[
\Lambda_{\pm}^{\alpha}(\varepsilon)\coloneqq   \left\{z\in\C: |\Im z|\leq 1+\varepsilon \text{ and }
 |\Re z|^2\geq \frac{m(1-|\Im z|^2)}{(1\pm \alpha)\varepsilon} \right\}.
\]
These are illustrated in Figure~\ref{Level set} for two different values of $\varepsilon$.

In the next corollary we establish sharp inclusions and exclusion zones for the pseudospectra of $\mathscr{L}_m$ in terms of these regions. Note that $\mathscr{L}_{0}$ has trivial pseudospectra although it is a non-normal operator.

\begin{corollary}
Let $m\geq 0$ and let $\mathscr{L}_{m}$ be the operator given by \eqref{Dirac Operator}.  
\begin{enumerate}[label=\textup{\textbf{\arabic*)}}]
\item If $m=0$, then the pseudospectrum of $\mathscr{L}_{0}$ is trivial and for all $\varepsilon>0$,
\[ \operatorname{Spec}_{\varepsilon}(\mathscr{L}_{0})=\{z\in \C: |\Im z|\leq 1+\varepsilon \}.\]
\item If $m>0$, then the pseudospectrum of $\mathscr{L}_{m}$ is non-trivial. Namely, for every fixed $\alpha\in (0,1)$, there exists $M>0$ such that, 
\begin{align}
\Lambda_-^{\alpha}(\varepsilon)\!\cap\! \{|\Re{z}|>M\} &\ \ \subset \ \ \left(\operatorname{Spec}_{\varepsilon}\mathscr{L}_{m}\right)\!\cap\! 
\{|\Re{z}|>M\}\label{eq in} \\ \left(\operatorname{Spec}_{\varepsilon}\mathscr{L}_{m}\right)\!\cap\! \{|\Re{z}|>M\}&\ \ \subset\ \ \Lambda_+^{\alpha}(\varepsilon)\!\cap\! \{|\Re{z}|>M\} \label{eq out}
\end{align}
for all $\varepsilon>0$.
\end{enumerate}
\end{corollary}

\begin{proof}
For $m=0$, the statement follows directly from the spectrum of $\mathscr{L}_{0}$ and \eqref{Eq Im E>1}. Here is the proof for $m>0$. According to Theorem~\ref{Theo Resolvent 2}, there exists a constant $C>0$ such that 
\begin{align*}
   \left(1-\frac{C}{|\Re z|^2}\right)\frac{|\Re z|^2}{m(1-|\Im z|^2)}< \Vert (\mathscr{L}_{m}&-z)^{-1} \Vert<\left(1+\frac{C}{|\Re z|^2}\right)\frac{|\Re z|^2}{m(1-|\Im z|^2)}
\end{align*}
for large $|\Re(z)|$ and $|\Im z|<1$. Thus, for fixed $\alpha\in(0,1)$ and $|\Re z|\geq \sqrt{\frac{C}{\alpha}}$ we have the following.  For $\varepsilon>0$, \eqref{eq in} is ensured whenever $|\Re z|$ is such that
\[\left(1-\alpha\right)\frac{|\Re z|^2}{m(1-|\Im z|^2)}\geq \frac{1}{\varepsilon}.\]
Similarly, if
\[
   \frac{1}{\varepsilon} \leq \Vert (\mathscr{L}_{m}-z)^{-1} \Vert,
\]
then
\[ \frac{1}{\varepsilon}< \left(1+\alpha\right)\frac{|\Re z|^2}{m(1-|\Im z|^2)}.\] Therefore \eqref{eq out} is valid. For $|\Im z|>1$, we recall \eqref{Eq Im E>1}. \end{proof}

\subsection{Long-range potential perturbations} \label{subsection2.2}
In Section~\ref{Sec L1 Perturbation} we examine perturbations of $\mathscr{L}_{m}$ by an integrable matrix-valued function $V$. We establish results about the localisation of eigenvalues and pseudospectrum, under general assumptions on the decay of $V$ at infinity. 

Our starting point is the construction of a specific closed densely defined extension of the differential operator \[\mathscr{L}_m+V\,:\,H^{1}(\R,\C^2) \cap \mathrm{Dom }(V)\longrightarrow L^{2}(\R,\C^2),\] via
a version of the classical selfadjoint framework of Kato \cite{Kato66}, obtained in \cite{Gesztesy-Latushkin-Mitrea-Zinchenko20}. Concretely, for almost\footnote{Here and elsewhere below we assume without further mentioning that all the identities involving point-wise evaluation are valid almost everywhere with respect to the Borel $\sigma$-algebra on $\mathbb{R}$.} all $x\in \R$, let
\begin{equation*}
V(x)=B(x)A(x) \quad \text{where} \quad B(x)\coloneqq  U(x) \vert V(x) \vert^{1/2} \quad\text{and}\quad A(x)\coloneqq  \vert V(x) \vert^{1/2},
\end{equation*}
in the polar decomposition $V(x)=U(x)\vert V(x) \vert$,
for $U(x)$ partial isometries and $|V|(x) = (V^*(x)V(x))^{1/2}$.
Consider the family of closed operators
\[ Q(z)\coloneqq \overline{AR_{0}(z)B}, \qquad \text{where} \qquad R_{0}(z)\coloneqq (\mathscr{L}_{m}-z)^{-1},\]
for $z\in \rho(\mathscr{L}_{m})$.
According to Theorem~\ref{Theo Closed Exten} quoted in Appendix~\ref{Appendix Abstract Birman Schwinger}, if 
\begin{enumerate}[label=\textup{\textbf{H\arabic*:}}, ref=\textup{\textbf{H\arabic*}}]
\item \label{Ass non empty} there exists $z_{0}\in \rho(\mathscr{L}_{m})$ such that $-1\in \rho(Q(z_{0}))$,
\end{enumerate}
 then there exists a closed extension $\mathscr{L}_{m,V}\supseteq \mathscr{L}_m+V$, whose resolvent coincides with the family of bounded operators
\[ R_{0}(z)-\overline{R_{0}(z)B}\left( I+Q(z)\right)^{-1}AR_{0}(z).\] Our analysis below refers to the closed operator $\mathscr{L}_{m,V}$.

We begin Section~\ref{Sec L1 Perturbation} by showing that the abstract condition \ref{Ass non empty} holds, for a measurable function $V$ satisfying either of the following hypotheses;
\begin{enumerate}[label=\textup{\textbf{H\arabic*:}}, ref=\textup{\textbf{H\arabic*}}]
\setcounter{enumi}{1}
\item \label{Ass norm less than 1} $\Vert V \Vert_{L^1}<1$,
\item \label{Ass in Lp}$V \in L^1(\R,\C^{2\times 2})\cap L^p(\R,\C^{2\times 2})$ for some $p\in(1,\infty]$.
\end{enumerate}
Once we settle that, we formulate the next spectral stability result.

\begin{proposition}\label{Prop Perturbation}
Let $m\geq 0$ and let $\mathscr{L}_{m}$ be as in \eqref{Dirac Operator}. Let $V\in L^1(\R,\C^{2\times 2})$ and assume that either \ref{Ass norm less than 1} or \ref{Ass in Lp} holds true. 
Then, there exists a closed densely defined extension $\mathscr{L}_{m,V}\supseteq \mathscr{L}_{m}+V$ with a non-empty resolvent set and its spectrum is as follows.\\ If $m>0$,
\begin{enumerate}[label=\textup{\textbf{\arabic*)}}]
\item \label{Prop3.2 1} all the essential spectra of $\mathscr{L}_{m,V}$ are stable and coincide, namely
\[ \operatorname{Spec}_{\mathrm{ej}}(\mathscr{L}_{m,V})=\operatorname{Spec}_{\mathrm{ej}}(\mathscr{L}_{m})=S_{m}, \quad \forall \mathrm{j}\in [[1,5]]. \]
\end{enumerate}
If $m=0$, 
\begin{enumerate}[label=\textup{\textbf{\arabic*)}}]
\addtocounter{enumi}{1}
\item\label{Prop3.2 2} all the first four notions of essential spectra are stable and coincide, namely
	\[ \operatorname{Spec}_{\mathrm{ej}}(\mathscr{L}_{0,V})=\operatorname{Spec}_{\mathrm{ej}}(\mathscr{L}_{0})=\{z\in \C: |\Im z|= 1 \}, \quad \forall \mathrm{j}\in [[1,4]]. \] 
\end{enumerate}
Moreover,
\begin{enumerate}[label=\textup{\textbf{\arabic*)}}]
\addtocounter{enumi}{2}
\item\label{Prop3.2 3} If $V$ satisfies the assumption \ref{Ass norm less than 1}, we have
  $$\operatorname{Spec}_{\mathrm{p}}(\mathscr{L}_{0,V})\subset S_{0}.$$
\item \label{Prop3.2 4}If $V$ satisfies the assumption \ref{Ass in Lp} for $p<\infty$, we have
  $$\operatorname{Spec}_{\mathrm{p}}(\mathscr{L}_{0,V})\subset \left\{z\in \C: |\Im z|\leq 1+ \frac{2(p-1)}{p}\Vert V \Vert_{L^p}^{\frac{p}{p-1}} \right\}.$$ 
 If $V$ satisfies the assumption \ref{Ass in Lp} for $p=\infty$, then
  $$\operatorname{Spec}_{\mathrm{p}}(\mathscr{L}_{0,V})\subset \left\{z\in \C: |\Im z|\leq 1+ \Vert V \Vert_{L^{\infty}} \right\}.$$
\end{enumerate}
\end{proposition}

In Remark~\ref{Remark Compact Q}, we expand further on the justification for the condition~\ref{Ass norm less than 1} and compare with the case of the dislocated Schr{\"o}dinger operator. This assumption has already been documented in  \cite{Cuenin-Laptev-Tretter14} for $\mathcal{D}_m$ in place of $\mathscr{L}_m$. Noticeably, it  gives regions of inclusion for the perturbed eigenvalues. Our next result shows that a control on $\|V\|_{L^1}$ also gives inclusions for the spectrum and pseudospectrum of $\mathscr{L}_{m,V}$.

\begin{theorem}\label{inclusion pseudospectrum}
For $m\geq 0$, let
\[
    D\coloneqq \left\{ z\in\C\,:\,|\Im z|<\frac{3}{2},\,|\Re z|<\frac52 m\right\}.
\]
There exists a constant $0<C(m)\leq 1$ such that, if $\|V\|_{L^1}<C(m)$, then
\[
    \operatorname{Spec}_{\varepsilon}(\mathscr{L}_{m,V})\subseteq \overline{D} \cup \left\{z\in\C\,:\,|\Im z|\leq 1+\varepsilon \left(1+\frac{\|V\|_{L^1}}{4C(m)\big(C(m)-\|V\|_{L^1}\big)}\right)\right\},
\]
for all $\varepsilon \geq 0$. Moreover, for $C(0)=1$ this statement is valid. 
\end{theorem}

\medskip

According to \cite[Cor. 4.6]{Cuenin-Laptev-Tretter14}, $\mathcal{D}_m+V$ has point spectrum lying on a strip covering the real axis for $m=0$. According to \cite[Thm. 2.1]{Cuenin-Laptev-Tretter14}, for $m>0$, it lies inside two disks centered near the endpoints of the essential spectrum. The radius of these disks decreases as $\|V\|_{L^1}\to 0$.
Our next result shows that, rather unexpectedly, the point spectrum of $\mathscr{L}_{m,V}$ behaves differently from that of $\mathcal{D}_m+V$, under fast rates of decay of the perturbation at infinity. Namely, it escapes to infinity through the channel created by the numerical range, rather than concentrating at the end-points of the essential spectrum, as $\|V\|_{L^1}\to 0$. The condition we impose is analogous to the one considered in \cite[Thm. 2.3]{Henry-Krejcirik17} and in turn allows refined versions of the conclusions \ref{Prop3.2 3}~--~\ref{Prop3.2 4} of Proposition~\ref{Prop Perturbation} for $m>0$. 

Here and everywhere below we write
\begin{align*}
\nu_1'(x)= \sqrt{1+x^2},& \qquad \qquad &&\dd \nu_1(x)=\nu_1'(x)\dd x, \\
     \nu_2'(x)= \nu_1'(x)^2=1+x^2,& \qquad \qquad &&\dd\nu_2(x)=\nu_2'(x)\dd x.
\end{align*}
\begin{theorem}\label{Theo Location m>0}
Let $m>0$ be fixed. Then, there exists a constant $C(m)>0$, ensuring the following. For all $V\in L^1(\R,\C^{2\times 2}; \,\dd\nu_1)$ such that
\[
    \|V\|_{L^1(\dd\nu_1)}<C(m),
\]
we have
\[
\operatorname{Spec}_{\mathrm{p}} (\mathscr{L}_{m,V}) \subseteq \left\{z \in \C: |\Im z|\leq 1,\, |\Re(z)| > \frac{C(m)}{\|V\|_{L^1}^{\frac{1}{2}}}\right\}.
\]
\end{theorem}

Succinctly, this theorem says that for potentials $V$ decaying sufficiently rapidly at $\infty$, the discrete spectrum escapes to $\infty$ inside the instability band, $\{|\Im z|<1\}$, at a rate proportional to $\|V\|_{L^1}^{-\frac{1}{2}}$ in the regime $\|V\|_{L^1}\to 0$. In the conclusion, the norm is taken with respect to the Lebesgue measure.

For Schr{\"o}dinger operators,
a similar result was established in \cite[Thm. 2.3]{Henry-Krejcirik17}. In the latter, the density considered was $\nu'_2$ instead of $\nu_1'$, and  upper bounds for the norm of $Q(z)$ were obtained via the Hilbert-Schmidt norm.
 Here we employ the Schur Test instead, which allows a better control of the upper bounds for the norm of the integral operator $Q(z)$. Therefore, we manage to weaken the condition on the decay of $V$.

\medskip

Theorem~\ref{Theo Location m>0} opens the question of whether there exist $V\in L^1(\R,\C^{2\times 2})$ with non-empty discrete spectrum for $m>0$. We give an answer to this question in Subsection~\ref{secsteppot}. Our concrete finding in this respect is the following. Consider the step potential 
\[ V(x)=V_{a,b}(x)\coloneqq (-i\sgn(x)-b)\chi_{[-a,a]}(x)  I, \]
where $a>0$ and $b\in \R$.  Set $\mathrm{Dom}(\mathscr{L}_{m,V_{a,b}})=H^{1}(\R,\C^2)$. According to Proposition~\ref{Prop Perturbation}, the essential spectra of $\mathscr{L}_{m,V_{a,b}}$ are \[ \operatorname{Spec}_{\mathrm{ej}}(\mathscr{L}_{m,V_{a,b}})=\{z\in \C: |\Im z|=1,\, |\Re z| \geq m\},\]
for all $\mathrm{j}\in [[1,4]]$, including $\mathrm{j}=5$ if $m>0$.

\begin{proposition}\label{Prop Steplike}
For $m\geq 0$, let $\mathscr{L}_{m,V_{a,b}}$ be as in the previous paragraph. 
\begin{enumerate} [label=\textup{\textbf{\arabic*)}}]
\item If $m=0$, then the point spectrum of $\mathscr{L}_{0,V_{a,b}}$ is the band
\[ \operatorname{Spec}_{\mathrm{p}}(\mathscr{L}_{0,V_{a,b}})=\{z\in \C: |\Im z|<1\}.\]
\item If $m>0$, then the operator $\mathscr{L}_{m,V_{a,b}}$ has infinitely many isolated real eigenvalues accumulating at $\pm\infty$. 
\end{enumerate}
\end{proposition}

\medskip

In the context of Theorem~\ref{Theo Location m>0}, it is also natural
to examine the behaviour of the family of operators $\mathscr{L}_{m,\epsilon V}$ for fixed $m>0$, fixed potential $V$ and a small moving complex parameter $\epsilon\in\mathbb{C}$, in the regime $|\epsilon|\to 0$. This is the so-called \emph{weakly coupled model}. For this model we establish the following in Subsection~\ref{Subsec Improved dynamics}.

Let $V\in L^1(\R,\C^{2\times 2};\, \dd \nu_1)$. From Proposition \ref{Prop Perturbation} and Theorem \ref{Theo Location m>0}, it follows that
\[\operatorname{Spec}_{\mathrm{ej}}(\mathscr{L}_{m,\epsilon V})=\{z\in \C:  |\Im z|=1,\, |\Re z|\geq m\},\qquad \forall j\in [[1,5]]\] and, for a suitable constant $C(m)>0$,  
\begin{equation}\label{Weak Couple 0}
 \operatorname{Spec}_{\mathrm{dis}} (\mathscr{L}_{m,\epsilon V}) \subseteq \left\{z \in \C: |\Im z|<1,\, |\Re z| >  \frac{C(m)}{|\epsilon|^{\frac12}\|V\|_{L^1}^{\frac12}}\right\},
\end{equation}
for all $|\epsilon|$ small enough. If $V$ satisfies additional conditions of regularity and decay at infinity, then  Theorem~\ref{Theo Location m>0} refines as follows. 

\begin{theorem}\label{Theo Weak Coupling}
Let $m>0$. For any $V\in L^1(\R,\C^{2\times 2}; \dd \nu_2)\cap W^{2,1}(\R,\C^{2\times 2})$, 
there exists a constant $C(m)>0$ independent of $\epsilon$, such that  
\begin{equation*}
\operatorname{Spec}_{\mathrm{dis}}\left(\mathscr{L}_{m,\epsilon V}\right)\subset \left\{z \in \C: |\Im z|< 1,\, |\Re z| > \frac{C(m)}{|\epsilon|}\right\},
\end{equation*}
for all $|\epsilon|$ small enough.
\end{theorem}

\section{Resolvent and spectrum of $\mathscr{L}_m$}\label{Sec Res and Sp}

In this section we settle the framework of the differential expression  \eqref{Dirac Operator} and the associated linear operator. We begin by determining $\operatorname{Num} (\mathscr{L}_m)$. We then compute the Green matrix associated with the resolvent $(\mathscr{L}_m-z)^{-1}$ and find $\operatorname{Spec}(\mathscr{L}_m)$.

Since $\mathcal{D}_m^*=\mathcal{D}_m$ and multiplication by $i\sgn(x)I$ is bounded, it follows that the domain of $\mathscr{L}_m^*$ is also $H^1(\R,\C^2)$.  Let $\mathcal{T}$ be the antilinear operator of complex conjugation, $\mathcal{T} g=\overline{g}$ and $\mathcal{P}$ be the parity operator,  $(\mathcal{P}g)(x)=g(-x)$.  For $f\in H^{1}(\R,\C^2)$, direct substitution of the operators involved, gives
\begin{equation}\label{T self adjoint}
\begin{aligned}
(\mathcal{T} \mathscr{L}_{m} \mathcal{T}f)(x) =&\overline{(-i\partial_{x}) \sigma_{2}\overline{f(x)} +m\sigma_{3} \overline{f(x)} +i\sgn(x) \overline{f(x)}}\\
=&(-i\partial_{x})\sigma_{2}f(x)+m\sigma_{3}f(x)-i\sgn(x)f(x)= (\mathscr{L}_{m}^{*} f)(x)
\end{aligned}
\end{equation}
however,
\[  \mathcal{PT} \mathscr{L}_{m}f \neq \mathscr{L}_{m}\mathcal{PT}f.\]
Therefore, $\mathscr{L}_{m}$ is $\mathcal{T}$-selfadjoint but it is not $\mathcal{PT}$-symmetric.

Now we show that 
\begin{equation} \label{the numerical range}
    \operatorname{Num}(\mathscr{L}_{m})=\{z\in\mathbb{C}:|\Im z|\leq 1\}.
\end{equation}
Indeed, let $f=\begin{pmatrix}
f_{1}\\
f_{2}
\end{pmatrix}\in \mathrm{Dom}(\mathscr{L}_{m})$ be such that $\Vert f \Vert=1$. Integration by parts gives
\begin{equation*}
\begin{aligned}
\left\langle \mathscr{L}_{m}f,f \right\rangle = &m \left( \Vert f_{1} \Vert^2 - \Vert f_{2} \Vert^2\right)+2 \Re \langle f_{1}',f_{2}\rangle+i \left( \Vert f \Vert_{L^2(\R_{+})}^2 -\Vert f \Vert_{L^2(\R_{-})}^2\right).
\end{aligned}
\end{equation*}
Since $\Vert f \Vert=1$,  it follows immediately that $\vert \Im \left\langle \mathscr{L}_{m}f,f \right\rangle \vert\leq 1$, which implies that
$$\operatorname{Num}(\mathscr{L}_{m})\subseteq \{z\in\mathbb{C}:|\Im z|\leq 1\}.$$ Conversely, for the reverse inclusion, we construct appropriate test functions as follows. Choose smooth $g_{1}$ and $g_{2}$ with supports in $\R_{+}$ such that $\Vert g_{1} \Vert^2+\Vert g_{2} \Vert^2=1$ and $\Re \langle g_{1}',g_{2} \rangle>0$. For $\lambda>0$, let 
$$f_{j}(x)=\lambda^{1/2} g_{j}(\lambda x).$$
Then, $\|f\|=1$ and
\begin{equation*}
\begin{aligned}
\left\langle \mathscr{L}_{m}f,f \right\rangle = &m \left( \Vert g_{1} \Vert^2 - \Vert g_{2} \Vert^2\right)+2 \lambda\Re \langle g_{1}',g_{2}\rangle+i.
\end{aligned}
\end{equation*}
Hence, by increasing $\lambda$, it follows that $\R_{+}+i\subseteq \operatorname{Num}(\mathscr{L}_{m})$. Similar arguments, choosing  $g_{1}$ and $g_{2}$ with support in $\R_{-}$ or changing the sign of $\Re \langle g_{1}',g_{2}\rangle$, give $\R \pm i\subseteq \operatorname{Num}(\mathscr{L}_{m})$. The convexity of the numerical range completes the proof of \eqref{the numerical range}.

Below, we often write $\Sigma=\operatorname{Int}\left(\operatorname{Num}\mathscr{L}_m\right)$ and call this set the \emph{instability band}.


\subsection{Integral form of the resolvent} \label{Subsec Resolvent}
In this subsection we compute explicitly the integral kernel associated to the resolvent of $\mathscr{L}_{m}$. We also establish preliminary estimates on its  norm. We begin by fixing a notation that will prove crucial for simplifying our expressions and subsequent analysis. 

For $z\in \C\setminus \{\pm m\pm i\}$, we write
\begin{equation}\label{Notation pm}
 \begin{aligned}
& \mu_{z}^{-}\coloneqq \sqrt{(m+i+z)(m-i-z)}, \qquad && \mu_{z}^{+}\coloneqq \sqrt{(m-i+z)(m+i-z)},\\
& w_{z}^{-}\coloneqq  \frac{m-i-z}{\mu_{z}^{-}}=\frac{\sqrt{m-i-z}}{\sqrt{m+i+z}},\qquad && w_{z}^{+}\coloneqq \frac{m+i-z}{\mu_{z}^{+}}=\frac{\sqrt{m+i-z}}{\sqrt{m-i+z}}.
\end{aligned}
\end{equation}
Here and throughout the paper, we pick the principal branch $z\mapsto \sqrt{z}=z^{\frac12}$ defined on $\C$, holomorphic on $\C\setminus (-\infty,0]$ and having a positive imaginary part on $(-\infty,0)$. This choice ensures that $w_{z}^{\pm}$ have the stated expressions in \eqref{Notation pm}.

Since
\begin{equation*}
\begin{aligned}
(m\pm i+z)(m\mp i-z)\in (-\infty,0]   & \quad\iff\quad  |\Re z| \geq m \text{ and } \Im z=\mp 1, 
\end{aligned}
\end{equation*}
then
\begin{equation}\label{Cond mu pm}
\Re \mu_{z}^{\pm}>0 \quad \iff \quad z\not\in\left\{\tau\pm i:\tau\in\mathbb{R},\, |\tau|\geq m\right\}. 
\end{equation}
Moreover,  $\overline{\mu^{\pm}_z}=\mu^{\mp}_{\overline{z}}$ and $\overline{w^{\pm}_z}=w^{\mp}_{\overline{z}}$ for all $z$ satisfying \eqref{Cond mu pm}. Below we will use these facts repeatedly, often without further explicit mention.

Recall $S_{m}$ given by \eqref{Set S_m}. In the next subsection, we will establish that $\operatorname{Spec}(\mathscr{L}_m)=S_{m}$. For $z\in \mathbb{C}\setminus S_m$, define the matrices
\begin{align*}
&N_{1,z}\coloneqq  \frac{1}{2}\frac{w_z^+ - w_z^-}{w_z^+ + w_z^-}\begin{pmatrix} 1/w_z^- & 1 \\ 1 & w_z^-  \end{pmatrix}, &&N_{6,z} \coloneqq  \frac{1}{2}\begin{pmatrix} -1/w_z^+ & 1  \\ -1 & w_z^+ \end{pmatrix},\\
&N_{2,z}\coloneqq  \frac{1}{2}\begin{pmatrix} -1/w_z^- & -1 \\ 1 & w_z^- \end{pmatrix}, &&N_{7,z}\coloneqq  N_{4,z},\\
&N_{3,z}\coloneqq  \frac{1}{w_z^+ + w_z^-}\begin{pmatrix} -1 & -w_z^- \\ w_z^+ & w_z^+ w_z^- \end{pmatrix}, &&N_{8,z}\coloneqq  \frac{1}{w_z^+ + w_z^-}\begin{pmatrix} -1 & w_z^+ \\ -w_z^- & w_z^+ w_z^- \end{pmatrix},\\
&N_{4,z}\coloneqq  \frac{1}{2}\frac{w_z^+ - w_z^-}{w_z^+ + w_z^-}\begin{pmatrix} -1/w_z^+ & 1 \\ 1 & -w_z^+ \end{pmatrix}, &&N_{9,z}\coloneqq  \frac{1}{2}\begin{pmatrix} -1/w_z^- & 1 \\ -1 & w_z^- \end{pmatrix},\\
&N_{5,z}\coloneqq  \frac{1}{2}\begin{pmatrix} -1/w_z^+ & -1 \\ 1 & w_z^+ \end{pmatrix}, &&N_{10,z}\coloneqq  N_{1,z}.
\end{align*}
For $x,y\in \mathbb{R}$, $x\neq y$, let the matrix kernel $\mathcal{R}_{z}(x,y)$ be given in regions of the plane $\mathbb{R}^2$ by
\begin{equation}\label{G}
\mathcal{R}_{z}(x,y)\coloneqq  \left\{
\begin{aligned}
			 & N_{1,z} \,e^{\mu^-_z(x+y)}+ N_{2,z}\, e^{-\mu^-_z(x-y)}, && \text{for } y< x < 0;\\ 
			 &  N_{3,z}\, e^{-\mu^+_zx+\mu^-_zy}, && \text{for } y < 0 < x; \\
			 & N_{4,z}\,e^{-\mu^+_z(x+y)} +N_{5,z}\,e^{-\mu^+_z(x-y)},  && \text{for } 0< y < x; \\
			 & N_{6,z}\,e^{\mu^+_z(x-y)} + N_{7,z}\,e^{-\mu^+_z(x+y)},  && \text{for } 0< x< y; \\
			 & N_{8,z} \, e^{\mu^-_zx-\mu^+_zy}, \qquad && \text{for } x< 0 < y; \\
			 & N_{9,z}\,e^{\mu^-_z(x-y)}  + N_{10,z}\,e^{\mu^-_z(x+y)}, && \text{for } x < y < 0.
		\end{aligned}
\right.
\end{equation}
Note that $\mathcal{R}_z$ is continuous across the lines $x=0$ and $y=0$, but discontinuous across the line $x=y$. Here it is not important what values we give to this matrix function on these lines, so we leave them unassigned. 

We also mention that, in contrast to the case of the resolvent kernel of the Schr{\"o}dinger operator $-\frac{\dd^2}{\dd x^2}+i\sgn(x)$, considered in \cite[Prop. 3.1]{Henry-Krejcirik17}, $\mathcal{R}_z$ is not symmetric with respect to the line $x=y$.  We will see in Lemma~\ref{Lemma Norm of G} that the induced operator norm of $\mathcal{R}_{z}$ is nonetheless symmetric. That is, $\vert \mathcal{R}_{z}(x,y)\vert_{\C^{2\times 2}}=\vert \mathcal{R}_{z}(y,x)\vert_{\C^{2\times 2}}$ for all $(x,y)\in \R^2$. 

As we shall see from the next statements, $\mathcal{R}_{z}$ is the (matrix) Green function associated to the operator $\mathscr{L}_{m}-z$. This explicit formula has a far reaching role in the analysis that we conduct below. Remarkably, we will see that it enables the systematic study of perturbation of $\mathscr{L}_m$ via a Birman-Schwinger principle, as described in Section \ref{Sec L1 Perturbation}.

\begin{proposition}\label{Prop Resolvent}
For $m\geq 0$, let $\mathscr{L}_{m}$ be the operator given by \eqref{Dirac Operator}. 
Then, $
\C\setminus S_{m} \subset \rho\left(\mathscr{L}_{m}\right)$. Moreover,
\begin{equation}\label{Integral Op}
\left[\left(\mathscr{L}_{m}-z\right)^{-1} f\right](x) = \int_{\R} \mathcal{R}_{z}(x,y) f(y)\, \dd y ,
\end{equation}
for all $z \in \C\setminus S_{m}$ and $f\in L^{2}(\R,\C^2)$,
\end{proposition}

Before proceeding to the proof, we highlight that the resolvent kernel of the free Dirac operator \cite[Proof of Thm. 2.1]{Cuenin-Laptev-Tretter14}, is a convolution kernel and it contains only exponential terms involving the difference $x-y$. By contrast, for $m>0$, $\mathcal{R}_z$ also includes mixed exponential terms and it cannot be expressed as a pure tensor product of two matrix functions. We will see below that this more complex structure of the Green function is responsible for the asymptotic growth of the resolvent norm inside the numerical range. 

\medskip

The remainder of this section is devoted to the proof of Proposition~\ref{Prop Resolvent}. Although this proof follows a routine argumentation, we include crucial details that will aid in the computation of sharp estimates for the resolvent norm. We split it into 4 steps.

\subsubsection*{Step 1} Fix $z \in \C\setminus S_{m}$ and $f\in L^{2}(\R,\C^2)$. Consider the eigenvalue equation \[(\mathscr{L}_{m}-z) u = f,\] which in matrix form reads
\begin{equation} \label{Spectral Eq}
\begin{pmatrix}
m+i\sgn(x)-z & -\partial_x\\
\partial_x & -m +i\sgn(x)-z
\end{pmatrix} \begin{pmatrix}
u_{1}(x)\\
u_{2}(x)
\end{pmatrix} = \begin{pmatrix}
f_{1}(x)\\
f_{2}(x)
\end{pmatrix}.
\end{equation}
Our final goal in the proof is to determine the explicit integral expression for the solution. Multiplying by the matrix $J\coloneqq  \begin{pmatrix}
0 & -1\\
1 & 0
\end{pmatrix}$ on the left of the two sides of \eqref{Spectral Eq}, gives the equivalent non-homogeneous linear system of first order, 
\begin{equation}\label{ODE}
\partial_{x} u(x) = A_{z}(x)u(x)+g(x),
\end{equation}
where
\[
	A_z(x) \coloneqq  \begin{pmatrix}
	0 & m-i \sgn(x) + z \\ m+i \sgn(x) - z & 0
\end{pmatrix} \qquad \text{ and } \qquad g(x)\coloneqq  -Jf(x).
\]
We now seek for the solution of this system.

Consider the equation \eqref{ODE} for $x>0$ and $x<0$, where the corresponding solutions are indicated as $u_{-}$ and $u_{+}$. Let
\[
	A_z^{-} \coloneqq  \begin{pmatrix}
	0 & m + i + z \\ m - i- z & 0
\end{pmatrix} \quad \text{and} \quad  A_z^{+} \coloneqq  \begin{pmatrix}
	0 & m - i + z \\ m + i- z & 0
\end{pmatrix},
\]
be the expressions for $A_z(x)$ for $x>0$ and $x<0$, respectively.
Let $\Phi_{z}^{\pm}(x)$ be fundamental matrices for $\partial_{x}u_{\pm}=A_{z}^{\pm} u^{\pm}$. The columns of $\Phi_{z}^{\pm}(x)$ ought to be linear independent for all $x\in\mathbb{R}$, thus, we can pick the choice
\begin{equation}\label{Der Phi}
\Phi_{z}^{\pm}(x)=\exp(A_{z}^{\pm} x) =\begin{pmatrix}
\cosh( \mu_{z}^{\pm} x) && \frac{\sinh( \mu_{z}^{\pm} x)}{w_{z}^{\pm}}\\
w_{z}^{\pm}\sinh( \mu_{z}^{\pm} x) && \cosh( \mu_{z}^{\pm} x)
\end{pmatrix}.
\end{equation}
Note that
\begin{gather} \label{semig}
\left(\Phi_{z}^{\pm}(x) \right)^{-1}=\Phi_{z}^{\pm}(-x)\qquad \text{and} \qquad\Phi_{z}^{\pm}(x)\Phi_{z}^{\pm}(y)=\Phi_{z}^{\pm}(x+y).
\end{gather}

We now determine the solutions $u_{\pm}$ of the inhomogeneous system
\begin{equation}\label{ODE pm}
\partial_{x}u_{\pm}(x) = A_{z}^{\pm} u_{\pm}(x) +g(x) \qquad\text{on} \quad \R_{\pm}
\end{equation}
in the form $u_{\pm}(x) = \Phi_{z}^{\pm}(x)v_{\pm}(x)$, where $v_{\pm}:\R_{\pm}\to \C^2$ are to be determined. Taking the derivative of $u_{\pm}$ and using \eqref{Der Phi}, we obtain
\begin{align*}
\partial_{x} u_{\pm}(x)=& A_{z}^{\pm} \Phi_{z}^{\pm}(x) v_{\pm}(x)+\Phi_{z}^{\pm}(x)\,\partial_{x}v_{\pm}(x)
= A_{z}^{\pm} u_{\pm}(x)+\Phi_{z}^{\pm}(x)\,\partial_{x}v_{\pm}(x) \qquad\text{on} \quad \R_{\pm}.
\end{align*}
Therefore, $u_{\pm}$ are solutions of \eqref{ODE pm} if and only if $ \Phi_{z}^{\pm}(x)\,\partial_{x}v_{\pm}(x)=g(x)$ on $\R_{\pm}$. Hence, we choose
\[ v_{\pm}(x)= \int_{0}^{x} \left(\Phi_{z}^{\pm}(y)\right)^{-1} g(y)\dd y+ \begin{pmatrix}
\alpha_{\pm}\\
\beta_{\pm}
\end{pmatrix},\]
where $\alpha_{\pm}, \beta_{\pm} \in \C$ are constants to be determined.
Then, the general solutions of \eqref{ODE pm} are given by
\[ u_{\pm}(x) = \Phi_{z}^{\pm}(x)\int_{0}^{x} \left(\Phi_{z}^{\pm}(y)\right)^{-1} g(y)\dd y+ \Phi_{z}^{\pm}(x)\begin{pmatrix}
\alpha_{\pm}\\
\beta_{\pm}
\end{pmatrix}.\]
As the condition 
\[
     \lim_{x\to 0^-}u_-(x)=\lim_{x\to 0^+}u_+(x),
\]
is needed for a solution to be in the domain of the operator, then necessarily
\[ \begin{pmatrix}
\alpha_{-}\\
\beta_{-}
\end{pmatrix}=\begin{pmatrix}
\alpha_{+}\\
\beta_{+}
\end{pmatrix}\eqqcolon\begin{pmatrix}
\alpha\\
\beta
\end{pmatrix}.\]
This gives the appropriate regularity to the solution for all $x\in\mathbb{R}$.

\subsubsection*{Step 2} Our next objective is to determine $\alpha$ and $\beta$ such that $u_{\pm}$ decay to zero at infinity. According to \eqref{semig}, we have 
\begin{align*}
 u_{\pm}(x) &= \int_{0}^{x} \Phi_{z}^{\pm}(x-y)  g(y)\, \dd y +\Phi_{z}^{\pm}(x)\begin{pmatrix}
\alpha\\
\beta
\end{pmatrix}\\
&=\int_{0}^{x} \begin{pmatrix}
\cosh( \mu_{z}^{\pm} (x-y)) && \frac{\sinh( \mu_{z}^{\pm} (x-y))}{w_{z}^{\pm}}\\
w_{z}^{\pm}\sinh( \mu_{z}^{\pm} (x-y)) && \cosh( \mu_{z}^{\pm} (x-y))
\end{pmatrix}g(y)\, \dd y +\begin{pmatrix}
\cosh( \mu_{z}^{\pm} x) && \frac{\sinh( \mu_{z}^{\pm} x)}{w_{z}^{\pm}}\\
w_{z}^{\pm}\sinh( \mu_{z}^{\pm} x) && \cosh( \mu_{z}^{\pm} x)
\end{pmatrix} \begin{pmatrix}
\alpha\\
\beta
\end{pmatrix}.
\end{align*}
Writing the right-hand side in exponential form, reduces to
\begin{equation} \label{u pm}
\begin{aligned}
u_{\pm}(x)= &\int_{0}^{x} e^{\mu_{z}^{\pm}(x-y)}S^{\pm}_{z}g(y) \dd y +  e^{\mu_{z}^{\pm}x}S^{\pm}_{z}\begin{pmatrix}
\alpha\\
\beta
\end{pmatrix}+\int_{0}^{x} e^{-\mu_{z}^{\pm}(x-y)}T_{z}^{\pm}g(y) \dd y +  e^{-\mu_{z}^{\pm}x}T_{z}^{\pm}\begin{pmatrix}
\alpha\\
\beta
\end{pmatrix},
\end{aligned}
\end{equation}
where
\begin{equation}\label{Matrixes ST}
S^{\pm}_{z}\coloneqq  \frac{1}{2}\begin{pmatrix}
1 && 1/w_{z}^{\pm}\\
w_{z}^{\pm} && 1
\end{pmatrix}\qquad \text{and} \qquad T_{z}^{\pm}\coloneqq  \frac{1}{2}\begin{pmatrix}
1 && -1/w_{z}^{\pm}\\
-w_{z}^{\pm} && 1
\end{pmatrix}.
\end{equation}
We seek for conditions ensuring $\lim_{x\to \pm\infty}u_{\pm}(x)=0$. Since the function $g$ belongs to $L^2$ but not necessarily to $L^1$, we cannot use directly the Dominated Convergence Theorem. Instead, we appeal to the density of $C_{c}^{\infty}(\R,\C^2)$ in $L^{2}(\R,\C^2)$ as follows. For $\epsilon>0$ small, let $g_\epsilon \in C_{c}^{\infty}(\R,\C^2)$ be such that $\|g-g_\epsilon\|<\epsilon$.
Then,
\begin{align*}
\left\vert \int_{0}^{x} e^{-\mu_{z}^{+}(x-y)} T^{+}_{z} g(y)\,\dd y\right\vert_{\C^2}\leq &\int_{0}^{x} e^{-\Re \mu_{z}^{+}(x-y)} \left\vert T^{+}_{z} \right\vert_{\C^{2\times 2}} \vert g(y)-g_{\epsilon}(y)\vert_{\C^2}\, \dd y\\
 &+\int_{0}^{x} e^{-\Re \mu_{z}^{+}(x-y)} \left\vert T^{+}_{z} \right\vert_{\C^{2\times 2}} \vert g_{\epsilon}(y)\vert_{\C^2}\, \dd y\\
 \leq & \left\vert T^{+}_{z} \right\vert_{\C^{2\times 2}} \sqrt{\frac{1-e^{-2\Re \mu_{z}^{+}x}}{2\Re \mu_{z}^{+}}} \Vert g-g_{\epsilon} \Vert\\
 &+\left\vert T^{+}_{z} \right\vert_{\C^{2\times 2}}\int_{0}^{x} e^{-\Re \mu_{z}^{+}(x-y)}  \vert g_{\epsilon}(y)\vert_{\C^2}\, \dd y.
\end{align*}
Recall \eqref{Cond mu pm}. Since $\Re \mu_{z}^{+}>0$, the integrand $e^{-\Re \mu_{z}^{+}(x-y)}  \vert g_{\epsilon}(y)\vert_{\C^2}$ is bounded by the $L^1$ function $\vert g_{\epsilon}\vert_{\C^2}$, then, according to the Dominated Convergence Theorem, we have
\[ \lim_{x\to+\infty} \left\vert \int_{0}^{x} e^{-\mu_{z}^{+}(x-y)} T^{+}_{z} g(y)\,\dd y\right\vert_{\C^2} \leq \left\vert T^{+}_{z} \right\vert_{\C^{2\times 2}} \frac{\epsilon}{\sqrt{\Re \mu_{z}^{+}}} .\]
Because $\epsilon$ is arbitrary, the limit on the left-hand side is zero. In turn, we gather that 
\begin{equation}\label{ab1}
\begin{aligned}
\lim_{x \to +\infty} u_{+}(x)=0 \quad &\iff \quad \lim_{x\to +\infty}
e^{\mu_{z}^{+}x}\left(\int_{0}^{x} e^{-\mu_{z}^{+}y}S_{z}^{+}g(y) \dd y +  S_{z}^{+}\begin{pmatrix}
\alpha\\
\beta
\end{pmatrix}\right)=0 \\
&\iff \quad S_{z}^{+}\begin{pmatrix}
\alpha\\
\beta
\end{pmatrix}=-\int_{0}^{+\infty} e^{-\mu_{z}^{+}y}S_{z}^{+}g(y)\, \dd y
\\ &\iff \quad \alpha+\frac{1}{w_{z}^{+}}\beta= -\int_{0}^{+\infty} \left(g_{1}(y)+\frac{1}{w_{z}^{+}}g_{2}(y)\right)e^{-\mu_{z}^{+}y}\, \dd y.
\end{aligned}
\end{equation}
For the latter note that $\det S_{z}^{+}=0$. 

Similarly, by considering the limit to $-\infty$, we conclude that
\begin{equation}\label{ab2}
\lim_{x\to-\infty} u_{-}(x) =0 \iff \alpha-\frac{1}{w_{z}^{-}}\beta= \int_{-\infty}^{0} \left(g_{1}(y)-\frac{1}{w_{z}^{-}}g_{2}(y) \right)e^{\mu_{z}^{-}y}\, \dd y.
\end{equation}
Now, \eqref{ab1}-\eqref{ab2} can be re-written as a single linear system of equations,
\begin{equation}\label{Eq of ab}
\begin{aligned}
\begin{pmatrix}
1 && 1/w_{z}^{+}\\
1 && -1/w_{z}^{-}
\end{pmatrix}\begin{pmatrix}
\alpha\\
\beta
\end{pmatrix}= &\int_{0}^{+\infty} e^{-\mu_{z}^{+}y} \begin{pmatrix}
-1 && -1/w_{z}^{+}\\
0 && 0
\end{pmatrix} g(y)\, \dd y
\\
&+\int_{-\infty}^{0} e^{\mu_{z}^{-}y} \begin{pmatrix}
0 && 0\\
1 && -1/w_{z}^{-}
\end{pmatrix} g(y)\, \dd y.
\end{aligned}
\end{equation}
Since \[\det \begin{pmatrix}
1 && 1/w_{z}^{+}\\
1 && -1/w_{z}^{-}
\end{pmatrix}\not=0\] for all $z\in \C\setminus S_{m}$, 
then
\begin{align*}
\begin{pmatrix}
\alpha\\
\beta
\end{pmatrix} = &\int_{0}^{+\infty} e^{-\mu_{z}^{+}y} M_{z}^{+} g(y) \, \dd y +\int_{-\infty}^{0} e^{\mu_{z}^{-}y}M_{z}^{-}g(y)\, \dd y,
\end{align*}
where
\[ M_{z}^{+}=\frac{-1}{w_z^+ + w_z^-} \begin{pmatrix}
			w_z^+ & 1 \\ w_z^+ w_z^- & w_z^-
		\end{pmatrix},\qquad M_{z}^{-}=\frac{1}{w_z^+ + w_z^-} \begin{pmatrix}
			w_z^{-} & -1 \\ -w_z^+ w_z^- & w_z^+
		\end{pmatrix}.\]

\subsubsection*{Step 3} We now verify the expression for the matrix integral kernel $\mathcal{R}_z$.
Replacing the constant vector $\begin{pmatrix}
\alpha\\
\beta
\end{pmatrix}$ into the formula for $u_{\pm}$, gives
		\begin{align*}
		u_{\pm}(x)= &\int_{0}^{x} e^{\mu_{z}^{\pm}(x-y)}S^{\pm}_{z}g(y) \dd y +  \int_{0}^{+\infty} e^{\mu_{z}^{\pm}x-\mu_{z}^{+}y}S^{\pm}_{z} M_{z}^{+} g(y) \, \dd y \\
&+\int_{-\infty}^{0} e^{\mu_{z}^{\pm}x+\mu_{z}^{-}y}S^{\pm}_{z} M_{z}^{-} g(y) \, \dd y+\int_{0}^{x} e^{-\mu_{z}^{\pm}(x-y)}T_{z}^{\pm}g(y) \dd y \\
&+  \int_{0}^{+\infty} e^{-\mu_{z}^{\pm}x-\mu_{z}^{+}y}T^{\pm}_{z} M_{z}^{+} g(y) \, \dd y+\int_{-\infty}^{0} e^{-\mu_{z}^{\pm}x+\mu_{z}^{-}y}T^{\pm}_{z} M_{z}^{-} g(y) \, \dd y.
\end{align*}
Since $S_{z}^{+}M_{z}^{+}=-S_{z}^{+}$, $S_{z}^{+}M_{z}^{-}=0$ and $g=-Jf$, then
\begin{equation*}
\begin{aligned}
u_{+}(x) = &\int_{x}^{+\infty} N_{6,z}e^{\mu_z^{+}(x-y)}f(y)\, \dd y+\int_{0}^{x} N_{5,z}e^{-\mu_{z}^{+}(x-y)} f(y) \dd y \\
&+\int_{0}^{+\infty} N_{4,z}e^{-\mu_{z}^{+}(x+y)}f(y) \, \dd y+\int_{-\infty}^{0}N_{3,z}e^{-\mu_{z}^{+}x+\mu_{z}^{-}y} f(y)\, \dd y.
\end{aligned}	
\end{equation*}
Similarly, since $T_{z}^{-}M_{z}^{-}=T_{z}^{-}$ and $T_{z}^{-}M_{z}^{+}=0$, then
\begin{align*}
u_{-}(x) = &\int_{-\infty}^{x} N_{2,z}e^{-\mu_z^{-}(x-y)}f(y)\, \dd y+\int_{x}^{0} N_{9,z}e^{\mu_{z}^{-}(x-y)} f(y) \dd y \\
&+\int_{-\infty}^{0} N_{1,z}e^{\mu_{z}^{-}(x+y)}  f(y) \, \dd y +\int_{0}^{+\infty} N_{8,z}e^{\mu_{z}^{-}x-\mu_{z}^{+}y} f(y)\, \dd y.
\end{align*}

\subsubsection*{Step 4} So far, we have obtained a solution to \eqref{Spectral Eq} on $\R$, which is continuous and decaying at infinity. This solution has the form
\[ u(x)\coloneqq  \left\{\begin{aligned}
&u_{+}(x) \qquad &&\text{if } x\geq 0,\\
&u_{-}(x) \qquad &&\text{if } x\leq 0,
\end{aligned} \right.\]
and has the integral representation
\begin{equation} \label{Integral rep u} u(x)=\int_{\R} \mathcal{R}_{z}(x,y) f(y)\, \dd y\end{equation}
as in the statement of Proposition~\ref{Prop Resolvent}. To complete the proof, it remains to check that $u \in L^{2}(\R,\C^2)$ and thus $\mathscr{L}_{m}u = zu+f \in L^{2}(\R,\C^2)$, so that $u\in H^{1}(\R,\C^2)$. This will be the consequence of the next two lemmas, which will be useful in their own right later on.

Firstly, we explicitly compute the matrix norm of the kernel.

\begin{lemma}\label{Lemma Norm of G}
Let $\mathcal{R}_{z}$ be as defined in \eqref{G}. Assume that $z \in \C\setminus S_{m}$ is fixed. Let
\begin{equation}\label{varphi}
\varphi_{z}(t)\coloneqq  \left\{
\begin{aligned}
&\frac{\sqrt{1+\vert w_{z}^{-}\vert^2}}{2}\left(\frac{1}{\vert w_{z}^{-}\vert^2}\left\vert k_{z}e^{2\mu_{z}^{-}t} -1\right\vert^2+\left\vert k_{z}e^{2\mu_{z}^{-}t} +1\right\vert^2 \right)^{1/2}, &&\text{ for }t\leq 0,\\
&\frac{\sqrt{1+\vert w_{z}^{+}\vert^2}}{2}\left(\frac{1}{\vert w_{z}^{+}\vert^2}\left\vert k_{z}e^{-2\mu_{z}^{+}t} +1\right\vert^2+\left\vert k_{z}e^{-2\mu_{z}^{+}t} -1\right\vert^2 \right)^{1/2},&& \text{ for }t\geq 0,
\end{aligned}
\right.
\end{equation}
where
\begin{equation}\label{Function k}
k_{z}\coloneqq \frac{w_z^+ - w_z^-}{w_z^+ + w_z^-}.
\end{equation}
Then, $\varphi_{z}:\R \longrightarrow \R_{+}$ is bounded and continuous, and
\begin{equation}\label{Norm of G}
\vert \mathcal{R}_{z}(x,y) \vert_{\C^{2\times 2}}= \left\{
\begin{aligned}
&\varphi_{z}(\max\{x,y\})\,e^{-\Re \mu_{z}^{-}\vert x-y \vert}, \qquad &&\text{for } x< 0, \,y < 0,\\
&\varphi_{z}(\min\{x,y\})\,e^{-\Re \mu_{z}^{+}\vert x-y \vert}, \qquad &&\text{for } x> 0, \,y > 0, \\
&\varphi_{z}(0)\, e^{- \Re \mu_{z}^{+}\vert x \vert-\Re \mu_{z}^{-} \vert y \vert}, \qquad && \text{for } x> 0, \,y < 0,\\
&\varphi_{z}(0)\, e^{- \Re \mu_{z}^{-}\vert x \vert-\Re \mu_{z}^{+} \vert y \vert}, \qquad && \text{for } x< 0, \,y > 0.\\
\end{aligned}
\right.
\end{equation} 
\end{lemma}
\begin{proof}
The function $\varphi_z$ is continuous, since
\[ \lim_{t\to 0_{-}}\varphi_{z}(t)=\frac{\sqrt{(1+|w^+_z|^2) (1+|w^-_z|^2)}}{|w^+_z+w^-_z|}=\lim_{t\to 0^{+}}\varphi_{z}(t).\]
It is bounded, since
\begin{equation}\label{varphi bounded}
\begin{aligned}
\sup_{t\in \R_{\pm}}\varphi_{z}(t)  \leq & \left(|k_{z}|+1\right) \frac{\sqrt{1+\vert w_{z}^{\pm}\vert^2}}{2}\sqrt{1+\frac{1}{|w_{z}^{\pm}|^2}}= \frac{1}{2}\left(|k_{z}|+1\right) \left( |w_{z}^{\pm}| +\frac{1}{|w_{z}^{\pm}|}\right).
\end{aligned}
\end{equation}
Now, recall that, for a singular matrix $A\in \C^{2\times 2}$, the operator norm and the Frobenius norm coincide
\begin{equation*}
\left\vert A \right\vert_{C^{2 \times 2}}=\displaystyle\sqrt{\sum_{1\leq i,j \leq 2} \vert A_{ij} \vert^2}.
\end{equation*}
Indeed, $\det (A^{*}A)=0$, therefore the two non-negative eigenvalues of $A^{*}A$ are $0$ and $\mathrm{Tr}\left( A^*\, A\right)$. Applying this, alongside a straightforward computation, gives 
\begin{equation*}
\vert \mathcal{R}_{z}(x,y) \vert_{\C^{2\times 2}}= \left\{
\begin{aligned}
&\frac{\sqrt{1+|w_{z}^{\pm}|^2}}{2} \left(\frac{1}{\vert w_{z}^{\pm} \vert^2} \left\vert k_{z}e^{-\mu_{z}^{\pm}\vert x+y \vert}\pm e^{-\mu_{z}^{\pm}\vert x-y \vert}\right\vert^2\right.\\
&\hspace{2.3 cm}+ \left.\left\vert k_{z} e^{-\mu_{z}^{\pm}\vert x+y \vert}\mp e^{-\mu_{z}^{\pm}\vert x-y \vert}\right\vert^2\right)^{1/2},&&\text{for }  \pm x> 0 \text{ and }\pm y>0, \\
&\frac{\sqrt{(1+|w^+_z|^2) (1+|w^-_z|^2)}}{|w^+_z+w^-_z|} e^{- \Re \mu_{z}^{\pm}\vert x \vert-\Re \mu_{z}^{\mp} \vert y \vert}, &&\text{for }  \pm x> 0 \text{ and }\pm y<0.
\end{aligned}
\right.
\end{equation*}
Factoring out suitable exponential terms and substitution of the expression \eqref{varphi} yields \eqref{Norm of G}.
\end{proof}

The next lemma directly implies that the function $u$ given by \eqref{Integral rep u} belongs to $L^{2}(\R,\C^2)$, and hence $z\in \rho(\mathscr{L}_m)$. This conclusion completes the proof of Proposition~\ref{Prop Resolvent}. Additionally, the lemma provides an initial estimate for the upper bound of the resolvent norm of the operator $\mathscr{L}_m$. This estimate is derived directly from the explicit expression of 
$|\mathcal{R}_z(x,y)|_{\C^{2\times 2}}$ in terms of the supremum of $\varphi_z(t)$, using the Schur Test. We will refine this upper bound in the next section.

\begin{lemma}\label{Lem Upper bound}
Let $z \in \C\setminus S_{m}$. Then, the integral operator on the right-hand side of \eqref{Integral rep u} is bounded on $L^{2}(\R,\C^2)$, and
\begin{equation}\label{Upper bound}
\Vert (\mathscr{L}_{m}-z)^{-1} \Vert \leq  \frac{2}{\min \left\{\Re \mu_{z}^{-}, \Re \mu_{z}^{+}\right\}}\sup_{t\in \R} \varphi_{z}(t).
\end{equation}
\end{lemma}
\begin{proof}
We will prove this lemma by the Schur Test. Our goal is to show that
\begin{equation}\label{Schur cond}
\begin{aligned}
&\int_{\R} \vert \mathcal{R}_{z}(x,y) \vert_{\C^{2\times 2}} \, \dd y \leq \frac{2}{\min \left\{\Re \mu_{z}^{-}, \Re \mu_{z}^{+}\right\}}\sup_{t\in \R} \varphi_{z}(t),\qquad \text{ for a.e. } x\in \R, \\
&\int_{\R} \vert \mathcal{R}_{z}(x,y) \vert_{\C^{2\times 2}} \, \dd x \leq \frac{2}{\min \left\{\Re \mu_{z}^{-}, \Re \mu_{z}^{+}\right\}}\sup_{t\in \R} \varphi_{z}(t),\qquad \text{ for a.e. } y\in \R.
\end{aligned}
\end{equation}
Since $\vert \mathcal{R}_{z}(x,y) \vert_{\C^{2\times 2}}$ is symmetric with respect to swapping $x$ and $y$, it suffices to show the first identity in \eqref{Schur cond}. 

For $x\leq 0$,  \eqref{Norm of G} yields
\begin{align*}
&\int_{\R} |\mathcal{R}_{z}(x,y)|_{\mathbb{C}^{2 \times 2}} \, \mathrm{d} y \\
			& = \int_{-\infty}^{x} \varphi_{z}(x)  e^{\Re \mu_{z}^{-}(y-x)}\, \dd y + \int_{x}^{0} \varphi_{z}(y)  e^{\Re \mu_{z}^{-}(x-y)}\, \dd y+ \int_{0}^{+\infty} \varphi_{z}(0) e^{\Re \mu_{z}^{-} x -\Re \mu_{z}^{+}y }\, \dd y.\\
			&\leq \varphi_{z}(x) \frac{1}{\Re \mu_{z}^{-}} +\left(\sup_{t\leq 0} \varphi_{z}(t)\right)\frac{1-e^{\Re \mu_{z}^{-} x}}{\Re \mu_{z}^{-}} +\varphi_{z}(0)\frac{e^{\Re \mu_{z}^{-} x}}{\Re \mu_{z}^{+}}\\
			&\leq \left(\sup_{x\leq 0} \varphi_{z}(x)\right)\left[ \frac{2}{\Re \mu_{z}^{-}}+\left(\frac{1}{\Re \mu_{z}^{+}}-\frac{1}{\Re \mu_{z}^{-}}\right)e^{\Re \mu_{z}^{-}x}\right]\\
&\leq \left(\sup_{x\leq 0} \varphi_{z}(x)\right)\frac{2}{\min \left\{ \Re \mu_{z}^{-},\Re \mu_{z}^{+}\right\}}.
\end{align*}
A similar estimate is obtained for $x\geq 0$, where now the supremum is taken for $x\geq 0$.
From this, \eqref{Schur cond} follows.
\end{proof}


\subsection{The spectrum of $\mathscr{L}_m$}\label{Subsec Essential Spectrum}
According to Proposition~\ref{Prop Resolvent}, 
\[
	\operatorname{Spec}(\mathscr{L}_{m}) \subseteq S_{m}.
\]
Since $\mathscr{L}_{m}$ is $\mathcal{T}$-selfadjoint, the residual spectrum is empty \cite[Sec. 5.2.5.4]{Krejcirik-Siegl15} and the first four essential spectra coincide \cite[Thm. IX.1.6(ii)]{Edmunds-Evans18}. In this subsection we show that the spectrum coincides exactly with $S_m$ and find $\operatorname{Spec}_{\mathrm{e5}}(\mathscr{L}_{m})$. We treat the case $m=0$ separately.

\subsubsection*{Case $m>0$} Let $z^{\pm}\coloneqq \tau\pm i$ for fixed $\tau \in \mathbb{R}$ such that $|\tau| \geq m$, so $z^{\pm} \in S_{m}$. Our goal is to construct singular Weyl sequences $(\Psi_n^\pm)_{ n\in\N}$ for $z^{\pm}$, to conclude that they are in the appropriate part of the spectrum.
With this purpose in mind, let 
	\begin{equation} \label{solution for Weyl}
		u(x) \coloneqq   \begin{pmatrix}
		\sqrt{|m+\tau|} \\
		 i \,\mathrm{sgn}(\tau)\sqrt{|m-\tau|}
	\end{pmatrix} e^{ i \sqrt{\tau^2-m^2} x},
	\end{equation}
which is one of the solutions of 
\begin{equation}\label{D0-tau}
\left(\mathcal{D}_{m}-\tau I\right)u = 0,
\end{equation}
for the free Dirac operator \eqref{Free Dirac}.
Let
\begin{equation}\label{Psi n}
\Psi_{n}^{\pm}(x) \coloneqq  c_{n}\varphi_n(\pm x)u(x),
\end{equation}
where 	\[
		c_n=\frac{1}{\sqrt{2| \tau |\left(n+\frac{2}{3}\right)}}
	\] 
and $\varphi_{n}$ is the continuous compact support cut-off function given by
	\[
	\varphi_n(x) \coloneqq \begin{cases}
		 x-n & \text{for } x \in [n, n+1], \\
		1& \text{for } x \in [n+1, 2n+1], \\
		2n+2-x & \text{for } x \in [2n+1, 2n+2],\\
		0 & \text{otherwise}.
	\end{cases}
	\]
Then, $\|\Psi_n^{\pm}\|=1$.

By considering the location of the supports and from \eqref{D0-tau}, it follows that
	\begin{align*}
\left( \mathscr{L}_{m} -z^{\pm} I\right)\Psi_{n}^{\pm}(x) 
	&= \left( \mathcal{D}_{m} -\tau I\right)\Psi_{n}^{\pm}(x)\\
	&=c_{n}\varphi_{n}(\pm x) \left( \mathcal{D}_{m} -\tau I\right)u(x) + c_{n}\left[ \mathcal{D}_{m}-\tau I, \varphi_{n}(\pm x) I \right]u(x)\\
	&=\pm c_{n}  \varphi_{n}'(\pm x) J u(x).
	\end{align*}
Then
\[ \left\Vert \left( \mathscr{L}_{m} -z^{\pm} I \right)\Psi_{n}^{\pm}(x) \right\Vert = \Vert c_{n}  \varphi_{n}'(\pm x) Ju(x) \Vert=\sqrt{\frac{2}{n+\frac{2}{3}}} \to 0\]
as $n\to \infty$.
Therefore, indeed, $(\Psi_n^\pm)_{ n\in\N}$ is a Weyl sequence for $z^{\pm}$. Thus, $z^{\pm} \in \operatorname{Spec}(\mathscr{L}_{m})$. This completes the proof that $\operatorname{Spec}(\mathscr{L}_{m})=S_{m}$ for all $m>0$, which we will use without further explicit mention. 

\medskip


Since the support of the Weyl sequence moves to infinity, $S_{m} \subset \operatorname{Spec}_{\mathrm{e2}}(\mathscr{L}_{m})$. Consequently, $\operatorname{Spec}_{\mathrm{ej}}(\mathscr{L}_{m})=S_{m}$ for $\mathrm{j}\in [[1,4]]$. To further deduce that
\[
\operatorname{Spec}_{\mathrm{e5}}(\mathscr{L}_{m})=S_{m},
\] it suffices to note that $\C\setminus \operatorname{Spec}_{\mathrm{e1}}(\mathscr{L}_{m})= \rho(\mathscr{L}_{m})$ is connected \cite[Prop. 5.4.4]{Krejcirik-Siegl15}. 

Finally, observe that
$\operatorname{Spec}_{\mathrm{p}}(\mathscr{L}_{m})=\varnothing.$
Indeed, the expression of the fundamental matrices \eqref{Der Phi} is also valid for $z\in S_m$. Applying \eqref{u pm} with $g=0$ and noting that $\Re \mu_{z}^{\pm}=0$ for $z\in \{\tau \pm i: |\tau|\geq m \}$, we conclude that there are no nontrivial solutions to the eigenvalue equation $(\mathscr{L}_{m}-z)u=0$ that belong to $L^{2}(\R,\C^2)$. Therefore, no point in $S_{m}$ can be an eigenvalue. This completes the description of the spectrum of $\mathscr{L}_m$ for $m>0$.

\subsubsection*{Case $m=0$} For $\mathscr{L}_0$ and $z^{\pm}=\tau \pm i$ where $\tau\not= 0$, the previous argumentation involving $(\Psi_{n}^{\pm})_{n\in \N}$ is still applicable, so we know that there exist singular Weyl sequences for $z^{\pm}\not=\pm i$. On the other hand, for $\tau=0$, we can choose constant $u(x)=\begin{pmatrix}
1\\0\end{pmatrix}$ instead of \eqref{solution for Weyl}, which is a solution to the eigenvalue equation \eqref{D0-tau}, and run a similar proof as above. Thus, we know up front that
\begin{equation}\label{Ess Spec m=0}
 \partial S_0=\{ z\in \C: |\Im z|=1\} \subset \operatorname{Spec}_{\mathrm{e2}}(\mathscr{L}_{0})=\operatorname{Spec}_{\mathrm{ej}}(\mathscr{L}_{0})
\end{equation}
for $\mathrm{j}=1,3,4$.

Our focus now is to show that all the points $z$ in the interior of the strip $S_{0}$ are eigenvalues with $\mathscr{L}_{0}-z$ being a Fredholm operator of index zero.

\begin{lemma}\label{Lem Point Spectrum}
For $m=0$, let $\mathscr{L}_{0}$ be defined as in \eqref{Dirac Operator}. Let $z\in \C$ be such that $|\Im z|<1$. Then the following holds true.
\begin{enumerate}[a)]
\item \label{lppa} $\Ker(\mathscr{L}_{0}-z)=\mathrm{span} \{v_{z}\} $ where
\[ v_{z}(x)\coloneqq \left\{\begin{aligned}
&\begin{pmatrix}
1\\
-i
\end{pmatrix}e^{(1-iz)x}\qquad &&\text{for }x\leq 0,\\
 &\begin{pmatrix}
1\\
-i
\end{pmatrix}e^{-(1+iz)x} \qquad &&\text{for }x\geq 0.
\end{aligned} 
\right.\]
\item \label{Kernel D*}$\Ker(\mathscr{L}_{0}^{*}-z)=\mathrm{span} \{\widetilde{v}_{z}\} $ where
\[ \widetilde{v}_{z}(x)\coloneqq \left\{\begin{aligned}
&\begin{pmatrix}
1\\
i
\end{pmatrix}e^{(1+iz)x}\qquad &&\text{for }x\leq 0,\\
 &\begin{pmatrix}
1\\
i
\end{pmatrix}e^{-(1-iz)x} \qquad &&\text{for }x\geq 0.
\end{aligned} 
\right.\]
\item \label{lppc} $\mathrm{Ran}(\mathscr{L}_{0}-z)$ is closed.
\item \label{algmult} The algebraic multiplicity of $z$ is infinite.
\end{enumerate}
\end{lemma}

\begin{proof}
Consider once again the equation $(\mathscr{L}_{0}-z)u=f$ as in Subsection~\ref{Subsec Resolvent}. When $m=0$ and $|\Im z|<1$, the quantities $\mu_{z}^{\pm}$, $w_{z}^{\pm}$  in \eqref{Notation pm} reduce to
\begin{equation}\label{Estimate m=0}
\begin{aligned}
&\mu_{z}^{+}=1+iz,\qquad &&\mu_{z}^{-}=1-iz, \qquad w_{z}^{+}=i, \qquad &&w_{z}^{-}=-i.
\end{aligned}
\end{equation}

We start with the proof of \ref{lppa}. Let $u\in H^{1}(\R,\C^2)$ be a solution to this eigenvalue equation.  We know that the restrictions of $u$ to $\R_{\pm}$, match \eqref{u pm} with $g=0$, namely
\[u_{\pm}(x)=   e^{\mu_{z}^{\pm}x}S^{\pm}_{z}\begin{pmatrix}
\alpha\\
\beta
\end{pmatrix}+ e^{-\mu_{z}^{\pm}x}T_{z}^{\pm}\begin{pmatrix}
\alpha\\
\beta
\end{pmatrix},\]
with $(\alpha,\beta)\in \C^2$ and
\begin{equation}\label{Matrix ST}
 S^{\pm}_{z}= \frac{1}{2}\begin{pmatrix}
1 && \mp i\\
\pm i && 1
\end{pmatrix},\qquad T_{z}^{\pm}=\frac{1}{2}\begin{pmatrix}
1 && \pm i\\
\mp i && 1
\end{pmatrix}.
\end{equation}
Since $\Re \mu_{z}^{\pm}>0$, the decay of $u$ at $\pm \infty$ should match the conditions \eqref{ab1} and \eqref{ab2} for $\alpha$ and $\beta$. Therefore, $\alpha=i\beta$. Hence \[u_{+}(x)= \alpha e^{-\mu_{z}^{+}x} \begin{pmatrix}
1\\
-i
\end{pmatrix} \qquad \text{and}\qquad u_{-}(x)= \alpha e^{\mu_{z}^{-}x} \begin{pmatrix}
1\\
-i
\end{pmatrix}.\]
Thus, indeed \ref{lppa} is valid. The proof of \ref{Kernel D*} is identical, so we omit it. 

To confirm \ref{lppc}, we will show that $\overline{\mathrm{Ran}(\mathscr{L}_{0}-z)}\subseteq \mathrm{Ran}(\mathscr{L}_{0}-z)$. Take $f\in \overline{\mathrm{Ran}(\mathscr{L}_{0}-z)}$ and seek for $u\in H^{1}(\R,\C^2)$ such that $(\mathscr{L}_{0}-z)u=f$. Now $u_{\pm}$ are given in the form \eqref{u pm} with $g=-Jf$. Since
\[ \overline{\operatorname{Ran}(\mathscr{L}_{0}-z)}=\left(\operatorname{Ker}(\mathscr{L}_{0}^{*}-\overline{z})\right)^{\perp},\]
according to \ref{Kernel D*}, it follows that
\begin{equation}\label{Cond Surjective}
\int_{-\infty}^{0}(f_{1}(y)-if_{2}(y)) e^{(1-iz)y}\, \dd y+ \int_{0}^{\infty} (f_{1}(y)-if_{2}(y)) e^{-(1+iz)y}\, \dd y=0.
\end{equation}
By proceeding as in the proof of Proposition~\ref{Prop Resolvent}, this is equivalent to the fact that both right-hand sides of the decaying conditions  \eqref{ab1} and \eqref{ab2} are equal. That is,
\begin{align*}
 -\int_{0}^{\infty} \left(g_{1}(y)-ig_{2}(y)\right)e^{-\mu_{z}^{+}y}\, \dd y=\alpha-i\beta = \int_{-\infty}^{0} \left(g_{1}(y)-ig_{2}(y)\right)e^{\mu_{z}^{-}y}\, \dd y.
\end{align*}
Set $\beta=0$. Then
\[ \begin{pmatrix}
\alpha\\
\beta
\end{pmatrix}= \begin{pmatrix}
-1 & i\\
0 & 0
\end{pmatrix} \int_{0}^{\infty} g(y)e^{-\mu_{z}^{+}y}\, \dd y = \begin{pmatrix}
1 & -i\\
0 & 0
\end{pmatrix} \int_{-\infty}^{0} g(y)e^{\mu_{z}^{-}y}\, \dd y .\]
Replacing these constants in \eqref{u pm} and writing the expression back in terms of $f_1$ and $f_2$, using $g=-Jf$, gives 
\begin{align*}
u_{+}(x)=&C_{1}^{+}\int_{x}^{+\infty} e^{\mu_{z}^{+}(x-y)}f(y)\, \dd y+C_{2}^{+}\int_{0}^{x} e^{-\mu_{z}^{+}(x-y)}f(y)\, \dd y\\
&+C_{3}^{+}\int_{0}^{\infty} e^{-\mu_{z}^{+}(x+y)}f(y)\, \dd y,
\end{align*}
and
\begin{align*}
u_{-}(x)=&C_{1}^{-}\int_{-\infty}^{x} e^{-\mu_{z}^{-}(x-y)}f(y)\, \dd y+C_{2}^{-}\int_{x}^{0} e^{\mu_{z}^{-}(x-y)}f(y)\, \dd y\\
&+C_{3}^{-}\int_{-\infty}^{0} e^{\mu_{z}^{-}(x+y)}f(y)\, \dd y,
\end{align*}
for suitable constant matrices $C_{j}^{\pm}$. Now, by invoking the Schur Test in identical manner as in the proof of Proposition~\ref{Prop Resolvent} - we omit the details, we find that $u\in L^{2}(\R,\C^2)$ and thus $u\in H^{1}(\R,\C^2)$. Hence, indeed $f\in \mathrm{Ran}(\mathscr{L}_{m}-z)$. This completes the proof of \ref{lppc}.

For the final statement, we invoke \cite[Thm. 2.16 in Ch.~2]{Locker_2000}, which establishes that only one of the following cases can occur. Either
\begin{equation*}
\Sigma = \big\{w\in\C\,:\, |\Im w|<1 \text{ and the algebraic multiplicity of }w \text{ is finite}  \big\}
\end{equation*}
or 
\begin{align*}
\Sigma = \big\{w\in\C\,:\, |\Im w|<1 \text{ and the algebraic multiplicity of }w \text{ is infinite} \big\}.
\end{align*}
From the above, the geometric multiplicity of all $z\in\Sigma$ is one. Therefore no point in $\Sigma$ has algebraic multiplicity zero. Thus, according to \cite[Thm. 2.15 in Ch.~2]{Locker_2000}, the set on the right-hand side of the first possibility is countable. But this is not the case for $\Sigma$, hence, the second possibility is the one that holds. This confirms \ref{algmult} and completes the proof of the lemma.
\end{proof}

According to Lemma \ref{Lem Point Spectrum}, we now know that 
\begin{equation}\label{Point Spec m=0}
\operatorname{int}S_0=\{ z\in \C: |\Im z|<1\} \subseteq \operatorname{Spec}_{\mathrm{p}}(\mathscr{L}_{0})\setminus \operatorname{Spec}_{\mathrm{e3}}(\mathscr{L}_{0}). 
\end{equation}
Combining this with \eqref{Ess Spec m=0}, we conclude that indeed
\[\operatorname{Spec}(\mathscr{L}_{0})=S_{0} \qquad \text{ and }\qquad \operatorname{Spec}_{\mathrm{ej}}(\mathscr{L}_{0})=\operatorname{Spec}_{\mathrm{e2}}(\mathscr{L}_{0})=\{ z\in \C: |\Im z|=1\}\]
for $\mathrm{j}=1,3,4$.
Since no eigenvalue on the left-hand side of \eqref{Point Spec m=0} is isolated, and $S_0$ is the closure of its interior, then $\operatorname{Spec}_{\mathrm{e5}}(\mathscr{L}_{0})=S_{0}$.
Finally, since no point on the lines $\{ \tau \pm i: \tau \in \R\}$ is an eigenvalue, these lines exactly form the continuous spectrum of $\mathscr{L}_{0}$. In summary, 
\[
     \operatorname{Spec}_{\mathrm{p}}(\mathscr{L}_0)=\{ z\in \C: |\Im z|<1\}.
\]
We have therefore completed the analysis of the spectrum of $\mathscr{L}_m$ for all $m\geq 0$.


\section{The resolvent norm of $\mathscr{L}_{m}$ inside the numerical range}\label{Sec Resolvent Estimate}
Let $m>0$. In the second part of this section we give the proof of the sharp estimate for the resolvent norm of $\mathscr{L}_m$ claimed in Theorem~\ref{Theo Resolvent 2}. Before that, and in order to visualise the asymptotic behaviour of the different parameters involved in this proof, we establish a sub-optimal version of the statement. Later on, this preliminary estimate will provide a link with the analysis of the pseudospectrum of the perturbations of $\mathscr{L}_m$. 

From now on, we write $\tau=\Re z$ and $\delta=\Im z$, often without explicit mention. 

\subsection{Preliminary estimate }\label{Subsec Rough Estimate}

A direct application of the Schur Test gives the next preliminary version of Theorem~\ref{Theo Resolvent 2}. 

\begin{theorem}\label{Theo resolvent}
For $m>0$, let $\mathscr{L}_{m}$ be given by \eqref{Dirac Operator}. Let $z=\tau+i \delta$ for $|\delta|<1$. Then,
\begin{equation}\label{Res est 2}
\begin{gathered} \big\|(\mathscr{L}_{m}-z)^{-1}\big\|\geq \frac{1}{2\sqrt{2}m\sqrt{1-\delta^2}}\left[\tau^2+\frac{m}{2}\tau+\frac{8(1+\delta^2)+4m(1+\delta)-m^2}{8}+\mathcal{O}(\tau^{-1})\right] \\ \text{and} \\ \big\|(\mathscr{L}_{m}-z)^{-1}\big\|\leq \frac{4}{m(1-\delta^2)}\left[\tau^2+\frac{1+\delta^2}{4}+m+\mathcal{O}(\tau^{-2})\right] \end{gathered}
\end{equation}
as $\vert \tau \vert \to +\infty$. Both limits inside the brackets converge uniformly\footnote{Here and elsewhere below, we adopt the common say that two function are related by $f(z)=g(z)+\mathcal{O}(\tau^{-p})$ as $|\tau|\to\infty$ and that the limit converges uniformly for parameters satisfying a condition, iff there exist a constant $K>0$ independent of $\tau$ and these parameters, such that $|f(z)-g(z)|\leq K|\tau|^{-p}$ for all $|\tau|\geq 1$ and all such parameters satisfying the condition.} for all $|\delta|<1$ and $m$ on compact set.
\end{theorem}

According to this statement, the resolvent norm of $\mathscr{L}_m$ has leading asymptotic behaviour $|\Re z|^2$, as $z\to\infty$ inside the band forming the numerical range. If we assume that Theorem~\ref{Theo Resolvent 2} has already been proved, then we see that the leading order coefficient of the upper bound in \eqref{Res est 2} is optimal, up to a factor~4. However, the leading coefficient of the lower bound is sub-optimal. The latter is a consequence of choosing below a sub-optimal pseudo-mode. We have included this calculation, in order to illustrate the construction of the optimal pseudo-mode for the full proof of the main Theorem~\ref{Theo Resolvent 2} in the next subsection, without too many technical details. 

Before proceeding to the asymptotic analysis of the coefficients involved, we highlight the connection with recent results of the same nature. 

According to \cite[Thm. 2.2]{Henry-Krejcirik17}, inside the numerical range the Schr{\"o}dinger operator with a dislocated potential, $-\frac{\dd^2}{\dd x^2}+i\sgn(x)$, has a resolvent norm with linear leading order  of growth at infinity in $|\Re z|$. This discrepancy is explained by the fact that $\mathscr{L}_m$ is an operator of order 1, therefore the perturbation affects the resolvent by a higher order of magnitude.  
  
When the potential is smooth and asymptotic to a dislocation, the effect of the perturbation on the Dirac operator appears to become substantially stronger. For instance, in the case of the operator $\widetilde{\mathscr{L}}=\mathcal{D}_{m}+i \frac{2}{\pi} \mathrm{arctan}(x) I $, by using the recent results of \cite[Thm. 3.11]{Krejcirik-Nguyen22}, it might be possible to show that there exists a WKB-pseudomode $\Psi_{\tau,n}$, such that \[\frac{\left\Vert\left(\widetilde{\mathscr{L}} -\tau I \right)\Psi_{\tau,n}\right\Vert}{\left\Vert \Psi_{\tau,n}\right\Vert}=\mathcal{O}(\tau^{-n})\] as $|\tau|\to \infty$ for all $n\in\mathbb{N}$. This maybe closely linked to the fact that interior singularities of a differential expression would create an exponential growth of the resolvent norm away from the spectrum. The phenomenon could be similar to the one reported recently in \cite{BM24} for the case of a singular Sturm-Liouville operator and is worthy of further exploration.

\medskip

We now give the proof of Theorem~\ref{Theo resolvent}. We split the proof into 3 steps.

\subsubsection*{Step 1} The first step is an asymptotic analysis of the parameters involved in the resolvent estimate. This result is useful on its own, and will be invoked repeatedly in later sections.

\begin{lemma}\label{Lem Asymptotic}
For $m>0$, let $\mu_{z}^{\pm}$ and $w_{z}^{\pm}$ be as in \eqref{Notation pm}, and let $k_z$ be as in \eqref{Function k}. For indices $p\in[[-4,2]]$ and parameters \[\mathsf{P}\in\left\{\mu_z^{\pm},\,\frac{\Re \mu^{\pm}_z}{1\mp \delta},\,w_z^{\pm},\,|w_z^{\pm}|,\,k_z\right\},\] let the coefficients $C_p(\mathsf{P})\in \mathbb{C}$ be prescribed in Table~\ref{Table coef}. There exists a constant $K_m>0$, independent of $\tau$ and $\delta$, such that
\[ \left| \mathsf{P}-\sum_{p=-4}^2C_p(\mathsf{P})\tau^{p}\right|\leq \frac{K_m}{\tau^5} \] 
for all $|\tau|\geq 1$ and $|\delta|<1$. This constant can be replaced by a uniform constant for all $m$ on a compact set. 
\end{lemma}

\begin{sidewaystable}
\centering
\rule{0ex}{120ex}
\begin{tabular}{c|rllllll}
$\mathsf{P}$ & 2&1&0&-1&-2&-3&-4 \\ \hline 
$\mu^{\pm}_z$ & 0 & $\pm i$ & $1\mp \delta$& $\mp\frac{im^2}{2}$ & $\frac{(1\mp \delta)m^2}{2}$ & $\frac{\mp im^2}{8}[m^2-4(1\mp\delta)^2]$ & $  \mp\frac{m^2}{8}(1\mp\delta)[3m^2-4(1\mp\delta)^2]$\\
$w^{\pm}_z$ & 0 &0& $\pm i$&$\mp im$ & $\frac{m}{2}[2(1\mp \delta)\pm im]$ & $\frac{m}{2}[im\mp (1+i)(1\mp\delta)]\times$&$\frac{m}{8}\left\{4(1\mp\delta)[3m^2+3i(1\mp \delta)m \right.$  \\ &&&&&& \qquad$[(1+i)(1\mp\delta)\mp m]$& $\qquad \left. -2(1\mp\delta)^2]\pm 3im^3\right\}$\\ 
$|w^{\pm}_z|$ & 0 &0& 1& $-m$ & $\frac{m^2}{2}$ & $\frac{m}{2}[2(1\mp \delta)^2-m^2]$ & $m^2\left[\frac38 m^2-(1\mp\delta)^2\right]$\\
$k_z$ & $\frac{i}{m}$ &$-\frac{2\delta}{m}$& $i\left(\frac{\delta^2-1}{m}-m\right)$ & 0 & $im$&$2\delta m$ & $im(1+3\delta^2-m^2)$ \\
$\frac{\Re \mu^{\pm}_z}{1\mp \delta}$& 0 & 0 & 1& 0 & $\frac{m^2}{2}$ & 0 & $  \mp\frac{m^2}{8}[3m^2-4(1\mp\delta)^2]$\\
\end{tabular}
\caption{Coefficients $C_p(\mathsf{P})$ for the parameters in Lemma~\ref{Lem Asymptotic}. \label{Table coef}}
\end{sidewaystable}

\begin{proof}
We only provide a proof sketch, because the details are lengthy routine calculations. 

The strategy for each of the parameters $\mathsf{P}$, is to find a suitable Taylor expansion from which the coefficients $C_p\in \mathbb{C}$ are  obtained from those of the different powers of $\frac{1}{\tau}$. The remainder at power~$5$ in this expansion, gives the uniform bound in terms of $\delta$ and $m$. The specific Taylor series and coefficients are determined from the following reductions,      
\begin{align*}
\mu_{z}^{\pm}&= i\tau \sqrt{1-\frac{m^2+(1\mp\delta)^2}{\tau^2}\mp i\frac{2(1\mp\delta)}{\tau}}\qquad \text{and} \\
w_{z}^{\pm} &= \sqrt{\frac{1-\frac{m\pm i(1\mp \delta)}{\tau}}{1+\frac{m\mp i(1\mp \delta)}{\tau}}} .
\end{align*}
For $\mathsf{P}=|w_z^{\pm}|$, we use \[|w_z^{\pm}|=\sqrt{w_z^{\pm}w_{\overline{z}}^{\mp}}=\sqrt{1+g(\tau)}\] where $g(\tau)=\mathcal{O}\big(\frac{1}{\tau}\big)$ as $|\tau|\to \infty$.  For $\mathsf{P}=k_z$, we use \eqref{Function k} and the expansions of $w_{\pm}$.

Concretely, for $\mathsf{P}=\mu_z^{+}$, we write  $\mu_{z}^{+}=i\tau f\left(\frac{1}{\tau}\right)$ where $f$ is analytic near $0$. Expanding in power series at that point, gives \[ i\frac{f(x)}{x}=\sum_{k=-1}^{4}C_{-k}(\mu_z^+)x^p+\frac{x^5}{6!} \int_{0}^{1}(1-s)^5f^{(6)}(xs)\, \dd s,\] where the coefficients $C_p(\mu_z^+)$ coincide with those in Table~\ref{Table coef}. Moreover, $f^{(6)}$ has a denominator involving a term of the form $\big(a(\delta,m)x^2+b(\delta,m)x+1\big)^{\frac{11}{2}}$ where $a$ and $b$ are bounded uniformly in both variables. Therefore, taking $|x|$ small enough,  ensures that $|f^{(6)}(x)|$ is bounded uniformly in $x$, $|\delta|<1$ and $m>0$ in compact sets.

For the expansion of $\frac{\Re \mu_{z}^{\pm}}{1\mp \delta}$, recall that $\Re \sqrt{\gamma}=\frac{|\Im \gamma|}{\sqrt{2(|\gamma|-\Re \gamma)}}$ for $\gamma\in \C\setminus \R$, and so
\[ \frac{\Re \mu_{z}^{\pm}}{1\mp \delta} = \sqrt{\frac{2}{1-\frac{m^2+(1\mp\delta)^2}{\tau^2}+\sqrt{\left(1-\frac{m^2+(1\mp\delta)^2}{\tau^2}\right)^2+\frac{4(1\mp \delta)^2}{\tau^2}}}}.\]
\end{proof}

\subsubsection*{Step 2} Consider the upper bound first. According to Lemma~\ref{Lem Upper bound}, we have that
\begin{align*}
\Vert \left(\mathscr{L}_{m}-z\right)^{-1} \Vert
&\leq  \frac{ \left( |k_{z}|+1\right)\max \left\{ |w_{z}^{\pm}|+\frac{1}{|w_{z}^{\pm}|},|w_{z}^{+}|+\frac{1}{|w_{z}^{+}|}\right\}}{\min \left\{\Re \mu_{z}^{-}, \Re \mu_{z}^{+} \right\}}\\
&\leq \left( |k_{z}|+1\right)\left(|w_{z}^{+}|+\frac{1}{|w_{z}^{+}|}+|w_{z}^{-}|+\frac{1}{|w_{z}^{-}|}\right)\left(\frac{1}{\Re \mu_{z}^{-}}+\frac{1}{\Re \mu_{z}^{+}}\right),
\end{align*}
for all $z\in \rho(\mathscr{L}_{m})$.
By Lemma~\ref{Lem Asymptotic}, the asymptotic behaviour of each of the multiplying terms above is
\begin{align}
 | k_{z} |+1 &= \frac{1}{m}\tau^2+\left(\frac{1+\delta^2}{m}-m+1\right)+\mathcal{O}(\tau^{-4}),\notag\\
|w_{z}^{+}|+\frac{1}{|w_{z}^{+}|}+|w_{z}^{-}|+\frac{1}{|w_{z}^{-}|}&= 2+m^2\tau^{-2}+\mathcal{O}(\tau^{-4}),\label{eqforR2}\\
 \frac{1}{\Re \mu_{z}^{-}}+\frac{1}{\Re \mu_{z}^{+}}&=\frac{1}{1-\delta^2}\left[2-m^2\tau^{-2}+\mathcal{O}(\tau^{-4})\right],\label{eq for R2}\end{align}
as $|\tau|\to\infty$. Note that the first two coefficients of the terms involving the maximum  above are identical. Multiplying and collecting terms, gives the upper bound in \eqref{Res est 2} as required. Observe that the fact that the limit inside the bracket is uniform in the parameters, carries over from the conclusion of  Lemma~\ref{Lem Asymptotic}.

\subsubsection*{Step 3} Now consider the lower estimate of the resolvent in \eqref{Res est 2}. Set a pseudo-mode
\[
		\psi_0(x) \coloneqq  \begin{pmatrix}
			1 \\ w_{\overline{z}}^{+} 
		\end{pmatrix} e^{\mu^+_{\overline{z}} x} \chi_{(-\infty,0)} (x), 
	\]
which belongs to $L^{2}(\R,\C^2)$ and has a norm $\Vert \psi_{0} \Vert = \frac{\sqrt{1+\vert w_{z}^{-} \vert^2}}{\sqrt{2\Re \mu_{z}^{-}}}$.
From the formula of the kernel in \eqref{G} and considering the location of the support of $\psi_{0}$, we have
\begin{align*}
\|(\mathscr{L}_{m}-z)^{-1}\psi_0 \|^2\geq &\int_{0}^{+\infty}\left\vert\int_{-\infty}^0 N_{3,z}\begin{pmatrix}
			1 \\ w_{\overline{z}}^{+} 
		\end{pmatrix}e^{-\mu_{z}^{+}x+2\Re\mu_{z}^{-}y} \, \dd y\right\vert_{\C^2}^2 \, \dd x\\
		=& \left\vert N_{3,z} \begin{pmatrix}
			1 \\ w_{\overline{z}}^{+} 
		\end{pmatrix}\right\vert_{\C^2}^2 \frac{1}{8 \left(\Re \mu_{z}^{-}\right)^2 \Re \mu_{z}^{+}}\\
=& \left(\frac{(1+\vert w_{z}^{-} \vert^2)\sqrt{1+\vert w_{z}^{+} \vert^2}}{\vert w_{z}^{+}+w_{z}^{-}\vert}\right)^2\frac{1}{8 \left(\Re \mu_{z}^{-}\right)^2 \Re \mu_{z}^{+}}.
\end{align*}
Then,\begin{align*}
\Vert (\mathscr{L}_{m}-z)^{-1} \Vert & \geq  \frac{\|(\mathscr{L}_{m}-z)^{-1}\psi_0 \|}{\Vert \psi_{0} \Vert}\\
& = \frac12 \left[\frac{(1+\vert w_{z}^{-} \vert^2)(1+\vert w_{z}^{+} \vert^2)}{\Re (\mu_{z}^{-})\Re (\mu_{z}^{+})}\right]^{\frac12}\ \vert w_{z}^{+}+w_{z}^{-}\vert^{-1}.\end{align*}
According to Lemma~\ref{Lem Asymptotic}, the first bracket is $\mathcal{O}(1)$ and the second term  is $\mathcal{O}(\tau^2)$. Multiplying and collecting coefficients, gives the lower estimate of the resolvent in \eqref{Res est 2}. We omit further details, but this confirms the validity of Theorem~\ref{Theo resolvent}.


\subsection{Proof of Theorem \ref{Theo Resolvent 2}}\label{Subsec Sharp Estimate}
We now establish improved estimates, which enable identifying the sharp constants for the resolvent norm asymptotic of $\mathscr{L}_m$ in the instability band. These will render  Theorem~\ref{Theo Resolvent 2} as a corollary. The strategy for the computation of the different constants has two parts. We first split the resolvent as the sum of two integral operators,
\[ (\mathscr{L}_{m}-zI)^{-1}= T_{1}(z)+T_{2}(z),\]
where the norm of the operator $T_{1}(z)$ carries the leading order behaviour and $\|T_{2}(z)\|$ is $\mathcal{O}(1)$ at infinity. We then proceed as in the proof of Theorem~\ref{Theo resolvent}, and compute the asymptotic coefficients of $\|T_{1}(z)\|$ and $\|T_{2}(z)\|$. In the case of the lower bound, we find a pseudo-mode for $T_{1}(z)$ that is sharp for the leading order, by optimising the contribution of the two components of the wave function in an adaptive manner as $z$ moves towards infinity. For a further discussion on the computation of other coefficients see Remark~\ref{Rem:more terms} at the end of this section.

The explicit formulas of the two integral operators is the following. For $j=1,2$,
\begin{equation}\label{R12}
T_j(z)f(x)\coloneqq \int_{\R} \mathcal{R}_{j,z}(x,y) f(y)\, \dd y
\end{equation}
where,
\begin{equation}\label{R1}
\mathcal{R}_{1,z}(x,y)\coloneqq  \left\{
\begin{aligned}
 & N_{1,z}\,e^{\mu^-_z(x+y)} , && \text{for } \{y < x < 0\};\\ 
 & N_{3,z}\, e^{-\mu^+_zx+\mu^-_zy}, && \text{for } \{y < 0 < x\}; \\
 & N_{4,z}\,e^{-\mu^+_z(x+y)}  , && \text{for } \{0 < y < x\}; \\
 & N_{7,z}\,e^{-\mu^+_z(x+y)}, && \text{for } \{0 < x < y\}; \\
 & N_{8,z} \, e^{\mu^-_zx-\mu^+_zy} , && \text{for } \{x< 0< y\}; \\
 & N_{10,z}\,e^{\mu^-_z(x+y)} , && \text{for } \{x< y< 0\};
\end{aligned}
\right.
\end{equation}
and
\begin{equation}\label{R2}
\mathcal{R}_{2,z}(x,y)\coloneqq  \left\{
\begin{aligned}
&  N_{2,z}\, e^{-\mu^-_z(x-y)}, && \text{for } \{y< x< 0\};\\ 
&  0, && \text{for } \{y< 0 < x\}; \\
&  N_{5,z}\,e^{-\mu^+_z(x-y)} , && \text{for } \{0< y < x\}; \\
& N_{6,z}\,e^{\mu^+_z(x-y)} , && \text{for } \{0< x < y\}; \\
& 0 , && \text{for } \{x< 0< y\}; \\
& N_{9,z}\,e^{\mu^-_z(x-y)}  , && \text{for } \{x< y< 0\}.
\end{aligned}
\right.
\end{equation}
Recall the definition of the matrices $N_{k,z}$ at the beginning of Section~\ref{Sec Res and Sp}.

A crucial reason for the split of the resolvent in this fashion is to observe that the term $\frac{1}{w_{z}^{+}+w_{z}^{-}}$ is the only one responsible for the quadratic growth of the norm and that this term appears only in  $N_{k,z}$ for $k \in \{1,3,4,7,8,10\}$. Moreover, now the kernel $\mathcal{R}_{1,z}$ is separable, therefore a sharp upper bound on its norm is obtained by means of H\"{o}lder's Inequality.
  
For $j\in [[1,10]]$, we denote $n_{j}=n_{j}(z)\coloneqq  |N_{j,z}|_{\C^{2\times 2}}$. The fact that $\det N_{j,z}=0$ and substitution give
\begin{equation}\label{nj}
\begin{array}{ll}
n_1 = n_{10}= \frac{1}{2}\left|\frac{w^+_z-w^-_z}{w^+_z+w^-_z}\right|\left(|w^-_z|+\frac{1}{|w^-_z|}\right),  & n_{2}=n_{9} = \frac{1}{2}\left(|w^-_z|+\frac{1}{|w^-_z|}\right), \\
 n_{3}=n_{8} =\frac{\sqrt{(1+|w^+_z|^2) (1+|w^-_z|^2)}}{|w^+_z+w^-_z|}, & n_{4}=n_{7} = \frac{1}{2}\left|\frac{w^+_z-w^-_z}{w^+_z+w^-_z}\right|\left(|w^+_z|+\frac{1}{|w^+_z|}\right),\\
n_{5}=n_{6} = \frac{1}{2}\left( |w^+_z|+\frac{1}{|w^+_z|}\right). &
\end{array}
\end{equation}
Note that the asymptotic of all the terms in the formulas for $n_j$ are available from the calculations performed in the previous subsection.  Also note that the operator norms of $\mathcal{R}_{j,z}$ are symmetric with respect to the axis $y=x$, as was the case for $\mathcal{R}_z$.

The next proposition gives explicit bounds for the norm of the leading operator $T_1(z)$. The match of all terms in these bounds except $\widetilde{B}_z$ and $B_z$ is the main ingredient that enables the sharp constants below. For $z\not=0$, we write 
\[
\operatorname{csgn}(z)=\frac{z}{|z|}.
\] 

\begin{proposition}\label{Prop Norm R1}
Let $m>0$. Let the operator $T_1(z)$ be as in \eqref{R12}-\eqref{R1}. For all $z\in \rho(\mathscr{L}_{m})$, let
\begin{align*}
&A_{z}\coloneqq  \frac{n_{1}^2}{4(\Re \mu_{z}^{-})^2}, \qquad C_{z}=\frac{n_{4}^2}{4(\Re \mu_{z}^{+})^2},\qquad D_{z}\coloneqq  \frac{n_{3}^2}{4(\Re \mu_{z}^{-})(\Re \mu_{z}^{+})},\\
&B_{z}\coloneqq  \frac{n_{3}}{2\sqrt{(\Re \mu_{z}^{+})(\Re \mu_{z}^{-})}}\left(\frac{n_{1}}{\Re \mu_{z}^{-}} +\frac{n_{4}}{\Re \mu_{z}^{+}}\right),\\
&\widetilde{B}_{z}\coloneqq  \frac{n_{3}}{2\sqrt{(\Re \mu_{z}^{+})(\Re \mu_{z}^{-})}} \left(\frac{n_{1}\Re \left( \mathrm{csgn}\left(\frac{w_{z}^{+}-w_{z}^{-}}{w_{z}^{-}} \right)\right)}{\Re \mu_{z}^{-}} -\frac{n_{4}\Re \left( \mathrm{csgn}\left(\frac{w_{z}^{+}-w_{z}^{-}}{w_{z}^{+}} \right)\right)}{\Re \mu_{z}^{+}}\right).
\end{align*}
Then, $T_1(z)$ is bounded for all $z\in\rho(\mathscr{L}_m)$. Moreover,
\begin{equation} \label{eq bounds for R1}
\begin{gathered}
\Vert T_1(z) \Vert \leq \frac{1}{\sqrt{2}} \sqrt{\sqrt{(A_{z}-C_{z})^2+B_{z}^2}+A_{z}+C_{z}+2D_{z}}\\ \text{and} \\
\Vert T_1(z) \Vert \geq \frac{1}{\sqrt{2}} \sqrt{\sqrt{(A_{z}-C_{z})^2+\widetilde{B}_{z}^2}+A_{z}+C_{z}+2D_{z}}.
\end{gathered}
\end{equation}
\end{proposition}

\begin{proof}

\emph{Upper bound}. Let $\Psi\in L^{2}(\R,\C^2)$ be such that $\Vert \Psi \Vert^2=1$. Then,
\begin{equation}\label{Norm ineq}
\Vert T_1(z) \Psi \Vert^2 \leq \int_{\R} \left( \int_{\R} |\mathcal{R}_{1,z}(x,y)|_{\C^{2\times 2}} |\Psi(y)|_{\C^2}\, \dd y \right)^2\, \dd x.
\end{equation}
The right-hand side of this expression is given explicitly as follows,
\begin{align*}
&\int_{\R} \left( \int_{\R} |\mathcal{R}_{1,z}(x,y)|_{\C^{2\times 2}} |\Psi(y)|_{\C^2}\, \dd y \right)^2\, \dd x\\
&=\frac{1}{2(\Re \mu_{z}^{-})}\left( n_{1}\int_{-\infty}^{0} e^{(\Re \mu_{z}^{-})y}|\Psi(y)|_{\C^2}\, \dd y+n_{8}\int_{0}^{+\infty} e^{-(\Re \mu_{z}^{+})y}|\Psi(y)|_{\C^2}\, \dd y \right)^2\\
&+\frac{1}{2(\Re \mu_{z}^{+})}\left( n_{3}\int_{-\infty}^{0} e^{(\Re \mu_{z}^{-})y}|\Psi(y)|_{\C^2}\, \dd y+n_{4}\int_{0}^{+\infty} e^{-(\Re \mu_{z}^{+})y}|\Psi(y)|_{\C^2}\, \dd y \right)^2.
\end{align*}
By applying Holder's Inequality to the inner integrals, we get
\begin{align*}
\Vert T_1(z) \Psi\Vert^2&\leq \frac{1}{2\Re \mu_{z}^{-}}\left( \frac{n_{1}}{\sqrt{2\Re \mu_{z}^{-}}}\Vert \Psi\Vert_{L^2(\R_{-})}+\frac{n_{8}}{\sqrt{2\Re \mu_{z}^{+}}}\Vert \Psi\Vert_{L^2(\R_{+})}\right)^2\\
&\qquad +\frac{1}{2\Re \mu_{z}^{+}}\left( \frac{n_{3}}{\sqrt{2\Re \mu_{z}^{-}}}\Vert \Psi\Vert_{L^2(\R_{-})}+\frac{n_{4}}{\sqrt{2\Re \mu_{z}^{+}}}\Vert \Psi\Vert_{L^2(\R_{+})} \right)^2\\
&\leq \frac{n_{1}^2}{4(\Re \mu_{z}^{-})^2}\Vert \Psi\Vert_{L^2(\R_{-})}^2+\frac{n_{4}^2}{4(\Re \mu_{z}^{+})^2}\Vert \Psi\Vert_{L^2(\R_{+})}^2+\frac{n_{3}^2}{4(\Re \mu_{z}^{-})(\Re \mu_{z}^{+})}\\
&\qquad +\frac{n_{3}}{2\sqrt{(\Re \mu_{z}^{+})(\Re \mu_{z}^{-})}}\left(\frac{n_{1}}{\Re \mu_{z}^{-}} +\frac{n_{4}}{\Re \mu_{z}^{+}}\right)\Vert \Psi\Vert_{L^2(\R_{-})}\Vert \Psi\Vert_{L^2(\R_{+})}.
\end{align*}
Here we have used that $n_1=n_{10}$, $n_3=n_8$ and that $\|\Psi\|=1$.

Now, for constants $A,B,C\in \R$, 
\begin{equation} \label{max in circle}
\max_{x^2+y^2=1} Ax^2+Bxy+Cy^2=\frac{A+C+\sqrt{(A-C)^2+B^2}}{2}.
\end{equation}
Then, taking $x=\Vert \Psi\Vert_{L^2(\R_{-})}$ and $y=\Vert \Psi\Vert_{L^2(\R_{+})}$, we obtain that, for all $\Psi\in L^{2}(\R,\C^2)$ such that $\Vert \Psi \Vert^2=1$,
\begin{align*}
\Vert T_1(z) \Psi \Vert^2 \leq \frac{1}{2}\left(A_{z}+C_{z}+\sqrt{(A_{z}-C_{z})^2+B_{z}^2}+2D_{z}\right),
\end{align*}
where $A_{z},\,B_{z},\,C_{z},\,D_{z}$ are as above. Moreover, since all these 4 terms are finite for all $z\in\rho (\mathscr{L}_{m})$, then indeed $T_1(z)$ is a bounded operator.\\
\emph{Lower bound}. In order to obtain the lower bound estimate on the norm of $T_1(z)$, consider a pseudo-mode  of the form
\begin{equation}\label{Psi0}
\Psi_{0}(y)=
\begin{cases}
\beta\begin{pmatrix}
-1 \\ w_{\overline{z}}^{-} 
\end{pmatrix} e^{-\overline{\mu^+_z} y}\qquad &\text{if }  y\geq 0,\\ \\
\alpha\begin{pmatrix}
1 \\  w_{\overline{z}}^{+} 
\end{pmatrix} e^{\overline{\mu^-_z} y} \qquad &\text{if }  y\leq 0.\end{cases}
\end{equation}
where $\alpha, \beta \in \R$ will be chosen later to normalise $\Vert \Psi_{0} \Vert$ and maximise $\Vert T_1(z) \Psi_{0}\Vert$. 

We first compute the square of $\Vert T_1(z) \Psi_{0}\Vert$. Indeed,
\begin{align*}
\Vert T_1(z) \Psi_{0} \Vert^2&= \int_{-\infty}^{0}\left\vert \int_{-\infty}^{+\infty} \mathcal{R}_{1,z}(x,y) \Psi_{0}(y)\, \dd y\right\vert_{\C^{2}}^2\, \dd x+\int_{0}^{+\infty}\left\vert \int_{-\infty}^{+\infty} \mathcal{R}_{1,z}(x,y) \Psi_{0}(y)\, \dd y\right\vert_{\C^{2}}^2\, \dd x\\ &=:K_1+K_2.
\end{align*}
By the definition of $\Psi_{0}$, we have, 
\begin{align*}
&\int_{-\infty}^{+\infty} \mathcal{R}_{1,z}(x,y) \Psi_{0}(y)\, \dd y=e^{\mu_{z}^{-}x}\left[\frac{\alpha \frac{w_{z}^{+}-w_{z}^{-}}{w_{z}^{+}+w_{z}^{-}}(1+|w_{z}^{-}|^2)\frac{1}{w_{z}^{-}}}{4\Re \mu_{z}^{-}} +\frac{\frac{\beta(1+|w_{z}^{+}|^2)}{w_{z}^{+}+w_{z}^{-}}}{2\Re\mu_{z}^{+}} \right]\begin{pmatrix}
1\\
w_{z}^{-}
\end{pmatrix} \qquad \text{for } x<0,
\end{align*}
and 
\begin{align*}
&\int_{-\infty}^{+\infty} \mathcal{R}_{1,z}(x,y) \Psi_{0}(y)\, \dd y=e^{-\mu_{z}^{+}x}\left[\frac{\frac{\alpha(1+|w_{z}^{-}|^2)}{w_{z}^{+}+w_{z}^{-}}}{2\Re\mu_{z}^{-}}+\frac{\beta \frac{w_{z}^{+}-w_{z}^{-}}{w_{z}^{+}+w_{z}^{-}}(1+|w_{z}^{+}|^2)\frac{-1}{w_{z}^{+}}}{4\Re \mu_{z}^{+}}\right] \begin{pmatrix}
-1\\
w_{z}^{+}
\end{pmatrix} \qquad \text{for } x>0.
\end{align*}
Now,
\begin{align*}
K_{1}=&\frac{1}{2\Re\mu_{z}^{-}}\left[\frac{|\alpha|^2 n_{1}^2(1+|w_{z}^{-}|^2)}{4 (\Re \mu_{z}^{-})^2} +\frac{|\beta|^2n_{8}^2(1+|w_{z}^{+}|^2)}{4(\Re \mu_{z}^{+})^2}\right.\\
&\left.+\frac{n_{1}n_{8}\alpha \beta\sqrt{(1+|w_{z}^{-}|^2)(1+|w_{z}^{+}|^2)}}{2(\Re \mu_{z}^{+})(\Re \mu_{z}^{-})}\Re \left( \mathrm{csgn}\left(\frac{w_{z}^{+}-w_{z}^{-}}{w_{z}^{-}} \right)\right)\right],
\end{align*}
and
\begin{align*}
K_{2}=&\frac{1}{2\Re\mu_{z}^{+}}\left[\frac{|\alpha|^2 n_{3}^2(1+|w_{z}^{-}|^2)}{4 (\Re \mu_{z}^{-})^2} +\frac{|\beta|^2n_{4}^2(1+|w_{z}^{+}|^2)}{4(\Re \mu_{z}^{+})^2}\right.\\
&\left.- \frac{n_{3}n_{4}\alpha \beta \sqrt{(1+|w_{z}^{-}|^2)(1+|w_{z}^{+}|^2)}}{2(\Re \mu_{z}^{+})(\Re \mu_{z}^{-})}\Re \left(\mathrm{csgn}\left(\frac{w_{z}^{+}-w_{z}^{-}}{w_{z}^{+}} \right)\right)\right].
\end{align*}
Moreover,
\begin{align*}
\Vert \Psi_{0} \Vert^2= |\alpha|^2\frac{1+|w_{z}^{-}|^2}{2\Re \mu_{z}^{-}}+|\beta|^2\frac{1+|w_{z}^{+}|^2}{2\Re \mu_{z}^{+}}.
\end{align*}
To normalise $\Psi_{0}$ properly, we re-write $\alpha=\frac{\sqrt{2\Re \mu_{z}^{-}}}{\sqrt{1+|w_{z}^{-}|^2}}x$ and $\beta=\frac{\sqrt{2\Re \mu_{z}^{+}}}{\sqrt{1+|w_{z}^{+}|^2}}y$, for $x,y\in \R$ such that $x^2+y^2=1.$
We seek for pairs $(x,y)$, such that the expression
\begin{align*}
&\Vert T_1(z) \Psi_{0} \Vert^2\\
=&\frac{n_{1}^2}{4(\Re \mu_{z}^{-})^2}x^2+\frac{n_{4}^2}{4(\Re \mu_{z}^{+})^2}y^2+\frac{n_{3}^2}{4(\Re \mu_{z}^{-})(\Re \mu_{z}^{+})}\\
&+\frac{n_{3}}{2\sqrt{(\Re \mu_{z}^{+})(\Re \mu_{z}^{-})}}\left(\frac{n_{1}\Re \left( \mathrm{csgn}\left(\frac{w_{z}^{+}-w_{z}^{-}}{w_{z}^{-}} \right)\right)}{\Re \mu_{z}^{-}} -\frac{n_{4}\Re \left( \mathrm{csgn}\left(\frac{w_{z}^{+}-w_{z}^{-}}{w_{z}^{+}} \right)\right)}{\Re \mu_{z}^{+}}\right)xy
\end{align*}
attains its maximum. Therefore, according to \eqref{max in circle} once again, we have
\[ \Vert T_1(z) \Vert^2 \geq \Vert T_1(z) \Psi_{0} \Vert^2 =\frac{1}{2}\left(A_{z}+C_{z}+\sqrt{(A_{z}-C_{z})^2+\widetilde{B}_{z}^2}+2D_{z}\right).\]
Thus, the lower bound of $T_1(z)$ is indeed confirmed.
\end{proof}

Our next goal is to determine an upper bound for the norm of $T_2(z)$.

\begin{proposition}\label{Prop Norm R2}
Let $m>0$. Let the operator $T_2(z)$ be as in \eqref{R12}-\eqref{R2}. Then, $T_2(z)$ is bounded for all $z\in \rho(\mathscr{L}_{m})$. Moreover, 
\begin{equation}
\Vert T_2(z) \Vert\leq \max \left\{\frac{n_{2}}{\Re \mu_{z}^-},\frac{n_{5}}{\Re \mu_{z}^{+}}\right\}.
\end{equation}
\end{proposition}
\begin{proof}
We use the Schur Test to prove this proposition. For $x< 0$, we have 
\begin{align*}
\int_{\R} |\mathcal{R}_{2,z}(x,y)|_{\mathbb{C}^{2 \times 2}} \, \mathrm{d} y &=\int_{-\infty}^{x} n_{2} e^{-\Re \mu_{z}^{-}(x-y)}\, \dd y +\int_{x}^{0} n_{9} e^{\Re \mu_{z}^{-}(x-y)}\, \dd y \\
&=\frac{n_{2}}{\Re \mu_{z}^{-}}\left(2-e^{\Re \mu_{z}^{-}x}\right)\leq \frac{n_{2}}{\Re \mu_{z}^{-}}.
\end{align*}
Note that $n_{2}=n_{9}$.  Similarly, for $x\geq 0$, we have
\begin{align*}
\int_{\R} |\mathcal{R}_{2,z}(x,y)|_{\mathbb{C}^{2 \times 2}}\, \mathrm{d} y &=\int_{0}^{x} n_{5}e^{-\Re \mu_{z}^{+}(x-y)}\, \dd y +\int_{x}^{+\infty} n_{6} e^{\Re \mu_{z}^{+}(x-y)}\, \dd y \\
&=\frac{n_{5}}{\Re \mu_{z}^{-}}\left(2-e^{-\Re \mu_{z}^{+}x}\right)\leq \frac{n_{5}}{ \Re \mu_{z}^{+}}.
\end{align*}
Note that $n_{5}=n_{6}$.
Now recall that $|\mathcal{R}_{2,z}(x,y)|_{\C^{2\times 2}}=|\mathcal{R}_{2,z}(y,x)|_{\C^{2\times 2}}$. Hence, by virtue of the Schur Test, the claimed conclusion follows.
\end{proof}

\medskip

With the bounds for the norms of the components $T_j(z)$ at hand from the previous two propositions, we are now in the position to tackle Theorem~\ref{Theo Resolvent 2}. This will be a direct consequence of the following theorem, which gives slightly more refined asymptotic estimates. This theorem is one of the main contributions of this paper and it directly implies Theorem~\ref{Theo Resolvent 2}.

\begin{theorem} \label{Theo sharper estimate}
For $m>0$, let $\mathscr{L}_{m}$ be the dislocated Dirac operator \eqref{Dirac Operator}. Let $z=\tau+i \delta$ for $|\delta|<1$. Then, as $|\tau|\to \infty$,
\begin{equation}\label{Res sharper}
\begin{gathered}
	\big\|(\mathscr{L}_{m}-z)^{-1}\big\|\leq\frac{1}{m(1-\delta^2)}\left[\tau^2+\mathsf{P}^{+}_0+\mathsf{P}_{-2}^{+}\tau^{-2}+\mathcal{O}(\tau^{-4})\right] \\ \text{and} \\
\big\|(\mathscr{L}_{m}-z)^{-1}\big\|\geq \frac{1}{m(1-\delta^2)}\left[\tau^2+\mathsf{P}^{-}_0+\mathsf{P}^{-}_{-2}\tau^{-2}+\mathcal{O}(\tau^{-4})\right], 
\end{gathered}
\end{equation}
where
\[
\mathsf{P}_{0}^{\pm}=\frac14(3 + 2\delta^2 - 2m^2)\pm 2m, \quad \mathsf{P}^{+}_{-2}=\frac14(3m^2 + \delta^2(9m^2-8)) \quad \text{and} \quad \mathsf{P}^{-}_{-2}=\frac18(5m^2 + \delta^2(19m^2-16)).
\]
The limits inside the brackets converge uniformly for all $|\delta|<1$ and $m$ on compact sets.
\end{theorem}
The remainder of this subsection is devoted to sketching the proof of this theorem. That involves computing the coefficients in the asymptotic expansions of the various terms in the expressions \eqref{eq bounds for R1}. We only give the intermediate results in the computation of these coefficients. 

\subsubsection*{Coefficients for $(A_z-C_z)^2$}
We take the square of the following.
\[A_z-C_z=\frac{-\delta}{m^2(1-\delta^2)^2}\left[\sum_{k=-4}^4 \mathsf{P}_k \tau^k+\mathcal{O}(\tau^{-6})\right]\]
were\footnote{Here and henceforth, $\mathsf{P}_k=0$ for $k$ odd.}
\begin{gather*}
\mathsf{P}_4    =1,\qquad 
\mathsf{P}_2    =2\left(\delta^2-m^2+1\right),\qquad 
\mathsf{P}_0    =\delta^4+2 \delta^2 \left(m^2-1\right)+m^4+1, \\
\mathsf{P}_{-2} =-2\left(\delta^4+6 \delta^2+1\right) m^2 ,\qquad 
\mathsf{P}_{-4} =m^2 \left(2 \delta^6-5 \delta^4 \left(m^2-6\right)-\delta^2 \left(22m^2-30\right)-6 m^2+2\right).
\end{gather*}

\subsubsection*{Coefficients for $B_z^2$}
We use
\[
    B_z^2=D_z\left(\frac{n_3}{\operatorname{Re}\mu_z^-}+\frac{n_4}{\operatorname{Re}\mu_z^+}\right)^2
\]
together with the following. On the one hand,
\[D_z=\frac{1}{m^2(1-\delta^2)}\left[\sum_{k=-4}^4 \mathsf{P}_k \tau^k+\mathcal{O}(\tau^{-6})\right]\] where 
\begin{gather*}
\mathsf{P}_4    = \frac14 ,\qquad   
\mathsf{P}_2    = \frac{1}{2} \left(\delta^2-m^2+1\right) ,\qquad  
\mathsf{P}_0    = \frac{1}{4} \left(\delta^4+\delta^2 \left(m^2-2\right)+m^4+1\right) ,\\ 
\mathsf{P}_{-2} = -4 \delta^2 m^2 ,\qquad
	\mathsf{P}_{-4} = -\frac14\left(m^2 \left(\delta^4 \left(m^2-32\right)+\delta^2 \left(31 m^2-32\right)+m^2\right)\right).
\end{gather*} On the other hand, 
\[\left(\frac{n_3}{\operatorname{Re}\mu_z^-}+\frac{n_4}{\operatorname{Re}\mu_z^+}\right)^2=\frac{1}{m^2(1-\delta^2)^2}\left[\sum_{k=-4}^4 \mathsf{P}_k \tau^k+\mathcal{O}(\tau^{-6})\right]\] where
\begin{gather*}
\mathsf{P}_4    = 4 ,\qquad 
\mathsf{P}_2    = 8 \left(\delta^2-m^2+1\right),\qquad
\mathsf{P}_0    = 4 \left(\delta^4+\delta^2 \left(3 m^2-2\right)+m^4-m^2+1\right),\\
\mathsf{P}_{-2} = -16 \delta^2 \left(\delta^2+3\right) m^2, \qquad
\mathsf{P}_{-4} = -4 m^2 \left(-4 \delta^6+8 \delta^4 \left(m^2-5\right)+5 \delta^2 \left(5
    m^2-4\right)\right).
\end{gather*}
\subsubsection*{Coefficients for $\sqrt{\sqrt{(A_z-C_z)^2+B_z^2}+A_z+C_z+2D_z}$}
We use the following,
\[A_z+C_z+2D_z=\frac{1}{m^2(1-\delta^2)^2}\left[\sum_{k=-2}^4 \mathsf{P}_k \tau^k+\mathcal{O}(\tau^{-4})\right]\] where 
\begin{gather*}
\mathsf{P}_4    =1 ,\qquad
\mathsf{P}_2    =2 \left(\delta^2-m^2+1\right) ,\\
\mathsf{P}_0    = \frac12\left(2 \delta^4+\delta^2 \left(5 m^2-4\right)+2 m^4-m^2+2\right),\qquad
\mathsf{P}_{-2} = -16 \delta^2 m^2,
\end{gather*} and
\[\sqrt{(A_z-C_z)^2+B_z^2}=\frac{1}{m^2(1-\delta^2)^2}\left[\sum_{k=-2}^4\mathsf{P}_k\tau_k+\mathcal{O}(\tau^{-4})\right]\]
where 
\begin{gather*}
\mathsf{P}_4    = 1,\qquad
\mathsf{P}_2    =2 (\delta^2-m^2+1),\\
\mathsf{P}_0    =\frac12\left(2 \delta^4+\delta^2 \left(5 m^2-4\right)+2 m^4-m^2+2\right),\qquad
\mathsf{P}_{-2} =-256\delta^2m^2.
\end{gather*}
Note that the first 3 coefficients of these two terms match exactly. This gives, 
\begin{equation} \label{eq asymp R1}
  \|T_1(z)\|\leq\frac{1}{m(1-\delta^2)}\left[\tau^2+\mathsf{P}_0+\mathsf{P}_{-2}\tau^{-2}+\mathcal{O}(\tau^{-4})\right],
\end{equation}
where
\[
\mathsf{P}_{0}=\frac14(3 + 2\delta^2 - 2m^2) \qquad \text{and} \qquad \mathsf{P}_{-2}=\frac14(3m^2 + \delta^2(9m^2-8)).
\]

\subsubsection*{Coefficients for $\sqrt{\sqrt{(A_z-C_z)^2+\widetilde{B}_z^2}+A_z+C_z+2D_z}$}
We use the fact that
\[
     \operatorname{Re}(\operatorname{csgn}(z))=\frac{\operatorname{Re}z}{|z|}
\]
to obtain the asymptotic
\[
    \operatorname{Re}\left(\operatorname{csgn}\frac{w_z^+-w_z^-}{w_z^\mp}\right)=\mp 1 \pm \frac{m^2}{2\tau^4}+\mathcal{O}(\tau^{-6}).
\]
Hence we get
\[\sqrt{(A_z-C_z)^2+\widetilde{B}_z^2}=\frac{1}{m^2(1-\delta^2)^2}\left[\sum_{k=-2}^4\mathsf{P}_k\tau_k+\mathcal{O}(\tau^{-4})\right]\]
where
\begin{gather*}
\mathsf{P}_4    = 1,\qquad
\mathsf{P}_2    =2 (\delta^2-m^2+1),\\
\mathsf{P}_0    = \delta^4+\delta^2 \left(3 m^2-2\right)+m^4-m^2+1,\qquad
\mathsf{P}_{-2} =-32 \delta^2 \left(\delta^2+7\right).
\end{gather*}
Note that only the two final coefficients are different from the ones of $B_z^2$.
Hence 
\begin{equation}   \label{eqasympR1}
  \|T_1(z)\|\geq\frac{1}{m(1-\delta^2)}\left[\tau^2+\mathsf{P}_0+\widetilde{\mathsf{P}}_{-2}\tau^{-2}+\mathcal{O}(\tau^{-4})\right]
\end{equation}
where
\[
\mathsf{P}_{0}=\frac14(3 + 2\delta^2 - 2m^2) \qquad \text{and} \qquad \widetilde{\mathsf{P}}_{-2}=\frac18(5m^2 + \delta^2(19m^2-16)).
\]
\subsubsection*{Completion of the proof}
According to \eqref{eqforR2}-\eqref{eq for R2}, we have
\[
\|T_2(z)\| \leq (n_{2}+n_{5})\left(\frac{1}{\Re \mu_{z}^{-}}+\frac{1}{\Re \mu_{z}^{+}}\right)=\frac{1}{1-\delta^2} \left(2+\mathcal{O}(\tau^{-4})\right),
\]
as $|\tau|\to \infty$. Thus, we complete the proof of Theorem~\ref{Theo sharper estimate}, by combining the asymptotic formulas \eqref{eq asymp R1}-\eqref{eqasympR1} with this and invoking the triangle inequality.
\begin{remark}  \label{Rem:more terms} We suspect that
with some more effort, it is possible to show that there exists $\mathsf{P}_0$ such that
\begin{equation*}\label{Accurate conjecture}
	\left\|(\mathscr{L}_{m}-z)^{-1}\right\| = \frac{1}{1-|\Im z|^2}\left(\frac{1}{m}|\Re z|^2+\frac{1}{2m}\big[|\Im z|^2 +\mathsf{P}_0\big]+\mathcal{O}\big(|\Re z|^{-2}\big)\right),
\end{equation*}
where $\mathsf{P}_0$ is a constant independent of $z$ such that
\[ |\mathsf{P}_0|\leq \frac{3}{2}+4m+m^2.\] 
Note that the only point left to show here is that the coefficient $\mathsf{P}_0$ actually exists. The bound on its magnitude follows directly from our Theorem~\ref{Theo sharper estimate}. Two comments on this are in place.
\begin{itemize}
\item To prove this claim, one possibility is to find sharp lower bounds for $\|T_2(z)\|$, then feed into the final step above. This might not be easy and a more refined construction of pseudomodes seems to be needed.
\item To confirm existence of $\mathsf{P}_0$ without control on $\|T_2(z)\|$, perhaps it is possible to begin by proving that the function
\[
     s_1(\tau,\delta)=\frac{1}{\|(\mathscr{L}_m-\tau-i\delta)^{-1}\|},
\] 
being the first singular value of the operator $\mathscr{L}_m-\tau-i\delta$, is real analytic in the variables $\tau\in\mathbb{R}$ and $\delta\in(-1,1)$. This might follow from the fact that $\mathscr{L}_m-z$ is an integral operator with coefficients analytic in $z$, but it is not automatic and an analysis of its multiplicity is required.
\end{itemize}
\end{remark}  


\section{Perturbations of $\mathscr{L}_m$ by a long-range potential}\label{Sec L1 Perturbation}
We now consider perturbations of the form
$\mathscr{L}_m+V$
where $V\in L^1(\mathbb{R},\mathbb{C}^{2\times 2})$.  To aid the reader through this section, we recall the hypotheses,
 \begin{itemize} 
\item[\ref{Ass norm less than 1}{\bf:}]  $\Vert V \Vert_{L^1}<1$,
\item[\ref{Ass in Lp}{\bf:}] $V \in L^1(\R,\C^{2\times 2})\cap L^p(\R,\C^{2\times 2})$ for some $p\in(1,\infty]$.
\end{itemize}
Either will be assumed throughout. 

Our main goal in the first part of the section is to identify a closed extension, $\mathscr{L}_{m,V}\supseteq \mathscr{L}_{m}+V$, which is well defined under either assumptions and preserves the essential spectra. Then, in the second part, we examine the spectrum and pseudospectrum of $\mathscr{L}_{m,V}$ inside the instability band, $\Sigma=\{|\Im z|<1\}$.

Subsection~\ref{Subsec Eigen estimate m>0} is devoted to the proof of the main Theorem~\ref{Theo Location m>0} about the localisation of the point spectrum for $m>0$. Subsection~\ref{secsteppot}  addresses confirmation that there can be potentials with infinitely many eigenvalues in the instability band.  Subsection~\ref{Subsec Improved dynamics} discusses the effect of weak coupling, when imposing further concrete hypotheses on the decay and regularity of $V$, and comprises the proof of Theorem~\ref{Theo Weak Coupling}. Finally the proof of Theorem~\ref{inclusion pseudospectrum} about the pseudospectrum is given in Subsection~\ref{Subsec pseudospectrum V}.


\subsection{Proof of Proposition~\ref{Prop Perturbation}}\label{Subsec L1 perturbation}
In order to show the existence of the closed extension $\mathscr{L}_{m,V}$,  and determine its spectrum and pseudospectrum, our main tool is the Birman-Schwinger-type principle formulated in \cite{Gesztesy-Latushkin-Mitrea-Zinchenko20}. We quote that result in Theorem~\ref{Theo Closed Exten} and apply it as follows. 

Set the Hilbert spaces $\mathcal{H}=\mathcal{K}=L^{2}(\R,\C^2)$. Set the unperturbed operator $H_{0}=\mathscr{L}_{m}$. The auxiliary factors $A$ and $B$ take the following form. Write
\[ V(x)=U(x)\vert V(x) \vert\]
in polar form, where $U(x)$ is the partial isometry and $|V(x)| \coloneqq  (V^*(x)V(x))^{1/2}$. We set
\begin{equation}\label{Decompose Matrix}
A(x)\coloneqq  \vert V(x) \vert^{1/2}	 \qquad \text{and} \qquad B(x)\coloneqq  U(x) \vert V(x) \vert^{1/2}       .
\end{equation}
Then, we have
\begin{equation}\label{Equalities norm}
\vert A(x) \vert_{\C^2\times \C^2} =\vert V(x) \vert_{\C^2\times \C^2}^{1/2}=\vert B(x) \vert_{\C^2\times \C^2}.
\end{equation}
The first equality follows from the property that the norm of the square root of a matrix equals the square root of the norm of the matrix itself. To establish the second equality in \eqref{Equalities norm}, observe that $\operatorname{Ker}(\vert V(x) \vert) \subset \operatorname{Ker}(\vert V(x)\vert^{1/2})$. Consequently, we have $\overline{\operatorname{Im}(\vert V(x)\vert^{1/2})}\subset  [\operatorname{Ker}(U(x))]^{\perp}$, because $\operatorname{Ker}(|V(x)|)=\operatorname{Ker}(U(x))$, \emph{cf.} \cite[Prop. 2.85]{Cheverry-Raymond21}). Hence, as $U(x)$ is a partial isometry, it follows that
\[ \vert B(x) \vert_{\C^2\times \C^2} = \sup_{|v|_{\C^2}=1} \left\vert U(x) |V(x)\vert^{1/2} v \right\vert_{\C^2}=\sup_{|v|_{\C^2}=1} \left\vert |V(x)\vert^{1/2} v \right\vert_{\C^2}= \vert A(x) \vert_{\C^2\times \C^2} =\vert V(x) \vert_{\C^2\times \C^2}^{1/2} .\]
The operators $A$ and $B$ denote the multiplication operators
\[
	Af(x) = A(x) f(x) \qquad \text{and} \qquad Bf(x)=B(x) f(x)
\]
with their maximal domains. By construction, they are both closed and densely defined.
For each $z \in \rho\left(\mathscr{L}_{m}\right)$, let 
\[
	R_{0}(z)\coloneqq  (\mathscr{L}_{m} - z)^{-1}\qquad \text{and}\qquad Q(z) \coloneqq  \overline{A R_{0}(z) B}. 
\]
This settles the framework of Theorem~\ref{Theo Closed Exten}.

\medskip

Now we give the proof of Proposition~\ref{Prop Perturbation}. We split it into three steps. 
In the first step we verify that the hypotheses
 \ref{Bounded Exten 1} and \ref{Bounded Exten 2} of Theorem~\ref{Theo Closed Exten} hold, for all $V\in L^1(\R,\C^{2\times 2})$. In the second step we verify that \ref{Nonempty} is implied by either of the two different assumptions, \ref{Ass norm less than 1} or \ref{Ass in Lp}. In the final step we derive the claims made about the spectrum.

\subsubsection*{Step~1: conditions \ref{Bounded Exten 1} and \ref{Bounded Exten 2}}  
Let $V\in L^1(\mathbb{R},\mathbb{C}^{2\times 2})$. Our goal is to show that,  for suitable $z\in \rho(\mathscr{L}_m)$,
\begin{itemize}
\item[\ref{Bounded Exten 1}] the operator $R_{0}(z)B$ is closable, and both $AR_{0}(z)$ and $\overline{R_{0}(z)B}$ are bounded;
\item[\ref{Bounded Exten 2}] the operator $AR_{0}(z)B$ has a bounded closure
$ Q(z)= \overline{AR_{0}(z)B}.$
\end{itemize}

Firstly, note that the Schur Test conditions,
\begin{equation}\label{G^n}
\begin{aligned}
&\int_{\R} \vert \mathcal{R}_{z}(x,y) \vert_{\C^{2 \times 2}}^n \, \dd y \leq K_{n}(z),\qquad \text{ for a.e. }x \in \R,\\
&\int_{\R} \vert \mathcal{R}_{z}(x,y) \vert_{\C^{2 \times 2}}^n \, \dd x \leq K_{n}(z),\qquad \text{ for a.e. }y \in \R,
\end{aligned}
\end{equation}
holds for all fixed $n\in \mathbb{N}$. The proof of this is the same as for the case of the power $n=1$, addressed in the verification of  \eqref{Schur cond}, as the integration only involves exponential functions. Here
\[ K_{n}(z)\coloneqq  \frac{2}{n\min \left\{\Re \mu_{z}^{-}, \Re \mu_{z}^{+}\right\}} \left(\sup_{t \in \R} \varphi_{z}(t)\right)^{n},\]
$\varphi_{z}$ is as in \eqref{varphi} and $K_{n}(z)$ is finite for each $z\in \rho(\mathscr{L}_{m})$. 

Now, since $H^{1}(\R,\C^2)\subset L^{\infty}(\R,\C^2)$ and $|V|^{\frac12}\in L^2(\mathbb{R},\mathbb{C}^2)$, it follows that $H^{1}(\R,\C^2)\subset \mathrm{Dom}(A)$. Given that \[R_{0}(z): L^{2}(\R,\C^2) \longrightarrow H^{1}(\R,\C^2),\] we have $\mathrm{Dom}(AR_{0}(z))=L^{2}(\R,\C^2)$. Therefore, $AR_{0}(z)$ is closed and hence bounded by the Closed Graph Theorem. This confirms the second requirement in \ref{Bounded Exten 1}. Moreover, $AR_{0}(z)$ is a Hilbert-Schmidt operator. Indeed, from \eqref{Equalities norm} and \eqref{G^n} with $n=2$, it follows that
\begin{equation} \label{bound AR}
\begin{aligned}
\int_{\R^2} \vert A(x)\mathcal{R}_{z}(x,y) \vert_{\C^{2 \times 2}}^2\, \dd x \, \dd y &\leq  \int_{\R^2} \vert V(x) \vert_{\C^{2 \times 2}} \vert \mathcal{R}_{z}(x,y) \vert_{\C^{2 \times 2}}^2\, \dd x\,\dd y\\
&\leq K_{2}(z) \Vert V \Vert_{L^1}.
\end{aligned}
\end{equation}
Similarly, we also have
\begin{equation} \label{bound RB}
\int_{\R} \vert \mathcal{R}_{z}(x,y)B(y) \vert_{\C^{2 \times 2}}^2\, \dd x \, \dd y\leq K_{2}(z)\Vert V \Vert_{L^1},
\end{equation}
and thus the integral operator 
\[
 f(\cdot) \longmapsto \int_{\R} \mathcal{R}_{z}(\cdot,y)B(y)f(y)\, \dd y,
\] 
is Hilbert-Schmidt. But the latter is an extension of the operator $R_{0}(z)B$, thus $R_{0}(z)B$ is closable. Furthermore, since  the domain $\mathrm{Dom}(R_{0}(z)B)=\mathrm{Dom}(B)$ is dense in $L^2(\mathbb{R};\mathbb{C}^2)$  and the integral operator above is continuous, $ \overline{R_{0}(z)B}$ coincides with it.
This completes the confirmation of \ref{Bounded Exten 1}. 

Now, we address \ref{Bounded Exten 2}. Once again we appeal to the Schur Test, in order to show that $AR_{0}(z)B$ has a bounded closure. But now we consider the weight,
\[ p(x)\coloneqq \left\{\begin{aligned}
&|V(x)|_{\C^{2\times 2}}^{1/2}\qquad &&\text{ if } |V(x)|_{\C^{2\times 2}}>0,\\
&1 \qquad &&\text{ if } |V(x)|_{\C^{2\times 2}}=0.
\end{aligned} \right. \]
According to \eqref{Equalities norm}, for almost every $x\in \R$, 
\begin{align*}
\int_{\R} \vert A(x)\mathcal{R}_{z}(x,y)B(y) \vert_{\C^{2 \times 2}}p(y)\, \dd y &\leq p(x) \int_{\R} \vert \mathcal{R}_{z}(x,y) \vert_{\C^{2\times 2}} p(y)^2\, \dd y\\
&=p(x) \int_{\R} \vert \mathcal{R}_{z}(x,y) \vert_{\C^{2\times 2}} \vert V(y)\vert_{\C^{2\times 2}}\, \dd y
\end{align*}
and, for almost every $y\in \R$,
\[ \int_{\R} \vert A(x)\mathcal{R}_{z}(x,y)B(y) \vert_{\C^{2 \times 2}}p(x)\, \dd x \leq p(y) \int_{\R} \vert \mathcal{R}_{z}(x,y) \vert_{\C^{2\times 2}}  \vert V(x)\vert_{\C^{2\times 2}}\, \dd x.\]
Since $V\in L^1(\R,\C^{2 \times 2})$ and $\vert \mathcal{R}_{z}(\cdot,\cdot) \vert_{\C^{2\times 2}}\in L^{\infty}(\R^2)$ (recall Lemma~\ref{Lemma Norm of G}), then by the Schur Test, the integral operator \[f(\cdot)\longmapsto  \int_{\R} A(\cdot)\mathcal{R}_{z}(\cdot,y)B(y) f(y)\, \dd y \]
is bounded. Thus, since $H^{1}(\R,\C^2)\subset \mathrm{Dom}(A)$,  $\mathrm{Dom}(AR_{0}(z)B)=\mathrm{Dom}(B)$, and since $\mathrm{Dom}(B)$ is dense, continuity implies that $Q(z)$ is equal to this integral operator. Moreover, 
\begin{equation}\label{Bound of Q(E)}
\Vert Q(z) \Vert \leq \sup_{y\in \R} \int_{\R} \vert V(x)\vert_{\C^{2\times 2}}\vert \mathcal{R}_{z}(x,y) \vert_{\C^{2\times 2}} \dd x.
\end{equation}
This ensures the validity of \ref{Bounded Exten 2}.

\subsubsection*{Step~2: hypothesis \ref{Nonempty} and existence of $\mathscr{L}_{m,V}$} We now need to verify that
\begin{itemize}
\item[\ref{Nonempty}] \qquad $\big\{z\in \rho(H_{0}): -1\in \rho(Q(z))\big\}\not=\varnothing$.
\end{itemize}
We split the proof into the two sub-cases depending on the hypothesis on $V$. 

First, assume that  \ref{Ass norm less than 1} is the one condition that holds. According to \eqref{Bound of Q(E)}, \eqref{Norm of G} and \eqref{varphi bounded}, for all $z\in \rho(\mathscr{L}_{m})$,
\begin{equation}\label{Upper bound norm Q}
\begin{aligned}
\Vert Q(z) \Vert \leq &\Vert V\Vert_{L^1} \sup_{(x,y)\in\R^2} |\mathcal{R}_{z}(x,y)|
\leq &\frac{\Vert V \Vert_{L^{1}}}{2} \left(|k_{z}|+1\right)\max \left\{ |w_{z}^{\pm}| +\frac{1}{|w_{z}^{\pm}|}\right\}.
\end{aligned}
\end{equation}
From \eqref{Notation pm} and \eqref{Function k}, it follows that $w_{i\delta}^{\pm} \to -i$ and $k_{i\delta}\to 0$ as $\delta\to \infty$. Thus, the terms inside the maximum become close to $2$ as $\delta\to \infty$. Since $\Vert V \Vert_{L^1}<1$, then indeed \ref{Nonempty} follows taking $z$ on the imaginary axis with sufficiently large modulus.
 
Now, if instead condition \ref{Ass in Lp} holds for $p\in (1,\infty]$, by Holder's inequality and  \eqref{G^n}, we get 
\begin{equation}\label{Upper bound norm Q Lp}
\Vert Q(z) \Vert \leq \Vert V(x) \Vert_{L^{p}} \left(K_{q}(z)\right)^{1/q}, \qquad \frac{1}{p}+\frac{1}{q}=1.\end{equation}
Therefore, a similar argument as the previous sub-case, yields  $K_{q}(z) \to 0$ for $z$ on the imaginary axis with sufficiently large modulus. Hence, \ref{Nonempty} is also valid under assumption \ref{Ass in Lp}. 

\medskip

As a conclusion to this step we now have that, by virtue of the first part of Theorem~\ref{Theo Closed Exten},  a closed extension $\mathscr{L}_{m,V}\supseteq \mathscr{L}_{m}+V$ indeed exists and
\[ (\mathscr{L}_{m,V}-z)^{-1}=R_{0}(z)-\overline{R_{0}(z)B}\left( I+Q(z)\right)^{-1}AR_{0}(z)\]
for all $z$ such that $-1 \in \rho(Q(z))$. Note that, under the hypotheses, the latter is non-empty.

\subsubsection*{Step~3: conclusions \ref{Prop3.2 1}~--~\ref{Prop3.2 4} of the proposition}
The crucial point in the proof is to observe that  $AR_{0}(z)$ and $\overline{R_{0}(z)B}$ are Hilbert-Schmidt operators, therefore the difference \[(\mathscr{L}_{m}-z)^{-1}-(\mathscr{L}_{m,V}-z)^{-1}=\overline{R_{0}(z)B}\left( I+Q(z)\right)^{-1}AR_{0}(z)\] is compact. The stability of the first four essential spectra follows directly from \cite[Thm. IX.2.4]{Edmunds-Evans18}. Moreover, when $m>0$, $\C\setminus \operatorname{Spec}_{\mathrm{e1}}(\mathscr{L}_{m,V})$ has one connected component, then $\operatorname{Spec}_{\mathrm{e5}}(\mathscr{L}_{m,V})$ in that case is also identical to the other essential spectra, \emph{cf.} \cite[Prop. 5.5.4]{Krejcirik-Siegl15}. This confirms the conclusions \ref{Prop3.2 1} and \ref{Prop3.2 2}. 

The rest of the proof assumes that $m=0$. In that case, notice that $k_{z}=0$ and $|w_{z}^{\pm}|=1$. For \ref{Prop3.2 3}, observe that by direct substitution into \eqref{Upper bound norm Q}, we have that, if $V$ satisfies \ref{Ass norm less than 1}, then $\Vert Q(z) \Vert\leq \Vert V\Vert_{L^1}<1$ for all \[z\in \rho(\mathscr{L}_{m})=\C \setminus \{z\in \C: |\Im z|\leq 1\}.\] Hence $-1 \in \rho(Q(z))$ for all $z\in \C$ such that $|\Im z|>1$. That is, any eigenvalue of $\mathscr{L}_{m,V}$ should all be located in the strip $\{z\in \C: |\Im z\leq 1\}$ as claimed.

Finally, for \ref{Prop3.2 4}, if $V$ satisfies \ref{Ass in Lp} instead, the estimate \eqref{Upper bound norm Q Lp}, the fact that  \[\min \Big\{\Re \mu_{z}^{+},\Re \mu_{z}^{-}\Big\} = |\Im z|-1 \qquad \text{for}\qquad |\Im z|>1,\] and that $\varphi_{z}(t)=1$ for all $t\in \R$, ensure that\[ \Vert Q(z) \Vert \leq \Vert V \Vert_{L^p}\left(\frac{2}{q(|\Im z|-1)}\right)^{1/q}, \qquad q= \frac{p}{p-1},\] for all $z\in \rho(\mathscr{L}_{0})=\{z\in\C: |\Im z|>1\}$.
Hence, for all $z\in \C$ such that $|\Im z|>1+\frac{2}{q}\Vert V \Vert_{L^p}^{q}$, the right-hand side is strictly smaller than 1. This gives \ref{Prop3.2 4}  for $p<\infty$. For the case $p=\infty$, note that $\operatorname{Spec}(\mathscr{L}_{0,V})\subset \operatorname{Spec}_{\|V\|_{\infty}}(\mathscr{L}_{0})$ (from the second equality in \eqref{Def Pseudospectrum}). This completes the proof of Proposition~\ref{Prop Perturbation}. 

\medskip

For later purposes, note that the operator $Q(z)$ is also Hilbert-Schmidt. Indeed, the identities \eqref{Equalities norm}, \eqref{Norm of G} and \eqref{varphi bounded}, yield
\begin{equation} \label{Hilbert-Schmidt norm}
\begin{aligned}
&\int_{\R^2} \vert A(x)\mathcal{R}_{z}(x,y)B(y) \vert_{\C^{2 \times 2}}^2\, \dd x \, \dd y\\
&\leq \Vert V \Vert_{L^{1}}^2 \sup_{(x,y)\in \R^2} \vert \mathcal{R}_{z}(x,y) \vert_{\C^{2\times 2}}^2\\
&\leq \frac{\Vert V \Vert_{L^{1}}^2}{4} \max \left\{\left(|k_{z}|+1\right)^2 \left( |w_{z}^{-}| +\frac{1}{|w_{z}^{-}|}\right)^2, \left(|k_{z}|+1\right)^2 \left( |w_{z}^{+}| +\frac{1}{|w_{z}^{+}|}\right)^2\right\},
\end{aligned}
\end{equation}
where the right-hand side is finite for each $z\in \rho(\mathscr{L}_{m})$.

We close this subsection with a remark that provides the context of our findings about perturbations of $\mathscr{L}_m$ by long-range potentials.

\begin{remark}\label{Remark Compact Q}
In the case of long-range perturbations of the dislocated Schr{\"o}dinger operator $-\frac{\dd^2}{\dd x^2}+i\sgn(x)$, considered in \cite[Subsec. 5.1]{Henry-Krejcirik17}, the condition $V\in L^1(\R)$ alone gives the existence of a closed extension directly, via the associated quadratic forms. The distinction with the present case and the need for the condition \ref{Ass norm less than 1} can be explained through the  behaviour of the resolvent kernel, as follows. 

The estimates around \cite[Est. (5.6)]{Henry-Krejcirik17}, ensure that the $L^{\infty}$ norm of the kernel of the unperturbed resolvent for the Schr{\"o}dinger model decays to zero as $z \to \infty$ along the imaginary axis. By contrast, the kernel $\mathcal{R}_{z}(x,y)$ in the present case does not have a decaying $L^{\infty}$ norm. Effectively, the right-hand sides in \eqref{Norm of G} and \eqref{varphi bounded} do not decay to zero. Note that this is analogous to the contrast between the free Schr{\"o}dinger kernel and the free particle Dirac kernel. For details, see \cite[Form. (13)--(14)]{Cuenin-Laptev-Tretter14}.
\end{remark}


\subsection{Proof of Theorem \ref{Theo Location m>0}} \label{Subsec Eigen estimate m>0}
The Proposition~\ref{Prop Perturbation}, which has been already proven, provides a general statement about exclusion regions (semi-planes) for the eigenvalues of $\mathscr{L}_{0,V}$.  We now proceed to give the proof of the first main result of this section, Theorem \ref{Theo Location m>0}. This theorem, in a similar vein, describes regions of exclusion for the eigenvalues of $\mathscr{L}_{m,V}$ in the case $m>0$.  Thus, for the remainder of this subsection, we assume that $m\not= 0$.

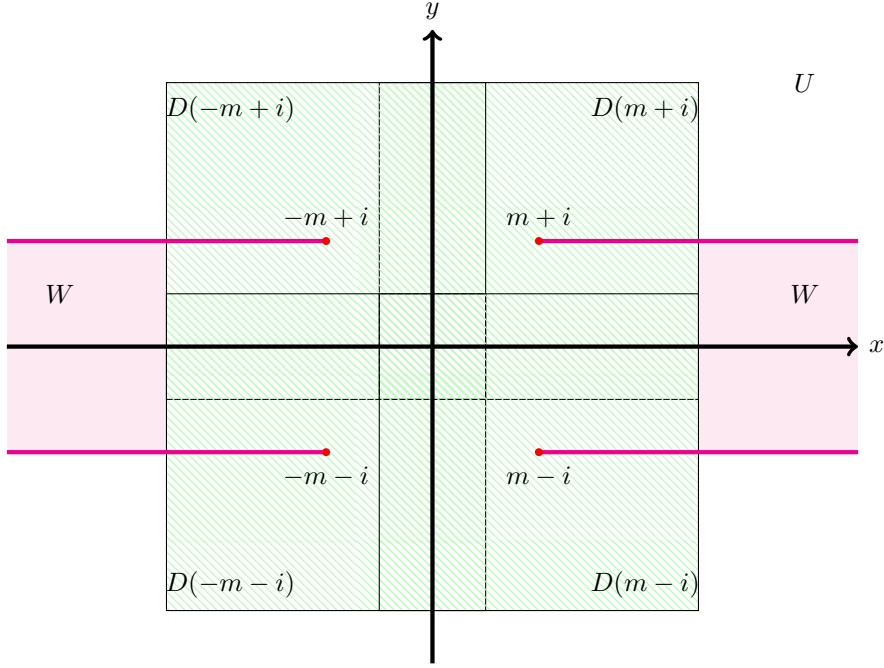
\begin{figure}[h]
\centering
\begin{tikzpicture}[scale=0.7]
\fill[ fill=magenta!10](-8,-2) rectangle (-5,2);
\fill[ fill=magenta!10](5,-2) rectangle (8,2);
\draw[draw=black, pattern=north west lines, pattern color=green!30] (-1,-1) rectangle (5,5);
\draw[draw=black, pattern=north west lines, pattern color=green!30] (-5,-1) rectangle (1,5);
\draw[draw=black, pattern=north west lines, pattern color=green!30] (-5,-5) rectangle (1,1);
\draw[draw=black, pattern=north west lines, pattern color=green!30]  (-1,-5) rectangle (5,1);
\draw[-,ultra thick, ultra thick, magenta] (-8,2)--(-2,2);
\draw[-,ultra thick, ultra thick, magenta] (-8,-2)--(-2,-2);
\draw[-,ultra thick, ultra thick, magenta] (2,2)--(8,2);
\draw[-,ultra thick, ultra thick, magenta] (2,-2)--(8,-2);
\draw[->,ultra thick] (-8,0)--(8,0) node[right]{$x$};
\draw[->,ultra thick] (0,-6)--(0,6) node[above]{$y$};
\node[circle,fill=red,inner sep=0pt,minimum size=3pt,label=above:{$m+i$}] () at (2,2) {};
\node[circle,fill=red,inner sep=0pt,minimum size=3pt,label=above:{$-m+i$}] () at (-2,2) {};
\node[circle,fill=red,inner sep=0pt,minimum size=3pt,label=below:{$-m-i$}] () at (-2,-2) {};
\node[circle,fill=red,inner sep=0pt,minimum size=3pt,label=below:{$m-i$}] () at (2,-2) {};
\node[black] at (7,5) {$U$};
\node[black] at (4,4.5) {$D(m+i)$};
\node[black] at (-3.8,4.5) {$D(-m+i)$};
\node[black] at (-3.8,-4.5) {$D(-m-i)$};
\node[black] at (4,-4.5) {$D(m-i)$};
\node[black] at (7,1) {$W$};
\node[black] at (-7,1) {$W$};
\end{tikzpicture}
\caption{Partition of the resolvent set of $\mathscr{L}_{m}$ in subdomains for the proof of Theorem~\ref{Theo Location m>0}.}\label{Figure Partition Resolvent}
\end{figure}

 Our main arguments involve deriving upper bound for the norm of the \emph{Birman-Schwinger operator}, $Q(z)$, introduced in \eqref{Bound of Q(E)}. To achieve this, we give sharp estimates on the matrix norm of the kernel in various regions of the complex $z$-plane. Specifically, we partition the resolvent set into three main disjoint regions as follows: \[ \rho(\mathscr{L}_{m})=D \cup W \cup U,\]
where
\begin{align*}
W \coloneqq  &\big\{z \in \C: |\Re(z)| \geq 5m/2, \, | \Im z | <1 \big\}, \\
D \coloneqq  &D(m+i)\cup D(-m+i)\cup D(-m-i)\cup D(m-i) \qquad \text{for}\\
&\qquad D(m\pm i)\coloneqq  \left\{z\in \C\setminus \rho(\mathscr{L}_{m}):\begin{aligned}
&\vert \Re(z-(m \pm i))\vert\leq 3m/2,\\
&\vert\Im(z-(m \pm i)) \vert\leq 3/2 
\end{aligned}\right\} \qquad \text{and}\\
&\qquad D(-m\pm i)\coloneqq  \left\{z\in \C\setminus \rho(\mathscr{L}_{m}):\begin{aligned}
&\vert \Re(z-(-m \pm i))\vert\leq 3m/2,\\
&\vert\Im(z-(-m\pm i)) \vert\leq 3/2 
\end{aligned}\right\}, \qquad \text{and}\\
U \coloneqq  &\rho(\mathscr{L}_{m})\setminus (D\cup W).
\end{align*}
See Figure~\ref{Figure Partition Resolvent}.

We split the proof of Theorem~\ref{Theo Location m>0} according to this partition, deriving suitable estimates for the norm of $Q(z)$ where $z$ is in each of the sub-regions.
In this proof, the parameters $c_j\equiv c_j(m)>0$ only depend on $m$.

\subsubsection*{Estimates on $D$} The region $D\subset \mathbb{C}$ is the union of four neighbourhoods of  singularities of the matrix kernel $\mathcal{R}_{z}$. We will see that these singularities are removable. That is, the matrix norm $|\mathcal{R}_{z}(x,y)|_{\mathbb{C}^{2\times 2}}$ is bounded uniformly for all $(x,y)\in \R^2$, whenever $z\in D(\pm m \pm i)$. We treat each different component separately.

Consider first that $z\in D( m + i)$. From the formulae of $\mu_{z}^{+}$, $w_{z}^{\pm}$ in \eqref{Notation pm} and $k_{z}$ in \eqref{Function k}, it follows that 
\begin{equation}\label{E to m+i}
\begin{aligned}
&\lim_{z \to m+i}w^+_z = 0 , \qquad\lim_{z \to m+i} w^-_z = \frac{\sqrt{-1-im}}{\sqrt{m^2+1}},\\
&\lim_{z \to m+i}\frac{\mu_{z}^{+}}{w_{z}^{+}}=2m,\qquad \lim_{z \to m+i} k_{z}=-1.
\end{aligned}
\end{equation}
Then, we gather the following. According to \eqref{Norm of G}, for all $(x,y)\in \R^2$ such that $xy\leq 0$,
\begin{equation}\label{E to m+i 1}
\vert \mathcal{R}_{z}(x,y) \vert_{\C^{2\times 2}} \leq \varphi_{z}(0)= \frac{\sqrt{(1+|w_{z}^{+}|^2)(1+|w_{z}^{-}|^2)}}{|w_{z}^{+}+w_{z}^{-}|} \leq c_1.
\end{equation}
According to \eqref{varphi bounded}, for all $(x,y)\in \R^2$ such that $x\leq 0, y\leq 0$,
\begin{equation}\label{E to m+i 2}
\vert \mathcal{R}_{z}(x,y) \vert_{\C^{2\times 2}} \leq \sup_{t\leq 0}\varphi_{z}(t)\leq \frac{1}{2}\left(|k_{z}|+1\right) \left( |w_{z}^{-}| +\frac{1}{|w_{z}^{-}|}\right)\leq c_2.
\end{equation}
According to \eqref{Norm of G}, for all $(x,y)\in \R^2$ such that $x\geq 0, y\geq 0$, 
\begin{align*}
\vert \mathcal{R}_{z}(x,y) \vert_{\C^{2\times 2}} &\leq  \varphi_{z}(\min\{x,y\}).
\end{align*}
Now, for all $t\geq 0$, we know that
\[\varphi_{z}(t)= \frac{\sqrt{1+\vert w_{z}^{+}\vert^2}}{2}\left(\frac{1}{\vert w_{z}^{+}\vert^2}\left\vert k_{z}e^{-2\mu_{z}^{+}t} +1\right\vert^2+\left\vert k_{z}e^{-2\mu_{z}^{+}t} -1\right\vert^2 \right)^{1/2}.\]
Because of \eqref{E to m+i}, the term that we should take care of is $\frac{1}{\vert w_{z}^{+}\vert}\left\vert k_{z}e^{-2\mu_{z}^{+}t} +1\right\vert$, as it might potentially exhibit a singularity. We show that this is not the case. Indeed, from \eqref{Function k} and \eqref{E to m+i}, for $z\in D(m+i)$ and $x \geq 0, y\geq 0$, we have
\begin{align*}
\frac{1}{\vert w_{z}^{+}\vert}\left\vert k_{z}e^{-2\mu_{z}^{+}t} +1\right\vert &= \frac{1}{\vert w_{z}^{+}\vert}\left\vert k_{z}(e^{-2\mu_{z}^{+}t}-1) +1+k_{z}\right\vert\\
&= \frac{1}{\vert w_{z}^{+}\vert}\left\vert k_{z}(e^{-2\mu_{z}^{+}t}-1) +\frac{2w_{z}^{+}}{w_{z}^{+}+w_{z}^{-}}\right\vert\\
&\leq \frac{| \mu_{z}^{+}|}{|w_{z}^{+}|}|k_{z}| |t|+\frac{2}{|w_{z}^{+}+w_{z}^{-}|}\\
&\leq c_3( |t|+1).
\end{align*}
Here we have used the inequality $|e^{-z}-1|\leq |z|$ for $\Re z\geq 0$. Thus, for all $z\in D(m+i)$ and $x\geq 0$, $y\geq0$, we have 
\begin{equation}\label{E to m+i 3}
\vert \mathcal{R}_{z}(x,y) \vert_{\C^{2\times 2}} \leq \varphi_{z}(\min\{x,y\}) \leq c_4 (1+|\min\{x,y\}|) \leq c_5 \nu_1'(x).
\end{equation}
Finally, by virtue of \eqref{Bound of Q(E)}, combining the above estimates we obtain,
\begin{equation}\label{Estimate on D(m+i)}
\Vert Q(z) \Vert \leq c_6 \int_{\R} |V(x)|_{\C^{2\times 2}}\nu_1'(x)\, \dd x
\end{equation}
for all $z\in D( m + i)$.

Now let $z\in D(m-i)$. We have
\begin{equation}
\begin{aligned}
&\lim_{z\to m-i} w_{z}^{-}=0, \qquad \lim_{z\to m-i} w_{z}^{+}=\frac{\sqrt{-1+im}}{\sqrt{m^2+1}},\\
&\lim_{z \to m-i} \frac{\mu_{z}^{-}}{w_{z}^{-}}=2m, \qquad\lim_{z \to m-i} k_{z}=1.
\end{aligned}
\end{equation}
By means of an identical proof to the one given in the previous region, we can show that
if $xy\leq 0$ then \eqref{E to m+i 1} holds true and if $x \geq 0, y\geq 0$ then \eqref{E to m+i 2} holds true. We omit the details of that. Moreover, in the case $x\leq 0$ and $y\leq 0$, the term $\frac{1}{\vert w_{z}^{-}\vert}\left\vert k_{z}e^{2\mu_{z}^{-}t} -1\right\vert$ is bounded by $c_7(|t|+1)$. Thus, we gather that
\begin{align*}
\vert \mathcal{R}_{z}(x,y) \vert_{\C^{2\times 2}} &\leq \varphi_{z}(\max\{x,y\}) \leq c_8( 1+|\max\{x,y\}|)\\
&\leq c_9 (1+\min\{|x|,|y|\})\leq c_{10} \nu_1'(x).
\end{align*}
Hence, by virtue of \eqref{Bound of Q(E)}, we also obtain the estimate \eqref{Estimate on D(m+i)} for all $z\in D(m-i)$, perhaps with a larger constant $c_6$.

Now, let $z\in D(-m+i)$. Then,
\begin{equation}
\begin{aligned}
&\lim_{z \to -m+i}w^+_z = \infty , \qquad &&\lim_{z \to -m+i} w^-_z = \sqrt{-1-im},\\
&\lim_{z \to -m+i}\mu_{z}^{+}w_{z}^{+}=2m,\qquad &&\lim_{z \to -m+i} k_{z}=1.
\end{aligned}
\end{equation}
We obtain the same estimates as \eqref{E to m+i 1} for all $(x,y)\in \R^2$ satisfying $xy\leq 0$ and \eqref{E to m+i 2} for all $(x,y)\in \R^2$ satisfying $x\leq 0$, $y\leq 0$,
by writing
\begin{equation}
\begin{aligned}
  \frac{\sqrt{(1+|w_{z}^{+}|^2)(1+|w_{z}^{-}|^2)}}{|w_{z}^{+}+w_{z}^{-}|} =\frac{\sqrt{\left(\frac{1}{|w_{z}^{+}|^2}+1\right)(1+|w_{z}^{-}|^2)}}{\left\vert 1+\frac{w_{z}^{-}}{w_{z}^{+}}\right\vert}.
\end{aligned}
\end{equation}
On the other hand, whenever $x \geq 0$, $y\geq 0$, we invoke the same estimate $\vert \mathcal{R}_{z}(x,y) \vert_{\C^{2\times 2}} \leq  \varphi_{z}(\min\{x,y\})$ as on $D(m+i)$, in which, we rewrite
\[
\varphi_{z}(t)= \frac{1}{2}\left(\frac{1+\vert w_{z}^{+}\vert^2}{\vert w_{z}^{+}\vert^2}\left\vert k_{z}e^{-2\mu_{z}^{+}t} +1\right\vert^2+ \left(1+\vert w_{z}^{+}\vert^2\right)\left\vert k_{z}e^{-2\mu_{z}^{+}t} -1\right\vert^2 \right)^{1/2}.
\]
It is readily seen that, for all $t\geq 0$ and $z\in D(-m+i)$, 
\[ \frac{1+\vert w_{z}^{+}\vert^2}{\vert w_{z}^{+}\vert^2}\left\vert k_{z}e^{-2\mu_{z}^{+}t} +1\right\vert^2+\left\vert k_{z}e^{-2\mu_{z}^{+}t} -1\right\vert^2\leq c_{11}.\]
For the remain term $ \vert w_{z}^{+}\vert^2\left\vert k_{z}e^{-2\mu_{z}^{+}t} -1\right\vert^2$, note that
\begin{align*}
\vert w_{z}^{+}\vert\left\vert k_{z}e^{-2\mu_{z}^{+}t} -1\right\vert &= \vert w_{z}^{+}\vert\left\vert k_{z}(e^{-2\mu_{z}^{+}t}-1) +k_{z}-1\right\vert\\
&= \vert w_{z}^{+}\vert\left\vert k_{z}(e^{-2\mu_{z}^{+}t}-1) -\frac{2w_{z}^{-}}{w_{z}^{+}+w_{z}^{-}}\right\vert\\
&\leq | \mu_{z}^{+} w_{z}^{+}||k_{z}| |t|+\frac{2|w_{z}^{-}|}{\left|1+ \frac{w_{z}^{-}}{w_{z}^{+}} \right|}\\
&\leq c_{12}( |t|+1).
\end{align*}
Thus, the estimate \eqref{E to m+i 3} is valid for all  $(x,y)\in \R^2$ such that $x\geq 0$, $y\geq 0$, also whenever $z\in D(-m+i)$. Hence, combining the estimates, once again we gather that  \eqref{Estimate on D(m+i)} also holds true for all $z\in D(-m+i)$.

Finally, the case $z\in D(-m-i)$ is similar to the case $z\in D(-m+i)$, so we omit further details. Hence, we then gather that the upper bound,
\begin{equation}\label{Estimate on D}
\Vert Q(z) \Vert \leq c_{13}(m) \|V\|_{L^1(\mathrm{d}\nu_1)}
\end{equation}
is valid for all  $z\in D$.

\subsubsection*{Estimate on $U=\C \setminus \overline{(D\cup W)}$} 
Split $\overline{U}=U^+\cup U^-$ for
\[
U^{\pm}\coloneqq \{ z\in \overline{U}: \pm \Im z>0\}.
\]
Then, according to \eqref{Notation pm} and \eqref{Function k}, 
\begin{align*}
\lim_{\substack{z\to \infty \\ z\in U^{+}}} w_{z}^{\pm}= -i,\qquad \lim_{\substack{z\to \infty \\ z\in U^{-}}} w_{z}^{\pm}=i \qquad \text{and} \qquad \lim_{\substack{z\to \infty \\ z\in\overline{U}}}k_z=0.
\end{align*}
 Hence, by continuity of each of the quantities $|w_{z}^{\pm}|$, $|w_{z}^{\pm}|^{-1}$ and $|k_{z}|$, in $z$, we have that
\begin{equation} \label{For the estimate on U}
\big(|k_z|+1\big)\max\left\{|w_{z}^\pm|+\frac{1}{|w_{z}^{\pm}|}\right\}\leq c_{14}(m)
\end{equation}
for all $z\in U$.
According to  \eqref{Bound of Q(E)}, \eqref{Norm of G} and \eqref{varphi bounded},  we then gather that
\begin{equation}\label{Estimate on U}
\Vert Q(z) \Vert \leq c_{14}(m)  \Vert V \Vert_{L^1}\leq c_{14}(m)\|V\|_{L^1(\mathrm{d}\nu_1)}
\end{equation}
for all $z\in U$.

\subsubsection*{Estimate on $W$} By invoking Lemma~\ref{Lem Asymptotic}, we obtain the following bounds for all $z\in W$;
\begin{align*}
&\frac{\sqrt{(1+|w^+_z|^2) (1+|w^-_z|^2)}}{|w^+_z+w^-_z|} \leq c_{15} (\Re z)^2,\qquad |w_{z}^{\pm}| \leq c_{16},\qquad |w_{z}^{\pm}|^{-1} \leq c_{17},\qquad |k_{z}|\leq c_{18}(\Re z)^2.
\end{align*}
Then, by virtue of \eqref{Bound of Q(E)}, \eqref{Norm of G} and \eqref{varphi bounded},
it follows that
\begin{equation}\label{Estimate on W}
 \Vert Q(z) \Vert \leq c_{19}(m) (\Re z)^2\Vert V \Vert_{L^1}
\end{equation}
 for all $z\in W$.

\medskip

 \subsubsection*{Completion of the proof of Theorem~\ref{Theo Location m>0}} According to the estimates \eqref{Estimate on D} on $D$, \eqref{Estimate on U} on $U$ and \eqref{Estimate on W} on $W$, we conclude that there exists a constant $C(m)>0$ such that
\begin{equation}\label{Bound out W}
\Vert Q(z) \Vert \leq \frac{1}{C(m)} \Vert V \Vert_{L^1(\mathrm{d}\nu_1)} \qquad \text{for all }z\in \rho(\mathscr{L}_{m})\setminus W,
\end{equation}
and
\begin{equation}\label{Bound on W}
\Vert Q(z) \Vert \leq \frac{1}{C(m)^2} (\Re z)^2 \Vert V \Vert_{L^1}\qquad \text{ for all }z \in W.
\end{equation}
Hence, if we assume $\Vert V \Vert_{L^1(\mathrm{d}\nu_1)}<C(m)$, it follows that $\Vert Q(z) \Vert<1$ for all $z\in \rho(\mathscr{L}_{m})\setminus W$. Hence, $-1$ is not in the spectrum of $Q(z)$ whenever $z\in \rho(\mathscr{L}_{m})\setminus W$. Theorem~\ref{Theo BS}, thus implies that all points in $\rho(\mathscr{L}_{m})\setminus W $ are not eigenvalues of $\mathscr{L}_{m,V}$. Moreover,  the points $z\in W$ satisfying $|\Re z|< C(m) \sqrt{\frac{1}{\Vert V \Vert_{L^1}}}$ cannot be eigenvalues of  $\mathscr{L}_{m,V}$. Therefore, we deduce that the eigenvalues of $\mathscr{L}_{m,V}$ (if there is any), should lie in the subset of $W$, corresponding to real part larger than or equal to $C(m)\sqrt{\frac{1}{\Vert V \Vert_{L^1}}}$. This completes the proof  of Theorem~\ref{Theo Location m>0}.


\subsection{Step potentials}\label{secsteppot}
In this subsection we construct a potential $V\in L^1$ such that $\mathscr{L}_{m,V}$ has infinitely many eigenvalues inside the instability band. We compute explicitly the eigenfunctions by means of an argument similar to the one employed in the proof of Proposition~\ref{Prop Resolvent}. We end the subsection by giving the proof of Proposition~\ref{Prop Steplike}. 

Set 
\[ V_{a,b}(x)\coloneqq (-i\sgn(x)-b)\chi_{[-a,a]}(x)  I, \]
for $a>0$ and $b\in \R$. If $u\not=0$ is such that $(\mathscr{L}_{m,V_{a,b}}-z)u=0$, then 
\[ u(x)=\left\{ \begin{aligned}
&u_{-}(x) \qquad && \text{ for }x\in(-\infty,-a),\\
&u_{0}(x)\qquad && \text{ for }x\in[-a,a],\\
&u_{+}(x)\qquad && \text{ for }x\in(a,+\infty),
\end{aligned}\right.\]
where $u_{\dagger}$ for $\dagger\in\{-,0,+\}$, are given similarly as in \eqref{u pm}, by
\begin{equation}\label{u0}
u_{\dagger}(x)=
\left(e^{\mu_{z}^{\dagger}x}S_{z}^{\dagger}+e^{-\mu_{z}^{\dagger}x}T_{z}^{\dagger}\right)\begin{pmatrix}
\alpha_{\dagger}\\
\beta_{\dagger}
\end{pmatrix}.
\end{equation}
Here, $\mu_z^{\pm}$ and $w_z^{\pm}$ are the same parameters from \eqref{Notation pm} that we considered above, 
\[
\mu_{z}^{0}\coloneqq  \sqrt{(m+b+z)(m-b-z)}\qquad\text{and} \qquad w_{z}^{0}\coloneqq   \frac{\sqrt{m-b-z}}{\sqrt{m+b+z}}.
\]
Also
\[
S^{\dagger}_{z}\coloneqq  \frac{1}{2}\begin{pmatrix}
1 && 1/w_{z}^{\dagger}\\
w_{z}^{\dagger} && 1
\end{pmatrix}\qquad\text{and} \qquad T_{z}^{\dagger}\coloneqq  \frac{1}{2}\begin{pmatrix}
1 && -1/w_{z}^{\dagger}\\
-w_{z}^{\dagger} && 1
\end{pmatrix}.
\]
Note that there is a singularity of $u_{0}(x)$ for $z=\pm m -b$, but taking the limits $z\to \pm m-b$ in \eqref{u0} for $\dagger=0$, we get
\[ u_{0}(x) = \left\{
\begin{aligned}
&\begin{pmatrix}
1 & 2mx\\
0 & 1
\end{pmatrix}\begin{pmatrix}
\alpha_{0}\\
\beta_{0}
\end{pmatrix},\qquad &&\text{if }z=m-b,\\
&\begin{pmatrix}
1 & 0\\
2mx & 1
\end{pmatrix}\begin{pmatrix}
\alpha_{0}\\
\beta_{0}
\end{pmatrix},\qquad &&\text{if }z=-m-b.
\end{aligned}
\right. \]
We seek for $z\in\mathbb{C}$, $\alpha_{\dagger}$ and $\beta_{\dagger}$ (non-vanishing simultaneously), such that $u\in H^1(\mathbb{R},\mathbb{C}^2)$, so it becomes a proper eigenfunction.

Since $u$ has to decay at $\pm \infty$ and must be continuous at $\pm a$, then 
\[
\begin{pmatrix}
-w_{z}^{-} & 1 & 0 & 0 & 0 & 0\\
0		   & 0 & 0 & 0 & w_{z}^{+} & 1\\
e^{-\mu_{z}^{-}a} & \frac{e^{-\mu_{z}^{-}a}}{w_{z}^{-}} &-(e^{\mu_{z}^{0}a}+e^{-\mu_{z}^{0}a}) & \frac{(e^{\mu_{z}^{0}a}-e^{-\mu_{z}^{0}a})}{w_{z}^{0}} &0 &0	\\
w_{z}^{-}e^{-\mu_{z}^{-}a} & e^{-\mu_{z}^{-}a} &w_{z}^{0}(e^{\mu_{z}^{0}a}-e^{-\mu_{z}^{0}a}) & -(e^{\mu_{z}^{0}a}+e^{-\mu_{z}^{0}a}) &0 &0\\
0 & 0 & (e^{\mu_{z}^{0}a}+e^{-\mu_{z}^{0}a}) & \frac{(e^{\mu_{z}^{0}a}-e^{-\mu_{z}^{0}a})}{w_{z}^{0}} & -e^{-\mu_{z}^{+}a} & \frac{e^{-\mu_{z}^{+}a}}{w_{z}^{+}}\\
0 & 0 & w_{z}^{0}(e^{\mu_{z}^{0}a}-e^{-\mu_{z}^{0}a}) & (e^{\mu_{z}^{0}a}+e^{-\mu_{z}^{0}a}) & w_{z}^{+} e^{-\mu_{z}^{+}a} &-  e^{-\mu_{z}^{+}a}
\end{pmatrix}\begin{pmatrix}
\alpha_{-}\\
\beta_{-}\\
\alpha_{0}\\
\beta_{0}\\
\alpha_{+}\\
\beta_{+}
\end{pmatrix}=0.
\]
The determinant of the matrix vanishes, if and only if,
\begin{equation}\label{Eigen Eq Step}
 e^{4a \mu_{z}^{0}}(w_{z}^{0}+w_{z}^{+})(w_{z}^{0}+w_{z}^{-})=(w_{z}^{0}-w_{z}^{+})(w_{z}^{0}-w_{z}^{-}) 
\end{equation}
for eigenvalues $z\not = \pm m -b$ such that \[z\not \in \{z\in \C: |\Re z|\geq m, |\Im z|=1\}.\] And it vanishes for $z=\pm m-b$, if and only if
\begin{equation}\label{Eigen Eq Point Step}
\begin{aligned}
&w_{m-b}^{+}+w_{m-b}^{-}+4a m\, w_{m-b}^{+}w_{m-b}^{-}=0 \qquad &&\text{when }z=m-b\text{ and}\\
&w_{-m-b}^{+}+w_{-m-b}^{-}+4a m =0 \qquad &&\text{when }z=-m-b.
\end{aligned}
\end{equation}
Note that, no points in $\{z\in \C: |\Re z|\geq m, |\Im z|=1\}$ can be an eigenvalue. 

\begin{proof}[Proof of Proposition~\ref{Prop Steplike}]
Let $m=0$. Then, directly from the definition it follows that $w_{z}^{0}$ takes one of two values, $\pm i$, depending on the location of $z\in \C\setminus \{\pm m -b\}$. Moreover, substitution shows that, $|\Im z|<1$, iff $w_{z}^{+}=i$ and $w_{z}^{-}=-i$. Also, we directly see that $z=m-b$ or $z=-m-b$ are eigenvalues, if $|\Im b|<1$. This gives the first conclusion of Proposition~\ref{Prop Steplike}.

Let $m>0$. We show the existence of infinitely many large real eigenvalues. Since, $(w_{z}^{0}+w_{z}^{+})(w_{z}^{0}+w_{z}^{-})\not=0$ whenever $z\in \R$, then \eqref{Eigen Eq Step} reduces to
\[ e^{i4a \sqrt{(z+b)^2-m^2}}=\frac{(w_{z}^{0}-w_{z}^{+})(w_{z}^{0}-w_{z}^{-})}{(w_{z}^{0}+w_{z}^{+})(w_{z}^{0}+w_{z}^{-})}.\]
Now, $\cot(w)=i \frac{e^{i2w}-1}{e^{i2w}+1}$ for $w\in \R$ and note that $w_{z}^{+}+w_{z}^{-}\neq 0$ for all $z\in \C$. Thus,  we can rewrite this equation as
\[
\cot(2a\sqrt{(z+b)^2-m^2})=-i\frac{(w_{z}^{0})^2+w_{z}^{+}w_{z}^{-}}{w_{z}^{0}(w_{z}^{+}+w_{z}^{-})}.
\]
The right-hand side of this expression is real for real $z$ with sufficiently large modulus. Indeed,
\[ -i\frac{(w_{z}^{0})^2+w_{z}^{+}w_{z}^{-}}{w_{z}^{0}(w_{z}^{+}+w_{z}^{-})}=\left\{\begin{aligned}
&S(z) \qquad &&\text{ if }z>m-b,\\
&-S(z)  \qquad &&\text{ if }z<-m-b,
\end{aligned} \right.\]
where
\[ S(z)=\sqrt{\frac{z+m+b}{z-m+b}}\frac{(z-m+b)\sqrt{(z+m)^2+1}-(z+m+b)\sqrt{(z-m)^2+1}}{2(z+m+b)\Re\left(\sqrt{z+m+i}\sqrt{m+i-z} \right)}.\]
Since, 
\[ \lim_{z\to+\infty} S(z) = \lim_{z\to-\infty} S(z)=b, \]
by periodicity, we have infinitely many positive and negative solutions. This is the second statement in Proposition~\ref{Prop Steplike}, and completes its proof. 
\end{proof}

\subsection{The weakly coupled model}\label{Subsec Improved dynamics}
At the end of this subsection we give the proof of Theorem~\ref{Theo Weak Coupling}, about localisation of the discrete spectrum for a potential of the form $\epsilon V$ in the asymptotic regime $|\epsilon| \to 0$, 
where $V\in L^1(\mathbb{R},\mathbb{C}^{2\times 2};\mathrm{d}\nu_1)$ satisfies additional hypotheses.

Recall  \eqref{Weak Couple 0},  \[ 
\operatorname{Spec}_{\mathrm{dis}} (\mathscr{L}_{m,\epsilon V}) \subseteq \left\{z \in \C: |\Im z|<1,\, |\Re z| >  \frac{C(m)}{|\epsilon|^{\frac12}\|V\|_{L^1}^{\frac12}}\right\} .
\
\]
Hence, for small enough $|\epsilon|$, the eigenvalues of $\mathscr{L}_{m,\epsilon V}$ lie in the instability band \[\Sigma \coloneqq  \{z\in \C: |\Im z|<1\},\] and they escape to infinity in the regime $|\epsilon|\to 0$ at a rate proportional to $|\epsilon|^{-\frac12}$ or faster. We now formulate a more precise statement about this, for $V$ satisfying a faster decay rate at infinity. The next result is a precursor of Theorem~\ref{Theo Weak Coupling} and it is one of the main contributions of this paper. We recall that the expression of the densities $\nu'_k$ are given in Section~\ref{Sec Main Results} in the paragraph above Theorem~\ref{Theo Weak Coupling}.

\begin{theorem}\label{Theo Weak Coupling 2}
Let $m>0$ and let $V\in L^1(\R,\C^{2\times 2}; \dd \nu_2)$. Let 
\begin{equation}\label{Upsilon}
\Upsilon_{z}\coloneqq \frac{1}{2(4|z|^2+m^2)}\begin{pmatrix}
-(2z+m)^2 && i( 4z^2-m^2)\\
i( 4z^2-m^2) && (2z-m)^2
\end{pmatrix}\in \mathbb{C}
\end{equation}
for all $z\in\mathbb{C}$. Then,
\begin{equation}\label{Weak Coupling 1}
\operatorname{Spec}_{\mathrm{dis}}\left(\mathscr{L}_{m,\epsilon V} \right)\subset\left\{z \in \C: |\Im z|< 1,\, \frac{4m}{4|z|^2+m^2} = W_z(\epsilon)\right\},
\end{equation}
where
\[
 W_z(\epsilon)=\epsilon \int_{\R} e^{-2\left(|x|+i\left(z-\frac{m^2}{2 \Re z}\right)x\right)}\langle V(x), \overline{\Upsilon_{z}} \rangle_{\mathrm{F}}\, \dd x +\mathcal{O}(\epsilon^2)
\]
in the regime $|\epsilon|\to 0$.
\end{theorem}

The verification of the validity of this statement occupies most of this subsection. The explicit dependence on $z$ of $W_z(\epsilon)$ and the constant in the limit, will be established below. The proof of \eqref{Weak Coupling 1}  will invoke Theorem~\ref{Theo BS} once again. Note that the Birman-Schwinger operator associated to $\epsilon V$ is $\epsilon Q(z)$. We test whether $z\in \Sigma$ does not belong to $\operatorname{Spec}_{\mathrm{p}}(\mathscr{L}_{m,\epsilon V})$, by testing whether $-1\not\in \operatorname{Spec}_{\mathrm{p}}(\epsilon  Q(z))$.  However, for this we will not appeal directly to estimates for the norm $\|\epsilon Q(z)\|$ as we did in Subsection~\ref{Subsec Eigen estimate m>0}. Instead, we will find a criterion, involving $W_z(\epsilon)$, to test whether the operator $\epsilon Q(z)+1$ is bijective. Since $Q(z)$ is compact (see the paragraph above Remark~\ref{Remark Compact Q}), this is a test on whether $-1\notin \operatorname{Spec}_{\mathrm{p}}(\epsilon  Q(z))$. 

The position of the discrete spectrum is  associated with sub-regions of $\Sigma$ where the norm of the resolvent is singular. In Subsection~\ref{Subsec Sharp Estimate}, we isolated the singular part of the unperturbed kernel through the term $\mathcal{R}_{1,z}$. Since the latter still carries the leading order information about the singularities of the perturbed kernel (under small enough perturbations), we now consider a decomposition of the kernel of the integral operator associated to $Q(z)$ that takes the singular contribution of $\mathcal{R}_{1,z}$ into account. For this we proceed as follows.  

Let
\[ Q(z) = L(z) +M(z),\]
where 
\begin{gather}
\label{expression for L_z} L(z)f(x)\coloneqq  \int_{\R} A(x)\mathcal{L}_{z}(x,y)B(y) f(y)\, \dd y, \qquad \mathcal{L}_{z}(x,y)\coloneqq  e^{-\eta_{z}(x)-\eta_{z}(y)} U_{z}, \\ \nonumber
\eta_{z}(w)\coloneqq  i \left( z-\frac{m^2}{2\Re z} \right)  w+|w| \qquad\text{and} \qquad U_{z} \coloneqq  \frac{1}{8m}\begin{pmatrix}
-(2z+m)^2 && i(4z^2-m^2)\\
i(4z^2-m^2) &&(2z-m)^2
\end{pmatrix}.
\end{gather}
Then, as we shall see below, $L(z)$ is a rank-one operator for all $z\in \Sigma$ and it carries the leading order contribution of the norm of $Q(z)$ in the asymptotic regime $|\Re z|\to \infty$.
Notably, the kernel $\mathcal{L}_{z}$ encodes the contribution of the singular part of $\mathcal{R}_{1,z}$ in this regime, while $\|M(z)\|$ is uniformly bounded with respect to $z$. Before proving all this, we motivate this decomposition and determine the asymptotic coefficients of the different terms involved in the argument. Recall the notation $z=\tau +i\delta$ for $\tau, \delta\in \R$ with $|\delta|<1$. 

On the one hand, we note that the expression of the powers $\eta_{z}(w)$, is motivated by the fact that
\[ -\eta_{z}(w)=\left\{\begin{aligned}
\eta_{z}^{-}w\qquad \text{if }w\leq 0,\\
-\eta_{z}^{+}w\qquad \text{if }w\geq 0,
\end{aligned} \right.\qquad \eta_{z}^{-}\coloneqq  (1+\delta)-i\left(\tau-\frac{m^2}{2\tau} \right),\qquad \eta_{z}^{+}\coloneqq  (1-\delta)+i\left(\tau-\frac{m^2}{2\tau} \right),\]
where $\eta_{z}^{\pm}$ are the leading order coefficients in the expansion of the terms $\mu_{z}^{\pm}$:
\begin{equation}\label{mu-eta}
\mu_{z}^{\pm}= \eta_{z}^{\pm} + \frac{m^2(1\mp \delta)}{2 \tau^2}+ \mathcal{O}\left(\frac{1}{|\tau|^3}\right)
\end{equation}
as $|\tau| \to +\infty$. See Lemma~\ref{Lem Asymptotic} and Table~\ref{Table coef}.

On the other hand, the expression of the matrix $U_{z}$, has a more involved motivation which is given after the following lemma.
 
\begin{lemma} \label{N-U^0}
Let the matrices appearing in the expression of the kernel $\mathcal{R}_{1,z}$ in \eqref{R1} be $N_{j,z}$, as defined at the beginning of Subsection~\ref{Subsec Resolvent}. Let 
\begin{equation*}
U^{0}_{\tau+i\delta} \coloneqq 
\begin{pmatrix}
-\frac{1}{2m}\tau^2-\left(\frac{1}{2}+\frac{i\delta}{m}\right) \tau &&  \frac{i}{2m}\tau^2-\frac{\delta}{m}\tau\\
\frac{i}{2m}\tau^2-\frac{\delta}{m}\tau && \frac{1}{2m}\tau^2-\left(\frac{1}{2}-\frac{i\delta}{m}\right) \tau
\end{pmatrix}.
\end{equation*}
Then, $N_{j,z}= U^{0}_{z} + \widehat{N}_{j,z}$  where 
\[
 \vert \widehat{N}_{j,\tau+i \delta} \vert_{\mathrm{F}}=\mathcal{O}(1)
\]
as $|\tau|\to \infty$ for all $j\in \{1,3,4,7,8,10\}$.
\end{lemma}
\begin{proof}
 By writing
\[
\frac{1}{w_{z}^{+}+w_{z}^{-}}=
\frac{-i}{4m}\left[(m+i+z)\mu_{z}^{+}-(m-i+z)\mu_{z}^{-} \right]
\]
we get that
\begin{align*}
\frac{1}{w_{z}^{+}+w_{z}^{-}}&=\frac{1}{2m}\tau^2 + \left(\frac{1}{2}+\frac{i\delta}{m} \right)\tau + \left(\frac{1-\delta^2}{2m}-\frac{m}{4}+\frac{i\delta}{2} \right)+\mathcal{O}\left(\frac{1}{\tau}\right)
\end{align*}
as $|\tau|\to 0$.
The claimed asymptotic formula is achieved by substituting  the expansions of $w_{z}^{\pm}$, $\frac{1}{w_{z}^{\pm}}$ and this expression, into the formulas of the different $N_{j,z}$. See Lemma~\ref{Lem Asymptotic}.
\end{proof}

Now we return to the motivation for the expression of $U_{z}$. We pick the latter, so that it approximates all the matrices $N_{j,z}$ up to order $\tau$. Indeed, writing $U_{z}$ in terms of $\tau$ and $\delta$, gives
\begin{equation}\label{U-U^0}
U_{\tau+i\delta}=U^{0}_{\tau+i\delta}+U^1_{\delta},
\end{equation}
for $U^0_z$ as in Lemma~\ref{N-U^0} and
\begin{align*}
U^1_{\delta}&\coloneqq  \begin{pmatrix}
\frac{\delta^2}{2m}-\frac{m}{8}-\frac{i\delta}{2} & -i\left(\frac{\delta^2}{2m}+\frac{m}{8}\right)\\
-i\left(\frac{\delta^2}{2m}+\frac{m}{8}\right) & -\frac{\delta^2}{2m}+\frac{m}{8}-\frac{i\delta}{2} 
\end{pmatrix}.
\end{align*}
Note that for all $z\in \Sigma$, $\det U_{z} =0$ and $U_z$ is of rank one. Then, the matrix norm of $U_{z}$ can be computed explicitly,
\begin{equation}\label{Norm of Phi}
\vert U_{z} \vert_{\C^{2\times 2}}=\vert U_{z} \vert_{\mathrm{F}}=\frac{|z|^2}{m}+\frac{m}{4}.
\end{equation}

We are now ready to formulate the crucial step in the proof of Theorem~\ref{Theo Weak Coupling 2}.

\begin{lemma}\label{Lem LM}
Let $V\in L^1(\R, \C^{2\times 2},\dd \nu_2)$. For  $A(x)$ and $B(y)$ as in \eqref{Decompose Matrix}, let $L(z)$ be the integral operator \eqref{expression for L_z} and let $M(z)=Q(z)-L(z)$. Let $\Upsilon_{z}=\frac{1}{|U_z|_{\mathrm{F}}}U_z$ be as in \eqref{Upsilon}. Then 
\[
L(z)f = \left(\frac{|z|^2}{m}+\frac{m}{4}\right) \langle f, \overline{\psi_{z}} \rangle\, \phi_{z},
\]
where
\begin{align*}
\psi_{z}(w) & \coloneqq  \left\{
\begin{aligned}
&e^{-\eta_{z}(w)} B(w)^{T} \begin{pmatrix}
(\Upsilon_{z})_{11}\\
(\Upsilon_{z})_{12}
\end{pmatrix}\qquad &&\text{if } z\neq \frac{-m}{2},\\
&e^{-\eta_{z}(w)} B(w)^{T} \begin{pmatrix}
(\Upsilon_{z})_{21}\\
(\Upsilon_{z})_{22}
\end{pmatrix} \qquad &&\text{if } z= \frac{-m}{2},
\end{aligned}
\right.\\
 \phi_{z}(w) & \coloneqq  
\left\{
\begin{aligned}
&e^{-\eta_{z}(w)} A(w) \begin{pmatrix}
1\\ (\Upsilon_{z})_{21}/(\Upsilon_{z})_{11}
\end{pmatrix} \qquad&&\text{if } z\neq \frac{-m}{2},\\
&e^{-\eta_{z}(w)} A(w) \begin{pmatrix}
(\Upsilon_{z})_{12}/(\Upsilon_{z})_{22}\\1
\end{pmatrix} \qquad &&\text{if } z= \frac{-m}{2}.
\end{aligned}
\right. 
\end{align*}
Moreover, there exists a constant $C(m)>0$ such that $\Vert M_{z} \Vert\leq C(m)$ for all $z \in \Sigma$.
\end{lemma}
\begin{proof}
Assume first that $z\neq \frac{-m}{2}$.
Since $U_{z}$ is singular and non-zero, then 
\[ U_{z}= \begin{pmatrix}
1\\ (U_{z})_{21}/(U_{z})_{11}
\end{pmatrix} \begin{pmatrix}
(U_{z})_{11} & (U_{z})_{12}
\end{pmatrix}.\]
Thus, we can rewrite the action of the operator $L_{z}$ in the  form
\begin{align*}
L_{z}f(x) &= |U_{z}|_{\mathrm{F}}\int_{\R} e^{-\eta_{z}(x)-\eta_{z}(y)} A(x) \begin{pmatrix}
1\\ (\Upsilon_{z})_{21}/(\Upsilon_{z})_{11}
\end{pmatrix} \begin{pmatrix}
(\Upsilon_{z})_{11} & (\Upsilon_{z})_{12}
\end{pmatrix} B(y)f(y)\, \dd y\\
&= |U_{z}|_{\mathrm{F}} \left( \int_{\R} e^{-\eta_{z}(y)} \begin{pmatrix}
(\Upsilon_{z})_{11} & (\Upsilon_{z})_{12}
\end{pmatrix} B(y)f(y)\, \dd y \right)\phi_{z}(x)\\
&=|U_{z}|_{\mathrm{F}} \left( \int_{\R} \langle f(y), \overline{\psi_{z}(y)} \rangle_{\C^{2}}\, \dd y \right)\phi_{z}(x)
\end{align*}
This is the claimed expression for the operator $L_{z}$. 

Both vector-valued functions $\psi_{z}$ and $\phi_{z}$ belong to $ L^{2}(\R,\C^2)$  because $\Re \eta_{z}^{\pm}=1\mp \delta>0$ and thus, for all $z\in \Sigma\setminus\left\{-\frac{m}{2}\right\}$,
\begin{equation}\label{Norm two vectors}
\Vert \psi_{z} \Vert \leq \Vert V \Vert_{L^1}^{1/2}\sqrt{(\Upsilon_{z})_{11}^2+(\Upsilon_{z})_{12}^2},\qquad \Vert \phi_{z} \Vert \leq \Vert V \Vert_{L^1}^{1/2}\sqrt{1+\frac{(\Upsilon_{z})_{21}^2}{(\Upsilon_{z})_{11}^2}}.
\end{equation}
Therefore, we can conclude that, for all $z\in \Sigma\setminus\left\{-\frac{m}{2}\right\}$, $L_{z}$ is an operator of rank-one and
\begin{equation}\label{Upper bound of L}
\Vert L_{z} \Vert = |U_{z}|_{\mathrm{F}}\Vert \psi_{z} \Vert \Vert \phi_{z} \Vert\leq |U_{z}|_{\mathrm{F}} \Vert V \Vert_{L^1}.
\end{equation}
Now,  to achieve the same conclusion for $z=-\frac{m}{2}$, we can just write the matrix $U_z$ in different form, as
\[ U_{z}= \begin{pmatrix}
(U_{z})_{12}/(U_{z})_{22}\\
1 
\end{pmatrix} \begin{pmatrix}
(U_{z})_{21} & (U_{z})_{22}
\end{pmatrix},\]
and follow the same argument as above. 

\medskip

Next, let us complete the proof of the lemma by showing the validity of the claim made about $M_{z}$. Let \begin{align*}
\mathcal{M}_{z}(x,y) \coloneqq  \mathcal{R}_{z}(x,y)- \mathcal{L}_{z}(x,y).
\end{align*}
Then,
\begin{equation}\label{Kernel M_E}
\mathcal{M}_{z}(x,y) = \left(\mathcal{R}_{1,z}(x,y)- \mathcal{L}_{z}(x,y)\right) + \mathcal{R}_{2,z}(x,y).
\end{equation}
By virtue of \eqref{R2}, \eqref{nj} and Lemma~\ref{Lem Asymptotic}, it follows that there exists a constant $c_1>0$ such that
\begin{equation}\label{Bound R2}
\left\vert \mathcal{R}_{2,z} \right\vert_{\C^{2\times 2}}\leq c_1,
\end{equation}
for all $|\tau|\geq 1$.

Let us now show that the difference of $\mathcal{R}_{1,z}-\mathcal{L}_{z}$ is also uniformly bounded. Let
\[ \mu_{z}(x)\coloneqq \left\{\begin{aligned}
&\mu_{z}^{-} x \qquad &\text{if } x\leq 0\\
-&\mu_{z}^{+} x \qquad &\text{if } x\geq 0
\end{aligned}\right..\]
Then
\[  \mathcal{R}_{1,z}(x,y)= N_{j,z} e^{\mu_{z}(x)+\mu_{z}(y)}\]
for $j$ chosen depending on the signs of $x$ and $y$ as in \eqref{R1}.
For all $(x,y)\in\mathbb{R}^2$, write
$$\Lambda_{z}(x,y)\coloneqq  \mu_{z}(x)+\eta_{z}(x) + \mu_{z}(y)+\eta_{z}(y).$$  Then, 
\begin{equation*}
 \mathcal{R}_{1,z}(x,y)-\mathcal{L}_{z}(x,y) =U_{z} \left( e^{\Lambda_{z}(x,y)}-1\right)e^{-\eta_{z}(x)-\eta_{z}(y)}+\left( N_{j,z}- U_{z}\right)e^{\mu_{z}(x)+\mu_{z}(y)}.
\end{equation*}
Since $\Re \Lambda_{z}(x,y)<0$ for all $(x,y)\in \R^2$ and by virtue of \eqref{mu-eta}, considering $|\tau|$ sufficiently large so that $\Re\left(\mu_{z}^{+}+\eta_{z}^{+}\right)>0$ and $\Re\left(\mu_{z}^{-}+\eta_{z}^{-}\right)>0$, gives 
\begin{align*}
\left\vert e^{\Lambda_{z}(x,y)} -1  \right\vert \leq \vert \Lambda_{z}(x,y)\vert \leq \frac{c_2}{\tau^2} \left( |x| + |y|\right).
\end{align*}
Therefore,
since $\Re \mu_{z}(x) <0$ and $\Re \eta_{z}(x)>0$ for all $x\in \R$, according to Lemma~\ref{N-U^0} and to \eqref{U-U^0}, we have
\begin{equation}\label{Bound R1-L}
\left\vert \mathcal{R}_{1,z}(x,y)-\mathcal{L}_{z}(x,y)\right\vert_{\C^{2\times 2}}
\leq   \vert U_{z} \vert_{\C^{2\times 2}}\left\vert e^{\Lambda_{z}(x,y)} -1 \right\vert+ \vert N_{j,z}-U_{z} \vert_{\C^{2\times 2}}\leq c_3( |x|+|y|+1),
\end{equation}
for all $|\tau|\geq \tau_0$.

From \eqref{Bound R2} and \eqref{Bound R1-L} it follows that, by taking the Hilbert-Schmidt norm in the decomposition \eqref{Kernel M_E},  
\begin{align*}
\Vert M(z) \Vert^2 &\leq c_4 \int_{\R^2} \vert V(x)\vert_{\C^{2\times 2}}(1+|x|+|y|)^2  \vert V(y)\vert_{\C^{2\times 2}} \dd x \, \dd y\\
&\leq c_5 \left(\int_{\R} \vert V(x)\vert_{\C^{2\times 2}} \nu_2'(x)\, \dd x\right)\left(\int_{\R} \vert V(y)\vert_{\C^{2\times 2}} \nu_2'(y)\, \dd y\right),
\end{align*}
where the right-hand side is finite and independent of $z$, for all $|\tau|\geq \tau_0$. On the other hand, according to \eqref{Upper bound of L} for $L(z)$, and  \eqref{Bound out W}-\eqref{Bound on W} for $Q(z)$, the triangle inequality yields,
\begin{align*}
\Vert M(z) \Vert &\leq \Vert Q(z) \Vert + \Vert L(z) \Vert\leq  c_6 \int_{\R} |V(x)|_{\C^{2\times 2}} \nu_1'(x)\, \dd x<\infty
\end{align*}
for all $|\tau|\leq \tau_0$, where the right-hand side is also independent of $z$.
This completes the proof of the lemma.
\end{proof}

With this lemma at hand, we complete the proof of Theorem~\ref{Theo Weak Coupling 2} as follows. 

\begin{proof}[Proof of Theorem~\ref{Theo Weak Coupling 2}] Since $\|M(z)\|$ is uniformly bounded for all $z\in \Sigma$, the operator $I+\epsilon M(z)$ is  invertible whenever 
	\[
		0 < \epsilon <\frac{1}{C(m)}. 
	\]
Assume this from now on. Then,
since
$$I+\epsilon Q(z)= I +\epsilon(L(z)+M(z)) = (I+\epsilon M(z))\big( I+\epsilon (I+\epsilon M(z))^{-1} L(z)\big),$$
we have
\begin{equation}\label{Equivalence 1}
z \in \operatorname{Spec}_\text{{\normalfont p}}(\mathscr{L}_{m,\epsilon V}) \iff -1 \in \operatorname{Spec}(\epsilon(I+\epsilon M(z))^{-1}L(z)).
\end{equation}

By virtue of Lemma~\ref{Lem LM}, $L(z)$ is a rank-one operator. Therefore, $\epsilon(I+\epsilon M(z))^{-1}L(z)$ is also a rank-one operator, such that
\[\epsilon(I+\epsilon M(z))^{-1}L(z) f = \epsilon |U_z|_{\mathrm{F}} \langle f, \overline{\psi_{z}} \rangle (I+\epsilon M(z))^{-1} \phi_{z} .\]
Hence, $\epsilon(I+\epsilon M(z))^{-1}L(z)$, being a compact operator on the infinite-dimension Hilbert space $L^{2}(\R,\C^2)$, has spectrum
\[ \operatorname{Spec}\left( \epsilon(I+\epsilon M(z))^{-1}L(z)\right)=\left\{0, \epsilon |U_z|_{\mathrm{F}}\left\langle (I+\epsilon M(z))^{-1} \phi_{z}, \overline{\psi_{z}} \right\rangle \right\}.\]
Thus, from \eqref{Equivalence 1} and by writing $ (I+\epsilon M(z))^{-1}=I-\epsilon M(z)(I+\epsilon M(z))^{-1}$, we have that
\begin{equation}\label{Equivalence 2}
\begin{aligned}
z\in \operatorname{Spec}_{\mathrm{p}}(\mathscr{L}_{m,\epsilon V}) &\iff -1 =\epsilon |U_{z}|_{\mathrm{F}}\left\langle (I+\epsilon M(z))^{-1} \phi_{z}, \overline{\psi_{z}} \right\rangle\\ &\iff \frac{1}{|U_{z}|_{\mathrm{F}}} = -\epsilon  \left\langle \phi_{z}, \overline{\psi_{z}} \right\rangle + \epsilon^2 \left\langle M(z)(I+\epsilon M(z))^{-1}\phi_{z}, \overline{\psi_{z}} \right\rangle \\ &\iff \frac{4m}{4|z|^2+m^2}=W_z(\epsilon) \end{aligned}
\end{equation}
where $W_z(\epsilon)$ is the right-hand side of the second line in this expression. This is \eqref{Weak Coupling 1}.

We complete the proof by confirming the claimed asymptotic formula for $W_z(\epsilon)$. Notice that, \begin{equation} \label{inner product expression}
\left\langle \phi_{z}, \overline{\psi_{z}} \right\rangle=-\int_{\R} e^{-2 \eta_{z}(x)} \langle V(x), \overline{\Upsilon_{z}} \rangle_{\mathrm{F}}\, \dd x,
\end{equation}
because $V(x)=B(x)A(x)$ and $(\Upsilon_{z})_{12}= (\Upsilon_{z})_{21}$. Now, \[ \left\vert  \left\langle M(z)(I+\epsilon M(z))^{-1}\phi_{z}, \overline{\psi_{z}} \right\rangle \right\vert\leq \frac{\Vert M(z) \Vert}{1-\epsilon \Vert M(z) \Vert} \Vert \phi_{z} \Vert \Vert \psi_{z} \Vert \leq \frac{\Vert M_{z} \Vert \Vert V \Vert_{L^1}}{1-\epsilon \Vert M_{z} \Vert}\leq c_7.\]
This is a consequence of \eqref{Norm two vectors} and the fact that $\Vert M(z) \Vert$ is uniformly bounded on $\Sigma$, and it completes the proof of the theorem.\end{proof}

\begin{proof}[Proof of Theorem~\ref{Theo Weak Coupling}] As setting in \eqref{Equivalence 2}, we write
\[
     W_z(\epsilon)=\epsilon a_z+\epsilon^2 b_{z}(\epsilon)
\]
where $a_z=-\left\langle \phi_{z}, \overline{\psi_{z}} \right\rangle$ is independent of $\epsilon$ and $|b_z(\epsilon)| \leq c_1$ for all $z\in\Sigma$ and sufficiently small $|\epsilon|$.

We aim to show that there exists a constant $c_2>0$ such that
\begin{equation} \label{bound bz}
    |a_{\tau+i\delta}|\leq \frac{c_2}{|U_{\tau+i\delta }|_{\mathrm{F}}}
\end{equation}
for all $|\tau|\geq 1$. To establish this, we first use integration by parts twice to obtain the following results:
\begin{align*}
\int_{\R} e^{-2 \eta_{z}(x)} v(x) \, \dd x = & \left(\frac{1}{4(\eta_{z}^{+})^2}-\frac{1}{4(\eta_{z}^{-})^2}\right)v'(0)-\left(\frac{1}{2\eta_{z}^{-}}+\frac{1}{2\eta_{z}^{+}}\right)v(0)\\
&+ \frac{1}{4 (\eta_{z}^{-})^2} \int_{-\infty}^{0} e^{2\eta_{z}^{-}x}v''(x)\, \dd x +\frac{1}{4 (\eta_{z}^{+})^2} \int_{0}^{+\infty} e^{-2\eta_{z}^{+}x}v''(x)\, \dd x,
\end{align*}
for any $v \in W^{2,1}(\R)$. Now, from the expressions for $\eta_{z}^{\pm}$, it follows that for suitable constants $c_k>0$,
\[ \left\vert \frac{1}{2\eta_{z}^{-}}+\frac{1}{2\eta_{z}^{+}}\right\vert=\frac{1}{|\eta_{z}^{-}\eta_{z}^{+}|}\leq \frac{c_3}{\tau^2}\leq \frac{c_4}{|U_{z}|_{\mathrm{F}}}\qquad \text{and} \qquad \frac{1}{\vert \eta_{z}^{\pm} \vert^2}\leq \frac{c_5}{|U_{z}|_{\mathrm{F}}}\]
for all $|\tau|\geq 1$. Since the Frobenius norm of $\Upsilon_z$ is one, setting $v(x)\coloneqq \langle V(x), \overline{\Upsilon_{z}} \rangle_{\mathrm{F}}$ gives \[ \vert v^{(j)}(x)\vert \leq  \vert V^{(j)}(x) \vert_{\mathrm{F}},\]
for $j\in \{0,1,2\}$ and all $x\in \R$. Hence, according to \eqref{inner product expression}, we have that \eqref{bound bz} is indeed valid under the hypothesis $V \in W^{2,1}(\R,\C^{2\times 2})$.
 
According to Theorem \ref{Theo Weak Coupling 2}, $z$ belongs to $\operatorname{Spec}_{\mathrm{dis}}(\mathscr{L}_{m,\epsilon V})$ if and only if 
\[
    \frac{1}{|U_z|_{\mathrm{F}}}=W_z(\epsilon).
\]
To satisfy this condition, $z$ must be in $\Sigma$ such that
\[
    \left|\frac{1}{|U_z|_{\mathrm{F}}}-\epsilon a_z \right|=\epsilon^2 |b_{z}(\epsilon)|.
\] 
Then, $z$ should meet the requirement
\[
    \frac{1}{|U_z|_{\mathrm{F}}}\left(1-\epsilon c_2\right)\leq \epsilon^2 c_1.
\]
Therefore, there must exists a constant $c_6>0$ such that
\[
     \frac{1}{|U_z|_{\mathrm{F}}}\leq c_6 \epsilon^2
\]
for all sufficiently small $|\epsilon|$. Consequently, for $z$ to be an eigenvalue of $\mathscr{L}_{m,\epsilon V}$ when $|\epsilon|$ is small, we require
\[
    |\Re z|^2\geq c_7 |U_z|_{\mathrm{F}} > c_8 \epsilon^{-2}.
\]
This condition ensures the validity of Theorem~\ref{Theo Weak Coupling}.
\end{proof}

\subsection{Proof of Theorem~\ref{inclusion pseudospectrum}\label{Subsec pseudospectrum V}}
Assume that $V$ satisfies the condition \ref{Ass norm less than 1}.
 According to Proposition~\ref{Prop Perturbation}, we know that $Q(z)=\overline{AR_0(z)B}$ is well defined for all $z\in \rho (\mathscr{L}_{m})$. Moreover, by virtue of \eqref{Hilbert-Schmidt norm}, 
\[
    \|Q(z)\|\leq \|V\|_{L^1} \frac{(|k_z|+1) \max\left\{|w_z^{\pm}|+\frac{1}{|w_z^{\pm}|}\right\}}{2}
\]
for all such $z$. Recall the set $U$ defined in Subsection \ref{Subsec Eigen estimate m>0} and recall that \eqref{For the estimate on U} is valid for all $z\in U$. Then, for $\|V\|_{L^1}<\frac{2}{c_{14}(m)}=:C(m)$, $\|Q(z)\|<1$ and so
\[
     \Vert (1+Q(z))^{-1} \Vert \leq \frac{C(m)}{C(m)-\|V\|_{L^1}},
\]
for all $z\in U$. Now, from \eqref{bound AR} and \eqref{bound RB}, we have
\[
    \max \left\{\|\overline{R_0(z)B}\|,\,\|AR_0(z)\|\right\} \leq \|V\|_{L^1}^{\frac12}
\frac{(|k_z|+1) \max\left\{|w_z^{\pm}|+\frac{1}{|w_z^{\pm}|}\right\}}{2\big(\min \Re \mu_z^{\pm}\big)^{\frac12}},
\] 
for all $z\in \rho (\mathscr{L}_m)$. 

Consider the minimum in the denominator of the right-hand side. Since $\Re \sqrt{\gamma}=\frac{1}{\sqrt{2}}\sqrt{|\gamma|+\Re \gamma},$
we have that
\[ \Re \mu_{z}^{\pm} = \frac{1}{\sqrt{2}}\sqrt{\sqrt{(\tau^2-m^2-(1\mp\delta)^2)^2+4(1\mp \delta)^2 \tau^2}+m^2+(1\mp \delta)^2-\tau^2}.\]
Hence, since $\tau^2\geq \tau^2-m^2$, then 
\begin{align*}
(\tau^2-m^2-(1\mp\delta)^2)^2+4(1\mp \delta)^2 \tau^2 &\geq (\tau^2-m^2-(1\mp\delta)^2)^2+4(1\mp \delta)^2 (\tau^2-m^2)\\
&=(\tau^2-m^2+(1\mp\delta)^2)^2.
\end{align*}
Thus, $\Re \mu_{z}^{\pm} \geq |1\mp \delta|$ for all $z\in \C$ and so $\min \Re \mu_{z}^{\pm} \geq |\delta|-1$ for all $|\delta|>1$. 

Gathering the estimates above, we then have that, if $\|V\|_{L^1}\leq C(m)$,
\begin{align*}
    \big\|(\mathscr{L}_{m,V}-z)^{-1}-(\mathscr{L}_{m}-z)^{-1}\big\|&\leq
\|\overline{R_0(z)B}\|\,\|(I+Q(z))^{-1}\|\,\|AR_0(z)\| \\
&\leq \frac{\|V\|_{L^1}}{(|\Im z|-1)C(m)(C(m)-\|V\|_{L^1})},
\end{align*}
for all $z\in U$. Moreover, since $U$ lies outside the numerical range of $\mathscr{L}_m$, then
\[
    \|(\mathscr{L}_{m}-z)^{-1}\|\leq \frac{1}{|\Im z|-1},
\]
for all $z\in U$. Therefore, according to the triangle inequality, whenever $\|V\|_{L^1}<C(m)$, we have that
\[
    \big\|(\mathscr{L}_{m,V}-z)^{-1}\big\|\leq \left(1+\frac{\|V\|_{L^1}}{C(m)(C(m)-\|V\|_{L^1})}\right)\frac{1}{|\Im z|-1}, 
\]
for all $z\in U$. This gives the claim made in Theorem~\ref{inclusion pseudospectrum}.


\appendix

%


\section{Abstract Birman-Schwinger principle}\label{Appendix Abstract Birman Schwinger}

In order to characterise the domain of the operator $\mathscr{L}_{m,V}$ in terms of the resolvent of $\mathscr{L}_{m}$, we follows the classical approach of \cite{Kato66}, which was extended to the non-selfadjoint regime in \cite{Gesztesy-Latushkin-Mitrea-Zinchenko20}. We include precise details of this construction in this appendix, as the abstract framework might be applicable to other families of operators related to perturbations of $\mathscr{D}_m$.

Let $\mathcal{H}$ and $\mathcal{K}$ be two Hilbert spaces. Let $\mathrm{Dom}(H_{0})\subset \mathcal{H},\,\mathrm{Dom}(A)\subset\mathcal{H}$ and $\mathrm{Dom}(B)\subset \mathcal{K}$ be three dense subspaces. Let
\[ H_{0}:\mathrm{Dom}(H_{0})\longrightarrow \mathcal{H},\qquad A: \mathrm{Dom}(A)\longrightarrow \mathcal{K}\quad \text{and}\quad B: \mathrm{Dom}(B)\longrightarrow \mathcal{H}, \]
be three closed operators, such that $H_{0}$ has a non-empty resolvent set. Let
\[R_{0}(z) \coloneqq  (H_{0}-z I_{\mathcal{H}})^{-1}, \qquad z\in \rho(H_{0}). \]
We consider the following hypotheses.

\begin{enumerate}[label=\textup{\textbf{(B\arabic*)}}]
\item \label{Bounded Exten 1} For some (and hence all) $z\in \rho(H_{0})$, \begin{itemize} \item the operator $R_{0}(z)B$ is closable, the closure $\overline{R_{0}(z)B}:\mathcal{K}\longrightarrow \mathcal{H}$ is bounded, \item and the operator $AR_{0}(z):\mathcal{H}\longrightarrow \mathcal{K}$ is bounded. \end{itemize}
\item \label{Bounded Exten 2} For some (and hence all) $z\in \rho(H_{0})$, the operator $AR_{0}(z)B$ has a bounded closure,
\[ Q(z)\coloneqq  \overline{AR_{0}(z)B}: \mathcal{K}\longrightarrow \mathcal{K}.\]
\item \label{Nonempty} The set \[\mathsf{R}\coloneqq \Big\{z\in \rho(H_{0}): -1\in \rho(Q(z))\Big\}\]
is non-empty.
\end{enumerate}
For all $z\in \mathsf{R}$, we let 
\begin{equation}\label{R}
R(z)= R_{0}(z)-\overline{R_{0}(z)B}\big( I_{\mathcal{K}}+Q(z)\big)^{-1}AR_{0}(z): \mathcal{H}\longrightarrow \mathcal{H}.
\end{equation}

\begin{theorem}\label{Theo Closed Exten}
Assume that the closed operators $H_0,\,A$ and $B$ are related in the above manner, and that \ref{Bounded Exten 1}--\ref{Nonempty} hold true. Then, there exists a closed densely defined extension\footnote{Here $\operatorname{Dom}(BA)=\{f\in\operatorname{Dom}(A)\,:\, Af\in\operatorname{Dom}(B)\}$.} \[H\supseteq H_{0}+BA:\mathrm{Dom}(H_{0})\cap \mathrm{Dom}(BA)\longrightarrow \mathcal{H},\] whose resolvent is given by 
\[ (H-zI_{\mathcal{H}})^{-1}=R(z), \qquad \text{ for all } z\in \mathsf{R}.\]
\end{theorem}

The conclusion of Theorem~\ref{Theo Closed Exten} can be found in \cite[Thm.~2.3]{Gesztesy-Latushkin-Mitrea-Zinchenko20}, assuming the slightly more stringent conditions \cite[Hyp. 2.1 (i)]{Gesztesy-Latushkin-Mitrea-Zinchenko20}; \[\mathrm{Dom}(H_{0})\subset \mathrm{Dom}(A) \quad \text{and} \quad \mathrm{Dom}(H_{0}^{*})\subset \mathrm{Dom}(B^{*}).\] In this appendix we show that the assumption \ref{Bounded Exten 1} also ensures the claimed conclusion.

\begin{proof}
We split the proof into four steps.

\underline{Step 1}. Let $z\in \rho(H_{0})$, then the product $AR_{0}(z)B$ is always closable. Indeed, we claim that $AR_{0}(z)B$ is densely defined and its adjoint is also densely defined, see \cite[Thm. 1.8 (i)]{Schmudgen12}. The former is an immediate consequence of \ref{Bounded Exten 1}, so that 
\begin{equation}\label{Dom B}
\mathrm{Dom}(AR_{0}(z)B) = \mathrm{Dom}(B).
\end{equation}
For the latter, according to \cite[Prop. 1.7]{Schmudgen12}, we know that
\[ \left(AR_{0}(z)B\right)^{*}\supset (R_{0}(z)B)^{*}A^{*}.\] 
Then, by virtue of \cite[Thm. 1.8 (ii)]{Schmudgen12}, $\left(R_{0}(z)B\right)^{*}=\left(\overline{R_{0}(z)B}\right)^{*}$, so by \ref{Bounded Exten 1} we have
\[ \mathrm{Dom}\left((R_{0}(z)B)^{*}A^{*}\right)=\mathrm{Dom}(A^{*}).\]
Since $A$ is closed, $\mathrm{Dom}(A^{*})$ and thus $\mathrm{Dom}(\left(AR_{0}(z)B\right)^{*})$ are dense in $\mathcal{K}$.

\underline{Step 2}. Let $z\in \mathsf{R}$. Our second task is to construct a closed densely defined operator, $H(z)$, such that
\begin{equation}\label{Resolvent}
R(z) = \left(H(z)-zI_{\mathcal{H}}\right)^{-1}.
\end{equation}

For this purpose, firstly observe that
\begin{equation}\label{Express Q}
\overline{R_{0}(z)B}\,(\mathcal{K})\subset \mathrm{Dom}(A) \qquad \text{ and }\qquad Q(z)=A \overline{R_{0}(z)B}.
\end{equation}
Indeed, for $f\in \mathrm{Dom}(B)$ and $g\in \mathrm{Dom}(A^{*})$, we have
\[ \left\langle R_{0}(z)B f,A^{*}g  \right\rangle_{\mathcal{H}}= \left\langle A R_{0}(z)B f,g  \right\rangle_{\mathcal{K}}=\left\langle Q(z) f,g  \right\rangle_{\mathcal{K}}.\]
Here the second equality is valid as a consequence of \eqref{Dom B}. Thus, the continuity of $\overline{R_{0}(z)B}$ and $Q(z)$, combined with the density of $\mathrm{Dom}(B)$ in $\mathcal{K}$, imply
\[ \left\langle \overline{R_{0}(z)B} f, A^{*}g \right\rangle_{\mathcal{H}}=\left\langle Q(z) f,g  \right\rangle_{\mathcal{K}}\]
for all $f \in \mathcal{K}$ and $g\in \mathrm{Dom}(A^{*})$. This gives \eqref{Express Q}. 

Now,  multiplying $R(z)$ by $A$ on the left, according to \eqref{R} and \eqref{Express Q}, we get that
\begin{equation*}
\begin{aligned}
AR(z)&= AR_{0}(z)- Q(z)\left( I_{\mathcal{K}}+Q(z)\right)^{-1}AR_{0}(z) = \left( I_{\mathcal{K}}+Q(z)\right)^{-1}AR_{0}(z).
\end{aligned}
\end{equation*}
This yields the following. If $R(z)f=0$, where $f\in \mathcal{H}$, then $AR_{0}(z)f=0$. Hence, by \eqref{R}, $R_{0}(z)f=0$ and so $f=0$. That is, $R(z)$ is necessarily injective.

The adjoint of $R(z)$ is given by
\begin{equation}\label{Adjoint R}
R(z)^{*}= R_{0}(z)^{*}-\left(AR_{0}(z)\right)^{*}\left( I_{\mathcal{K}}+Q(z)^{*}\right)^{-1}\left(R_{0}(z)B \right)^{*}.
\end{equation}
By virtue of \cite[Thm. 1.8 (ii) \& Prop. 1.7 (ii)]{Schmudgen12}, we gather that 
$ Q(z)^{*}=B^{*}(AR_{0}(z))^{*}$. In conjunction with \eqref{Adjoint R}, this yields
\begin{equation}\label{B*R*}
\begin{aligned}
B^{*}R(z)^{*} &= \left(R_{0}(z)B\right)^{*}-Q(z)^{*}\left( I_{\mathcal{K}}+Q(z)^{*}\right)^{-1}\left(R_{0}(z)B \right)^{*}\\
&=\left( I_{\mathcal{K}}+Q(z)^{*}\right)^{-1}\left(R_{0}(z)B \right)^{*}.
\end{aligned}
\end{equation}
Reasoning out as proving the injection of $R(z)$ as above by using \eqref{Adjoint R} and \eqref{B*R*}, we obtain
  \[ \mathrm{Ker}(R(z)^{*})=\{0\}.\]
Hence $R(z)(\mathcal{K})$ is dense in $\mathcal{H}$. 

For each $z\in \mathsf{R}$, we can then set $H(z)\coloneqq  R(z)^{-1}+z I_{\mathcal{H}}$. By construction, $H(z)$ is a closed densely defined operator satisfying \eqref{Resolvent}. Note that $H(z)$ is ensured to be closed as a consequence of \cite[Thm. 1.8 (vi)]{Schmudgen12}.

\underline{Step 3}. We now show that the operator $H(z)$ from the previous step, does not depend on $z \in \mathsf{R}$.

 Let $z_{1},z_{2} \in \mathsf{R}$. Then, $H(z_{j})$ are such that
\[ R(z_{j})=(H(z_{j})-z_{j}I_{\mathcal{H}})^{-1}.\] 
We claim that $H(z_{1})=H(z_{2})$. To prove this claim, note that $R(z)$ satisfies the resolvent identity
\begin{equation}\label{Second Resolvent}
R(z_{1})-R(z_{2})=(z_{1}-z_{2})R(z_{2})R(z_{1}).
\end{equation}
Indeed, by \eqref{R}, \eqref{RB Resolvent} and \eqref{Q1-Q2}, it follows that
\begin{align*}
(z_{1}-z_{2})R(z_{2})R(z_{1})
=&R_{0}(z_{1})-R_{0}(z_{2})-\left[\overline{R_{0}(z_{1})B}-\overline{R_{0}(z_{2})B} \right]\left( I_{\mathcal{K}}+Q(z_{1})\right)^{-1}AR_{0}(z_{1})\\
&-\overline{R_{0}(z_{2})B}\left( I_{\mathcal{K}}+Q(z_{2})\right)^{-1}A\left[ R_{0}(z_{1})-R_{0}(z_{2})\right]\\
&+ \overline{R_{0}(z_{2})B}\left( I_{\mathcal{K}}+Q(z_{2})\right)^{-1}\left[ Q(z_{1})-Q(z_{2})\right]\left( I_{\mathcal{K}}+Q(z_{1})\right)^{-1}AR_{0}(z_{1})\\
=&R(z_{1})-R(z_{2}) + \overline{R_{0}(z_{2})B}\left[\left( I_{\mathcal{K}}+Q(z_{1})\right)^{-1}-\left( I_{\mathcal{K}}+Q(z_{2})\right)^{-1}\right]AR_{0}(z_{1})\\
&+\overline{R_{0}(z_{2})B}\left( I_{\mathcal{K}}+Q(z_{2})\right)^{-1}\left[ Q(z_{1})-Q(z_{2})\right]\left( I_{\mathcal{K}}+Q(z_{1})\right)^{-1}AR_{0}(z_{1})\\
=&R(z_{1})-R(z_{2}).
\end{align*}
By multiplying both sides of \eqref{Second Resolvent} by $H(z_{2})-z_{2}I_{\mathcal{H}}$ on the left, we get
\[ \left(H(z_{2})-z_{2}I_{\mathcal{H}}\right)R(z_{1})-I_{H}=(z_{1}-z_{2})R(z_{1}).\]
This is equivalent to
\[ \left(H(z_{2})-z_{1}I_{\mathcal{H}}\right)R(z_{1})=I_{\mathcal{H}}.\]
Multiplying this identity by $H(z_{1})-z_{1}I_{\mathcal{H}}$ on the right, ensures that $H(z_{1})=H(z_{2})$ as claimed above.
We can therefore write $H\coloneqq H(z)$.

\underline{Step 4}. We complete the proof by showing that $H$ is an extension of the operator $H_{0}+BA$.
 
Let $z\in \mathsf{R}$ and let $f\in \mathrm{Dom}(H_{0}+BA)=\mathrm{Dom}(H_{0})\cap \mathrm{Dom}(BA)$. We set $g=(H_{0}-zI_{\mathcal{H}})f$. It follows from \eqref{R} and \eqref{B*R*} that
\[R(z)=R_{0}(z)-\overline{R(z)B}AR_{0}(z), \]
and thus that $R(z)g=f- R(z)BAf$. This completes the proof of the theorem.
\end{proof}

The next comments about the validity of the hypotheses \ref{Bounded Exten 1} and \ref{Bounded Exten 2}, when moving $z$, as now in place. The next observations was used in the step~3 in the proof of Theorem~\ref{Theo Closed Exten} and will be useful in the proof of Theorem \ref{Theo BS}. 
Here, $z_{1},z_{2}\in \rho(H_{0})$.

\subsubsection*{If $AR_{0}(z_{1})\in L(\mathcal{H},\mathcal{K})$, then $AR_{0}(z_{2})\in L(\mathcal{H},\mathcal{K})$} Indeed, we have 
\begin{equation*}
\begin{aligned}
AR_{0}(z_{2})&=A \left[ R_{0}(z_{1})+(z_{2}-z_{1})R_{0}(z_{1})R_{0}(z_{2})\right]\\
&\supseteq\underbrace{AR_{0}(z_{1})}_{\in L(\mathcal{H},\mathcal{K})}+(z_{2}-z_{1})\underbrace{AR_{0}(z_{1})}_{\in L(\mathcal{H},\mathcal{K})}\underbrace{R_{0}(z_{2})}_{\in L(\mathcal{H})}.
\end{aligned}
\end{equation*}
 In particular, since the domain of the right-hand side is total $\mathcal{H}$, we obtain 
\[ AR_{0}(z_{2})=AR_{0}(z_{1})+(z_{2}-z_{1})AR_{0}(z_{1})R_{0}(z_{2}).\]
Thus, $AR_{0}(z_{2}) \in L(\mathcal{H},\mathcal{K})$.

\subsubsection*{If $R_{0}(z_{1})B$ is closable and $\overline{R_{0}(z_{1})B}\in L(\mathcal{K},\mathcal{H})$, then $R_{0}(z_{2})B$ is also closable and $\overline{R_{0}(z_{2})B}\in L(\mathcal{K},\mathcal{H})$} Indeed, since $\mathrm{Dom}(R_{0}(z_{2})B)=\mathrm{Dom}(B)$ which is densely define, we have
\begin{align*}
\left( R_{0}(z_{2})B\right)^{*} &= \left[ R_{0}(z_{1})B+(z_{2}-z_{1})R_{0}(z_{2})R_{0}(z_{1})B\right]^{*}\\
&\supseteq \left[ R_{0}(z_{1})B \right]^{*}+(\overline{z_{2}}-\overline{z_{1}})\left[R_{0}(z_{2})R_{0}(z_{1})B\right]^{*}\\
&= \underbrace{\left(\overline{R_{0}(z_{1})B}\right)^{*}}_{\in L(\mathcal{H},\mathcal{K})}+(\overline{z_{2}}-\overline{z_{1}})\underbrace{\left(\overline{R_{0}(z_{1})B}\right)^{*}}_{\in L(\mathcal{H},\mathcal{K})}\underbrace{R_{0}(z_{2})^{*}}_{\in L(\mathcal{H})},
\end{align*}
where we used \cite[Prop. 1.6 (vi)]{Schmudgen12} for the inclusion and \cite[Prop. 1.7(ii) \& Theo. 1.8 (ii)]{Schmudgen12} for the last equality. Since the domain of the right-hand side is total $\mathcal{H}$, we deduce that 
\begin{equation*}
\left( R_{0}(z_{2})B\right)^{*} =\left(\overline{R_{0}(z_{1})B}\right)^{*}+(\overline{z_{2}}-\overline{z_{1}})\left(\overline{R_{0}(z_{1})B}\right)^{*}R_{0}(z_{2})^{*}\in L(\mathcal{H},\mathcal{K}).
\end{equation*}
Therefore, from \cite[Thm. 1.8 (i) \& (ii)]{Schmudgen12}, $R_{0}(z_{2})B$ is closable and $\overline{R_{0}(z_{2})B}=\left(R_{0}(z_{2})B\right)^{**}$ belongs to $ L(\mathcal{K},\mathcal{H})$, as claimed. Moreover, the following identity holds true,
\begin{equation}\label{RB Resolvent}
\overline{R_{0}(z_{1})B}-\overline{R_{0}(z_{2})B}=(z_{1}-z_{2})R_{0}(z_{2})\overline{R_{0}(z_{1})B}.
\end{equation}

\subsubsection*{If $\overline{AR_{0}(z_{1})B} \in L(\mathcal{K})$, then $\overline{AR_{0}(z_{2})B} \in L(\mathcal{K})$} This can be proved in the same manner as above:
\begin{align*}
\left( AR_{0}(z_{2})B\right)^{*} &= \left[ AR_{0}(z_{1})B+(z_{2}-z_{1})AR_{0}(z_{2})R_{0}(z_{1})B\right]^{*}\\
&\supset \left[ A R_{0}(z_{1})B \right]^{*}+(\overline{z_{2}}-\overline{z_{1}})\left[AR_{0}(z_{2})R_{0}(z_{1})B\right]^{*}\\
&= \underbrace{\left(\overline{AR_{0}(z_{1})B}\right)^{*}}_{\in L(\mathcal{K})}+(\overline{z_{2}}-\overline{z_{1}})\underbrace{\left(\overline{R_{0}(z_{1})B}\right)^{*}}_{\in L(\mathcal{H},\mathcal{K})}\underbrace{\left(AR_{0}(z_{2}\right)^{*}}_{\in L(\mathcal{K},\mathcal{H})}.
\end{align*}
Once again, since the domain of the right-hand side is total $\mathcal{K}$, \[ \left( AR_{0}(z_{2})B\right)^{*}=\left(\overline{AR_{0}(z_{1})B}\right)^{*}+(\overline{z_{2}}-\overline{z_{1}})\left(\overline{R_{0}(z_{1})B}\right)^{*}\left(AR_{0}(z_{2}\right)^{*} \in L(\mathcal{K}),\]
and the conclusion follows. Moreover, we get the formula
\begin{equation}\label{Q1-Q2}
Q(z_{1})-Q(z_{2})=(z_{1}-z_{2})AR_{0}(z_{2})\overline{R_{0}(z_{1})B}.
\end{equation}

\begin{theorem}\label{Theo BS}
Assume that the closed operators $H_0,\,A$ and $B$ are related in the above manner, and that \ref{Bounded Exten 1}--\ref{Nonempty} hold true. Additionally, assume that $\mathsf{R}$ has at least two elements. Then, for all $z\in \rho(H_{0})$,
\[ z\in \operatorname{Spec}_{\mathrm{p}}(H) \quad \text{if and only if} \quad-1 \in \operatorname{Spec}_{\mathrm{p}}(Q(z)).\] 
\end{theorem}
\begin{proof}
We first show that the left-hand side implies the right-hand side. 
Let $f\in \mathrm{Dom}(H)\setminus \{0\}$ be such that $Hf=zf$. Let $z_{0}\in \mathsf{R}\setminus \{z\}$. The equation $Hf=zf$ is equivalent to $f= (z-z_{0})R(z_{0})f$. By virtue of \eqref{R}, the latter happens if and only if
\[(H_{0}-zI_{\mathcal{H}})R_{0}(z_{0})f = -(z-z_{0})\overline{R_{0}(z_{0})B}\left( I_{\mathcal{K}}+Q(z_{0})\right)^{-1}AR_{0}(z_{0})f.\]
By setting $v=\left( I_{\mathcal{K}}+Q(z_{0})\right)^{-1}AR_{0}(z_{0})f\neq 0$ and applying $\left( I_{\mathcal{K}}+Q(z_{0})\right)^{-1}AR_{0}(z)$ on both sides, we obtain
\begin{equation*}
v=-(z-z_{0})\left( I_{\mathcal{K}}+Q(z_{0})\right)^{-1}AR_{0}(z)\overline{R_{0}(z_{0})B}v.
\end{equation*}
Then, according to \eqref{Q1-Q2}, 
\[ v= \left( I_{\mathcal{K}}+Q(z_{0})\right)^{-1}(Q(z_{0})-Q(z))v \qquad\iff \qquad Q(z)v=-v.\]
In other words, $-1$ is the eigenvalue of $Q(z)$ with the eigenfunction $v$, as the right-hand side of the conclusion states.

We now show the converse. Assume that $Q(z)v=-v$ for some $v\in \mathcal{K} \setminus \{0\}$. Take $z_{0}\in \mathsf{R}\setminus \{z\}$ and write
\begin{align*}
 v=&v-\left(I_{\mathcal{K}}+Q(z_{0})\right)^{-1}\left(I_{\mathcal{K}}+Q(z)\right)v=\left(I_{\mathcal{K}}+Q(z_{0})\right)^{-1}\left(Q(z_{0})-Q(z)\right)v\\
 =&(z_{0}-z)\left(I_{\mathcal{K}}+Q(z_{0})\right)^{-1}AR_{0}(z_{0})\overline{R_{0}(z)B}v.
\end{align*}
Here, in the last step, we applied \eqref{Q1-Q2}. By setting $f=-\overline{R_{0}(z)B}v \neq 0$, we get
\[ v= (z-z_{0})\left(I_{\mathcal{K}}+Q(z_{0})\right)^{-1}AR_{0}(z_{0})f.\]
Then, applying $\overline{R_{0}(z_{0})B}$ on both sides, and using \eqref{RB Resolvent} and \eqref{R}, we gather that
\[ -f+(z-z_{0})R_{0}(z_{0})f=(z-z_{0})\left(R_{0}(z_{0})-R(z_{0})\right)f.\]
This is equivalent to $ (z-z_{0})(H-z_{0}I_{\mathcal{H}})^{-1}f=f$ and thus, $f\in \mathrm{Dom}(H)$ and $Hf=zf$.
\end{proof}


\section*{Acknowledgement}
D.K. and T.N.D. were supported by the EXPRO grant number 20-17749X of the Czech Science Foundation (GA\v{C}R).


\bibliographystyle{abbrv}
\bibliography{Ref}
\end{document}